\documentclass[a4paper,11pt]{article}

\usepackage[utf8]{inputenc}
\usepackage{hyperref}
\usepackage{graphicx}
\usepackage{amsmath}
\usepackage{amsfonts}
\usepackage{titlesec}
\usepackage{amsthm}
\usepackage{enumitem}
\usepackage{amssymb}
\usepackage{epic}
\usepackage{tikz}

\usepackage{geometry}
\def\noopsort#1{}

\usepackage[
linesnumbered
,inoutnumbered
,commentsnumbered
,ruled
,lined
,vlined
]{algorithm2e}
\SetStartEndCondition{ }{}{}%
\SetKwProg{Fn}{def}{\string:}{}
\SetKwFunction{Range}{range}%%
\SetKw{KwTo}{in}\SetKwFor{For}{for}{\string:}{}%
\SetKwIF{If}{ElseIf}{Else}{if}{:}{elif}{else:}{}%
\SetKwFor{While}{while}{:}{fintq}%
\SetKwFor{Foach}{for}{:}{}%
\SetKwFor{ForEach}{for}{:}{}%
\SetKwInput{KwData}{data}%
\SetKwHangingKw{KwHData}{data\normalfont{:}}%
\SetKw{AddData}{\phantom{data\normalfont{:}}}%
\SetKwComment{tcc}{\#\# }{}%
\SetKwComment{tcp}{\# }{}%
\SetKw{Solve}{solve}%
\SetNlSty{}{\color{gray}\tiny}{}
\AlgoDontDisplayBlockMarkers\SetAlgoNoEnd\SetAlgoNoLine%

\newcommand{\quotes}[1]{``#1''}

\newcommand{\R}{\mathbb{R}}
\newcommand{\T}{\mathcal{T}}

\newcommand{\LOD}{\mbox{\tiny{LOD}}}
\newcommand{\PG}{\mbox{\tiny\rm PG}}
\newcommand{\Id}{\mbox{\rm Id}}

\newcommand{\coarse}{\mbox{\tiny{coarse}}}
\newcommand{\fine}{\mbox{\tiny{fine}}}
\newcommand{\corr}{\mbox{\tiny{corr}}}

\newcommand{\gmin}{\gamma_{\operatorname{min}}}
\newcommand{\gmax}{\gamma_{\operatorname{max}}}

\def\mathrlap{\mathpalette\mathrlapinternal}

\def\mathrlapinternal#1#2{\rlap{$\mathsurround=0pt#1{#2}$}}

\newcommand{\Vhell}{V_{h,\ell}}
\newcommand{\ttell}{{\mathrlap{\ell}\phantom{l}}}

\newcommand{\Qh}{\mathcal{Q}_h}

\newcommand{\QFh}{\mathcal{Q}_{\mathcal{F}^{\text{s}},h}}

\newcommand{\Ic}{I_H}

\newcommand{\supp}{\operatorname*{supp}}

\newcommand{\rows}{\operatorname*{rows}}
\newcommand{\VhGD}{V_{h,\Gamma_D}}
\newcommand{\VHGD}{V_{H,\Gamma_D}}
\newcommand{\bigsigma}{{\mbox{\LARGE$\boldsymbol{\sigma}$}}}

\newtheorem{theorem}{Theorem}[section]

\newtheorem{example}[theorem]{Example}
\theoremstyle{definition}
\newtheorem{definition}[theorem]{Definition}
\newtheorem{remark}[theorem]{Remark}

\begin{document}

\begin{center}
{\LARGE Efficient implementation of the Localized Orthogonal
  Decomposition method}\\[2em]
\end{center}

\renewcommand{\thefootnote}{\fnsymbol{footnote}}
\renewcommand{\thefootnote}{\arabic{footnote}}

\begin{center}
{\large Christian Engwer\footnote[1]{Institute for Computational and Applied Mathematics, University of M\"unster, Germany}, Patrick Henning\footnote[2]{Department of Mathematics, KTH Royal Institute of Technology, Stockholm, Sweden},
Axel M\r{a}lqvist\footnote[3]{Department of Mathematics, Chalmers University of Technology and University of Gothenburg, Sweden.},
Daniel Peterseim\footnote[4]{Institute for Mathematics, University of Augsburg, Germany}}\\[2em]
\end{center}

\begin{center}
{\large{\today}}
\end{center}

\begin{abstract}
  In this paper we present algorithms for an efficient implementation
  of the Localized Orthogonal Decomposition method (LOD). The LOD is a
  multiscale method for the numerical simulation of partial
  differential equations with a continuum of inseparable scales. We
  show how the method can be implemented in a fairly standard Finite
  Element framework and discuss its realization for different types of
  problems, such as linear elliptic problems with rough coefficients
  and linear eigenvalue problems.
\end{abstract}

\section{Introduction}
By now, the Localized Orthogonal Decomposition (LOD) of a subspace $V \subset H^1(\Omega)$ into a coarse space and a detail space is a well established method for the numerical homogenization of partial differential equations. So far it has been successfully applied to linear elliptic multiscale problems in the context of continuous finite elements \cite{MaP14,HeP13,HeM14}, discontinuous finite elements \cite{EGM13,EGM13b,Elf14},
mixed finite elements \cite{Mal11,HHM15},
partition of unity methods \cite{HMP13b} and reduced basis simulations \cite{AbH15}. The range of applications covers linear and quadratic eigenvalue problems \cite{MaP15,2015arXiv151005792M}, problems in perforated domains \cite{BrP14} and high-contrast media \cite{2016arXiv160106549P,HeM17}, stochastic homogenization \cite{2017arXiv170208858G,2018arXiv180701741F},  semilinear elliptic problems \cite{HMP14}, the wave equation \cite{AbH14c,Peterseim.Schedensack:2016,Maier2018},
parabolic and coupled problems \cite{MaPr15,refId0,2018arXiv180100615A},
the Buckley-Leverett equation \cite{EGH15}, fractional diffusion problems \cite{doi:10.1137/17M1147305}, 
Helmholtz problems \cite{Pet14,GaP15,Brown.Gallistl.Peterseim:2015,2016arXiv160804243B}
and the simulation of Bose-Einstein condensates \cite{HMP14b}. An introductionary general overview is given in \cite{Peterseim:2015}.

Initially inspired by the Variational Multiscale Method \cite{Hug95,HFM98,HuS07,LaM07,LaM09,Mal11}, the LOD in its present form was first proposed and rigorously justified in \cite{MaP14}. Further basic modifications of the method were suggested in \cite{HeP13}. The LOD is constructed to handle discrete problems that involve a high-dimensional solution space (also referred to as the 'fine space'). This typically takes place in two steps. In the first step the full fine space is decomposed into a low-dimensional space with good approximation properties and a high-dimensional remainder space. In the second step, this decomposition is localized in the sense that the low-dimensional space is approximated by constructing suitable locally supported basis functions that are the solutions of small patch problems. Due to their size, the patch problems are cheap to solve. Furthermore, they can be solved independently from each other and are hence perfect for parallelization. This strategy is particularly useful to reduce/distribute the computational cost of solving large systems of equations (arising from finite element discretizations). The method can be linked to conceptually very different techniques of mathematical modeling and scientific computing, e.g., it recovers the mathematical theory of homogenization \cite{2016arXiv160802092G} in periodic diffusion problems and even bridges to the theory of iterative solvers and subspace decomposition methods \cite{2016arXiv160804081K,2018arXiv181106319P}. Moreover, the method may be interpreted as a stabilization technique that coincides with the streamline upwind Petrov-Galerkin method SUPG \cite{HuS07}. In the last 5 years it has inspired numerous new developments in the field of multiscale partial differential equations including rough polyharmonic splines \cite{Owhadi2014}, iterative numerical homogenization \cite{KornhuberYserentant2015}, and gamblets \cite{Owhadi2017}.
While previous works focused on the
numerical analysis of the method, this paper aims at the detailed explanation of how the method can be algorithmically realized. We give detailed explanations on how the method works on an algebraic level. The results may as well be useful for implementing related multiscale methods.
% and we present numerical experiments that demonstrate
\section{Preliminaries}
\label{section-lod-initial-defs}
In this section we recall the Localized Orthogonal Decomposition (LOD) for finite element spaces. The decomposition is always with respect to a linear elliptic part of the differential operator.
\subsection{Computational domain and boundary}
For the rest of the paper, we consider a bounded polygonal domain $\Omega\subset\R^{d}$.
 %Lipschitz domain $\Omega\subset\R^{d}$ with a piecewise polygonal boundary $\partial \Omega$.
% We assume that $\partial \Omega$ is divided into two pairwise disjoint Hausdorff measurable submanifolds $\Gamma_{D}$ and $\Gamma_N$ with $\Gamma_D \cup \Gamma_N = \partial \Omega$ and $\Gamma_D$ being a closed set of non-zero Hausdorff measure of dimension $d-1$.
The boundary $\partial \Omega$ is divided into two parts $\Gamma_{D}$ and $\Gamma_N$. On $\Gamma_D$ we prescribe a Dirichlet boundary condition and on $\Gamma_N$ we prescribe a Neumann boundary condition. We have $\Gamma_D \cup \Gamma_N = \partial \Omega$ and we assume $\Gamma_D\not= \emptyset$. With that, we define the space $H^1_{\Gamma_D}(\Omega):=\{ v \in H^1(\Omega)| \hspace{2pt} v_{\vert \Gamma_D}=0\}$, where $v_{\vert \Gamma_D}=0$ is understood in the sense of traces.
%Later, we prescribe Dirichlet boundary conditions on $\Gamma_{D}$ and Neumann boundary conditions on $\Gamma_N$.
%By $n$ we denote the outward-pointing normal on $\partial \Omega$.
%Let $T_D : H^1(\Omega) \rightarrow H^{\frac{1}{2}}(\Gamma_D)$ denote trace operator with respect to $\Gamma_D$. We define the space $H^1_{\Gamma_D}(\Omega):=\{ v \in H^1(\Omega)| \hspace{2pt} T_D(v)=0\}$.
\subsection{Elliptic differential operator}
Subsequently we consider the following differential operator. Let $\kappa\in L^\infty(\Omega,\R^{d\times d})$ denote a matrix-valued, symmetric,
possibly highly varying and heterogeneous coefficient with uniform
spectral bounds $\gmin>0$ and $\gmax\geq\gmin$,
\begin{equation*}%\label{e:spectralbound}
  \sigma(\kappa(x))\subset [\gmin,\gmax]\quad\text{for almost all }x\in \Omega.
\end{equation*}
This coefficient defines a scalar product $\mathcal{A}(\cdot,\cdot)$ on $H^1_{\Gamma_D}(\Omega)$ that is given by
\begin{align*}
 \mathcal{A}\left(  v,w\right):=\int_{\Omega}\kappa \nabla v\cdot \nabla w \qquad \mbox{for } v,w \in H^1_{\Gamma_D}(\Omega).
\end{align*}
\subsection{Meshes and spaces}
  We wish to discretize a problem that is associated with $\mathcal{A}(\cdot,\cdot)$. Then the discretization is constrained by the diffusion coefficient $\kappa$, in the sense that variations of $\kappa$ must be resolved by the computational mesh. We call such a discretization a {\it fine scale discretization}. In addition to this, we have a second discretization on a {\it coarse scale}. The coarse mesh is arbitrary and no more related to $\mathcal{A}(\cdot,\cdot)$. It contains elements of maximum diameter $H>0$. The fine mesh consists
of elements of maximum diameter $h<H$.
Let $\T_H$, $\T_h$ denote the corresponding simplicial or
quadrilateral subdivisions of $\Omega$ into (closed)
conforming shape regular simplicial elements
or conforming shape regular quadrilateral elements, i.e.,
$\bar{\Omega}=\underset{{K_h\in\T_h}}{\bigcup}K_h=\underset{{K\in\T_H}}{\bigcup}K$. We
assume that $\T_h$ is a regular, possibly non-uniform, mesh refinement
of $\T_H$. Furthermore we also assume that $\T_H$ and $\T_h$ are
shape-regular in the sense that there exists a positive constant $c_0$
such that
\begin{align*}
  \max\left\{ \max_{K_h \in \T_h} \frac{\mbox{diam}(K_h)^d}{|K_h|}, \max_{K \in \T_H} \frac{\mbox{diam}(K)^d}{|K|} \right\} \le c_0
\end{align*}
and regular in the sense that any two elements are either disjoint,
share exactly one face, share exactly one edge, or share exactly one
vertex.  For $\T=\T_H,\T_h$, let
\begin{subequations}
\begin{align}
\label{notation-triangles}  P_1(\T) &= \{v \in C^0(\Omega) \;\vert \;\forall K\in\T,v\vert_K \text{ is a polynomial of total degree}\leq 1\} \enspace \mbox{and}\\
\label{notation-quadrilaterals}  Q_1(\T) &= \{v \in C^0(\Omega) \;\vert \;\forall K\in\T,v\vert_K \text{ is a polynomial of partial degree}\leq 1\}
\end{align}
\end{subequations}
denote the typical $p1$ degree and bi-$p1$ degree Finite Element
Spaces for triangular and quadrilateral partitions
respectively. We set $V_h:=P_1(\T_h)$ if $\T_h$ is a triangulation and $V_h:=Q_1(\T_h)$ if
$\T_h$ is a quadrilateration. The 'coarse space' (i.e. low dimensional space) $V_H \subset V_h$ is defined analogously. Furthermore, we set $\VhGD:=V_h \cap H^1_{\Gamma_D}(\Omega)$ and $\VHGD:=V_H \cap H^1_{\Gamma_D}(\Omega)$.  For simplicity we assume that $\VHGD$ is aligned with $\Gamma_D$ (in $2d$ this means that $\overline{\Gamma_D}\cap\overline{\Gamma_N}$ is a subset of coarse grid nodes).
The full sets of fine nodes in $V_h$, respectively coarse nodes in $V_H$, are given by
$$\mathcal{N}_H=\{ Z_i | \hspace{2pt}
0\le i \le N_H - 1\}
\qquad \mbox{and} \qquad
\mathcal{N}_h=\{ z_j| \hspace{2pt} 0\le j \le
N_h - 1\},$$
where $N_H$ and $N_h$ denote the number of vertices in the fine and
the coarse mesh. Accordingly we introduce $N_{\T_H} = |\T_H|$ and $N_{\T_h} = |\T_h|$
as the number of cells in the mesh.
The coarse nodal basis function that is associated with a node $Z_i \in \mathcal{N}_H$ shall be denoted by $\Phi_i \in V_H$ and the fine nodal basis function associated with $z_j \in \mathcal{N}_h$ shall be denoted by $\phi_j \in V_h$.
\subsection{Two-scale decompositions}
In order to introduce an $\mathcal{A}$-orthogonal decomposition of the space $\VhGD$, we require a projection $\Ic : \VhGD \rightarrow \VHGD$ (i.e. $(\Ic \circ \Ic)=\Ic$) that maps a fine-scale function into the coarse fine element space $\VHGD$. The chosen projection will help us to characterize the \quotes{details} in $\VhGD$ and it is desirable that $I_H$ is $L^2$- and $H^1$-stable. Before we introduce a decomposition based on $I_H$, we state examples of possible choices for $I_H$.
\begin{remark} \label{possible-choices-I_H}
Examples for projections $I_H$ that fulfill the desired stability properties on quasi-uniform meshes.
\begin{itemize}
\item The operator $I_H : \VhGD \rightarrow \VHGD$ can be chosen as the global $L^2$-projection onto finite elements given by
$$
(I_H(v_h),\Phi_H)_{L^2(\Omega)} = (v_h,\Phi_H)_{L^2(\Omega)} \qquad \mbox{for all } \Phi_H \in \VHGD.
$$
\item The operator $I_H : \VhGD \rightarrow \VHGD$ can be also constructed from a local $L^2$-projection. Given a coarse-node $Z_i$ and corresponding nodal patch $\omega_i:=\supp(\Phi_i)$, we let $P_{H,\omega_i}$ denote the $L^2$-projection onto the standard $P_1$ finite element space on $\omega_i$. Exploiting this, we define $\Ic$ for $v_h \in \VhGD$ by $\Ic(v_h):= \sum_{i=0}^{N_H-1} \alpha_i(v_h) \Phi_i$ where $\alpha_i(v_h)=0$ if $Z_i \in \mathcal{N}_H \cap \Gamma_D$ and $\alpha_i(v_h)=P_{H,\omega_i}(v_h)(Z_i)$ in any other case.
\item A similar construction is obtained by projecting locally into the space of discontinuous finite elements. Given $Z_i \in \mathcal{N}_H$ with corresponding nodal patch $\omega_i:=\supp(\Phi_i)$, we let $\tilde{P}_{H,\omega_i}$ denote the $L^2$-projection onto the space of functions on $\omega_i$ that are affine on each coarse grid element (discontinuous $P_1$ finite elements on $\omega_i$). For $v_h \in \VhGD$, we can now define $\Ic(v_h):= \sum_{i=0}^{N_H-1} \alpha_i(v_h) \Phi_i$ where $\alpha_i(v_h)=|\omega_i|^{-1}\int_{\omega_i}\tilde{P}_{H,\omega_i}(v_h)$ for all active nodes and $\alpha_i(v_h)=0$ if $Z_i \in \mathcal{N}_H \cap \Gamma_D$.
\item An example for a projection $\Ic : \VhGD \rightarrow \VHGD$ that is {\it not suitable} because it lacks  the desired stability properties is the Lagrange (nodal) interpolation.
\end{itemize}
There are also many other choices for $\Ic$, e.g., the orthogonal projection onto $V_H$ with respect to the $H^1$ inner product and quasi-interpolation operators of Cl\'ement or Scott-Zhang type as they are well-established in the finite element community in the context of fast solvers and a posteriori error estimation \cite{Car99,MR0400739,MR1011446,MR1706735,2015arXiv150506931E}. For some problems, it can be advantageous to equip $\Ic$ with information about the problem, e.g., $\kappa$-weighted $L^2$ averaging for high-contrast problems \cite{2016arXiv160106549P}.
As we see next, in practice we only require the kernel of the projection $\Ic$ for an implementation of the method. This simplifies the computations significantly. For instance, for typical choices of $\Ic$, there exist sets of functionals that can be used to decide if a function is in the kernel of $\Ic$ or not (cf. \cite[Section 2.3]{HuS07}).
\end{remark}
Once we decided for a suitable projection operator $I_H$, we can define the detail space
$$
W_h := \{ v_h \in \VhGD | \hspace{2pt} \Ic(v_h) = 0 \}.
$$
This detail space contains fine-scale functions in $\VhGD$ that cannot be expressed in the coarse space $\VHGD$. In terms of the LOD we wish to correct classical nodal basis functions by an appropriate \quotes{detail function} from the space $W_h$. This can be achieved in a natural way by introducing the following elliptic decomposition of $\VhGD$. We refer to this decomposition as the $\mathcal{A}$-orthogonal splitting of $\VhGD$ (cf. \cite{MaP14} for more details).
\begin{definition}[$\mathcal{A}$-orthogonal decomposition]
We define the $\mathcal{A}(\cdot,\cdot)$-orthogonal complement of $W_h$ in $\VhGD$ by
$$
V_{\LOD} := \{ v_h \in \VhGD| \hspace{2pt} \mathcal{A}(v_h,w_h) = 0 \enspace \mbox{for all } w_h \in W_h\}.
$$
This is well-defined since $\mathcal{A}(\cdot,\cdot)$ is a scalar product on $W_h$. We obtain the (ideal) splitting
$$\VhGD = V_{\LOD} \oplus W_h,$$
where $\mbox{dim}(V_{\LOD})=\mbox{dim}(\VHGD)$.
\end{definition}
We wish to use an approximation of $V_{\LOD}$ as a discrete solution space for Galerkin approximations. Observe that $V_{\LOD}$ is low dimensional, but practically expensive to assemble. Therefore we introduce a localized decomposition.
\subsection{Localization to patches}
To localize the splitting $\VhGD = V_{\LOD} \oplus W_h$, we first need
to localize the space $W_h$ to patches $U(K)\subset
\Omega$. We therefore introduce coarse-layer patches:
\begin{definition}[Coarse-layer patch]
For any positive $k\in \mathbb{N}$ and a coarse element $K \in \T_H$, we define patches $U_k(K)$ that consist of $K$ itself and $k$-surrounding layers of coarse elements, i.e. we define $U_k(K)$ iteratively by
\begin{equation}\label{def-patch-U-k}
    \begin{aligned}
      U_0(K) & := K, \\
      U_k(K) & := \cup\{T\in \T_H\;\vert\; T\cap U_{k-1}(K)\neq\emptyset\}\quad k=1,2,\ldots .
    \end{aligned}
\end{equation}
\end{definition}
The restriction of $W_h$ to a patch $U(K)$ is
defined by $W_h(U(K)):=\{ v_h \in W_h| \hspace{2pt}
v_h=0 \enspace \mbox{in } \Omega \setminus
U(K) \}$.
The localized decomposition can be now characterized using local correction operators.
\begin{definition}[Correction Operators]%[Localized Orthogonal Decomposition]
  \label{def-loc-lod}
  For a given positive $k\in \mathbb{N}$ and for a coarse
  function $\Phi_H \in \VHGD$, the correction operator $\Qh: \VHGD
  \rightarrow \VhGD$ is given by
  \begin{align*}
    \Qh(\Phi_H):=\sum_{K\in \T_H} \Qh^{K}(\Phi_H),
  \end{align*}
  where $\Qh^{K}(\Phi_h) \in W_{h}(U_k(K))$ (for $K\in \T_H$) is the solution of
  \begin{align}
    \label{local-corrector-problem}
    \int_{U_k(K)} \kappa \nabla \Qh^{K}(\Phi_H)\cdot \nabla w_h = - \int_K \kappa
    \nabla \Phi_H \cdot \nabla w_h \qquad \mbox{for all } w_h \in
    W_{h}(U_k(K)).
  \end{align}
%We denote $\Rh(\Phi_H):=\Phi_H+\Qh(\Phi_H)$.
\end{definition}
%\textcolor{red}{Die definition heisst zerlegung, es gibt aber keine zerlegung, sondern nur Korrektoren!}]
We obtain the space $\{ \Phi_H+\Qh(\Phi_H)| \hspace{2pt} \Phi_H \in \VHGD\}$ as a localized approximation of $V_{\LOD}$. Here, localized is to be understood in the sense that there exists a local (nodal) basis $\{\Phi_i+\Qh(\Phi_i)\,\vert\,Z_i\in\mathcal{N}_H\setminus\Gamma_D\}$ of $V_{\LOD}$ where the support of a basis function is restricted to $k+1$ layers of coarse elements around the corresponding node. The dimension of the new space is low (it is of the same dimension as $\VHGD$) and
and it can be constructed by solving the small problems \eqref{local-corrector-problem}, potentially in parallel. This generalized finite element space may be used in a Galerkin approximation of a prototypical linear elliptic model problem.
\begin{example}\label{ex:lodsource}
Let $f\in L^{2}( \Omega)$ and let $u\in H^1_{\Gamma_D}(\Omega)$ solve
\begin{equation*}
  \mathcal{A}\left(  u,v\right) =\int_{\Omega}fv \quad\text{for all } v\in H^1_{\Gamma_D}(\Omega).
\end{equation*}
Then the corresponding LOD approximation is given by $u_{\LOD}=u_H+\Qh(u_H)$, where $u_H \in \VHGD$ solves
  \begin{eqnarray*}
    \mathcal{A}\left( u_H+\Qh(u_H),\Phi_H+\Qh(\Phi_H) \right)&=&
    \int_{\Omega} f  (\Phi_H+\Qh(\Phi_H)) \quad \mbox{for all } \Phi_H \in \VHGD.
  \end{eqnarray*}
\end{example}

\begin{remark}[Alternative iterative localization]
	There is an alternative characterization of localized correctors $\Qh$ via a preconditioned iterative solver that is based on a domain decomposition preconditioner as proposed in \cite{2016arXiv160804081K} (which in turn is based on \cite{KornhuberYserentant2015}).
To explain this alternative strategy we define nodal patches $\Omega_i$ as the union of all elements $T\in \mathcal{T}_H$ that share the vertex $Z_i\in\mathcal{N}_H$ and let
\begin{equation}    \label{Wi} 
W_{h,i}:=\{v-I_H v\,|\,v\in \VhGD:v\vert_{\Omega\setminus\Omega_i}\equiv 0\}.
\end{equation}
The functions in $W_{h,i}$ are supported in a small 
neighbourhood of the vertex $Z_i$ depending on the choice of $I_H$ (typically within two layers of coarse elements). The $W_{h,i}$ are closed subspaces of the kernel $W_h$ of $I_H$, see \cite{2016arXiv160804081K}.
Let $P_i$ be the $\mathcal{A}$-orthogonal projection from 
$\VhGD$ to $W_{h,i}$, defined via the 
equation
\begin{equation}    \label{Pi}
\mathcal{A}(P_iv,w)=\mathcal{A}(v,w),\quad \forall w\in W_{h,i}.
\end{equation}
With this, we introduce an operator $P$ as
\begin{equation}    \label{P}
P=P_0+P_1+\cdots+P_{N_H-1}.
\end{equation}
The operator $P$ is symmetric with respect to the bilinear form $\mathcal{A}(\cdot,\cdot)$ and in \cite{KornhuberYserentant2015,2016arXiv160804081K} it is shown that it is a quasi-optimal preconditioner for the ideal corrector $\Qh^\infty$, i.e. for the $\mathcal{A}$-orthogonal projection from $\VHGD$ onto $W_h$. 
Starting from $\Qh^0(\bullet)=0$, localized approximations $\Qh^j$ of $\Qh^\infty$ can be defined, for any $\Phi_H\in\VHGD$, via the iteration
\begin{equation}  \label{correclocal}
\Qh^{j}(\Phi_H)=\Qh^{j-1}(\Phi_H)+\vartheta P(\Phi_H-\Qh^{j-1}(\Phi_H)),\quad j=1,2,\ldots,k.
\end{equation}
Note the information is spread at most by a fixed number of layers (typically two) of coarse elements in each step, that is, the support of $\Qh^{k}(\Phi_H)$ is at most $\mathcal{O}(k)$ layers larger than the support of $\Phi_H$. If the scaling factor $\vartheta>0$ is chosen appropriately as discussed in \cite[p. 2771]{2016arXiv160804081K}, the iteration converges linearly, thereby producing a localized corrector $\Qh=\Qh^{k}$ that is close to the one given in Definition~\ref{def-loc-lod}, not equal in general though. Please note that the parameter $k$ has a slightly different meaning in the two variants. To avoid the educated guess of $\vartheta$ and to achieve a more accurate global approximation one may consider the enlarged coarse space 
$$\{ \Phi_H+\Qh^j(\Phi_H)| \hspace{2pt} \Phi_H \in \VHGD,\;j=1,2,\ldots,k\}$$ 
which underlies the error analysis of \cite{2016arXiv160804081K}. Here, the choice $\vartheta=1$ is appropriate. This space provides the sharpest global error bounds at the price that the dimension of the final multiscale space is $k$-times larger compared to the dimension of the LOD multiscale space $(\Id+\Qh)(\VHGD)$ used in Example~\ref{ex:lodsource}. 
%in dimension of the coarse problem when compared to the variant of Example~\ref{ex:lodsource}. 

In terms of computational cost, the iterative computation of the corrector(s) is comparable to the original variant. Since only local problems of type  \eqref{Pi} need to be solved the possible degree of parallelism is slightly larger. However, in this paper we trade a slight reduction of offline efficiency for a better online performance which is achieved by the original variant as it typically produces more localized functions for given fixed accuracy.

We shall emphasize that this preference for the online efficiency pays off only if sufficiently many online computations are to be performed. If a source problem as given in Example~\ref{ex:lodsource} is to be solved for a very small number of right-hand sides the offline cost becomes relevant. In this case the iterative corrector computation seems favorable, in particular the variant introduced in \cite{KornhuberYserentant2015}. Here, the correctors are computed on the fly during a preconditioned iterative solution of the full problem in the spirit of \eqref{correclocal} (the corrector is only applied to the current approximation of the solution rather than precomputing it for all basis functions in the coarse space). For a more detailed discussion of this variant we refer to \cite{KornhuberYserentant2015} and to  \cite{10.1007/978-3-319-52389-7_21} for a comparison to the original variant and some numerical results. 
\end{remark}

In the following sections we quantify the approximation properties of $u_{\LOD}$ depending on the choice of the localization parameter $k$.
\section{The algebraic realization of the correctors $\Qh$}
\label{subsection-assembling-solving-local-problems}
Before we start to give a first example, we need to discuss how the
local problems (\ref{local-corrector-problem}) can be assembled and
solved practically. In particular we show how to interpret the
corrector $\Qh$ on an algebraic level. We employ the notation of matrices as we believe it eases reading but
want to point out that in an optimized implementation, most or all
linear operators can be implemented matrix-free.

\subsection{Analytic preliminaries}
We start by introducing a general terminology that we use subsequently in the context of localization. Every patch $U_k(K_{\ell})$ with $K_{\ell} \in \T_H$ and $k \in \mathbb{N}$ is in the following directly associated with the index $\ell$ (where $0 \le \ell < N_{\T_H}$). In particular, we consider the localization parameter $k$ to be fixed and hence drop it when defining $\mathcal{U}_{\ell}:=U_k(K_{\ell})$.
%Recall the set
%s of coarse nodes $\mathcal{N}_H=\{ Z_i | \hspace{2pt}
%0\le i \le N_H - 1\}$ and
%of fine nodes $\mathcal{N}_h=\{ z_j| \hspace{2pt} 0\le j \le
%N_h - 1\}$.
For a given patch $\mathcal{U}_{\ell}$ we denote the sets of active
coarse and
fines nodes in $\mathcal{U}_{\ell}$
respectively
by
\begin{align*}
\mathcal{N}_{\ell,H}&=\{ Z_{\ell,i} \in \mathcal{N}_{H} | \hspace{2pt} Z_{\ell,i} \in \overline{\mathcal{U}_{\ell}} \setminus \Gamma_D \} \qquad \hspace{25pt}\mbox{and}\\
\mathcal{N}_{\ell,h}&=\{ z_{\ell,j} \in \mathcal{N}_{h} | \hspace{2pt} z_{\ell,j} \in \overline{\mathcal{U}_{\ell}} \setminus \overline{\left(\partial \mathcal{U}_{\ell} \setminus \Gamma_N \right)} \}.
\end{align*}
%\comment[inline]{PH}{Check if nodes on $\partial \mathcal{U}_{\ell}$ should included for $\mathcal{N}_{\ell,h}$ or not (from the computation point of view)}
%{\Huge \textcolor{red}{In der zweiten Zeile sollte man den Abschluss von $\left(\partial \mathcal{U}_{\ell} \setminus \Gamma_N \right)$ abziehen}}
Furthermore, we set
$N_{\ell,H}:=|\mathcal{N}_{\ell,H}|$ and
$N_{\ell,h}:=|\mathcal{N}_{\ell,h}|$ the number of nodes. The corresponding
coarse and
fine Lagrange basis functions in the patch $\mathcal{U}_{\ell}$ are denoted
respectively
by
$\Phi_{\ell,i}$ (i.e. $\Phi_{\ell,i}$ is coarse nodal basis function for node $Z_{\ell,i} \in \mathcal{N}_{\ell,H}$) and
$\phi_{\ell,j}$ (i.e. $\phi_{\ell,j}$ is fine nodal basis function for node $z_{\ell,j} \in \mathcal{N}_{\ell,h}$). With this notation, we define the corresponding local basis function sets.
\begin{definition}[Local basis sets]
%%%%
For each patch $\mathcal{U}_{\ell}$ we define
%$$\Vhell:=V_h \cap H^1_0(\mathcal{U}_{\ell})$$
$$\Vhell:=\{ v_h \in V_h| \hspace{2pt} v_h(z_{\ell,j})=0 \mbox{ for } z_{\ell,j} \in \mathcal{N}_{h} \setminus \mathcal{N}_{\ell,h}\}$$
and we let
  \begin{equation*}
    \{\phi_{\ell,j}| \hspace{2pt} 0 \le j  < N_{\ell,h} \}\subset \Vhell
  \end{equation*}
  be the set of fine Lagrange basis functions that belong to the
  active fine nodes in $\mathcal{U}_{\ell}$ (i.e. to $\mathcal{N}_{\ell,h}$).
  Accordingly we let
  \begin{equation*}
    \{\Phi_{\ell,i}| \hspace{2pt} 0 \le i  < N_{\ell,H} \} \subset V_H
  \end{equation*}
  denote the ordered set of all active coarse Lagrange basis functions in $\mathcal{U}_{\ell}$, i.e. the coarse basis
  functions associated with the nodes in $\mathcal{N}_{\ell,H}$.
\end{definition}
Since each $K_{\ell} \in \T_H$ contains $c_d$ coarse nodes, we can order the global indices of these nodes by $p_{0}(\ell)< p_{1}(\ell) < \cdots <  p_{c_d-1}(\ell)$.
This implies
$\Qh^{K_{\ell}}(\Phi_{j})=0$ for all $j \not\in \{ p_0(\ell), \cdots, p_{c_d-1}(\ell)\}$, hence we only need to compute $\Qh^{K_{\ell}}(\Phi_{p_i(\ell)})$ for $i \in \{ 0, \cdots, c_d-1\}$. For arbitrary $\Phi_H \in \VHGD$ we can hence write
\begin{align}
\label{relation-loc-cor-glob-cor}  \Qh(\Phi_H)=\sum_{K_{\ell} \in \T_H} \sum_{i=0}^{c_d-1} \Phi_H(Z_{p_i(\ell)}) \Qh^{K_{\ell}}(\Phi_{p_i(\ell)}).
\end{align}
\subsection{Algebraic preliminaries}
In this section, we introduce some algebraic preliminaries. As a general notation in this paper we denote for any matrix $\mathcal{M} \in \R^{m \times n}$ the transposed of the $i$'th row of $\mathcal{M}$ by $\mathcal{M}[i] \in \R^n$ (for $0 \le i < m$), i.e.
\begin{align*}
\mathcal{M}[i]
=
\left( \begin{matrix}
\mathcal{M}_{i,0} \\
\vdots\\
\mathcal{M}_{i,n-1}
\end{matrix} \right).
\end{align*}
The entry at position $(i,j)$ is denoted by $\mathcal{M}[i][j]$.
\begin{definition}[Local-to-global-mapping]
  \label{local-to-global-mapping}
Let $c_d$ denote the number of nodes in a grid element, i.e. $c_d = d+1$ if $\T_H$ and $\T_h$ consist of simplicial elements and $c_d=2^d$ if $\T_H$ and $\T_h$ consist of quadrilateral elements. Then we call
$$\sigma : \{ 0, \ldots, N_{\T_h}-1\} \times \{ 0, \ldots, c_d -1 \} \rightarrow \{ 0, \ldots, N_h -1 \}$$
the local-to-global-maping for the fine grid, if it maps the local index of a node $m$ in an element $T_t \in \T_h$ to its global index $j$. Here, $t \in \{ 0, \ldots, N_{\T_h}-1\}$ denotes the index of the element $T_t \in \T_h$. We write $\sigma(t,m)=j$. An example is given in Figure~\ref{figure-sigma}.
For $T_t \in \T_h$ the algebraic version of $\sigma(t,\cdot)$ is
given by the matrix $\bigsigma_{\!t} \in \R^{N_h \times c_d}$ where for $0\le m < c_d$ and $0 \le j < N_h$
\begin{align}
\label{loc-to-glob-map}\bigsigma_{\!t}[m][j]:=
\begin{cases}
1 \qquad  &\mbox{if } \sigma(t,m)=j,\\
0 \qquad  &\mbox{else.}
\end{cases}
\end{align}
\end{definition}
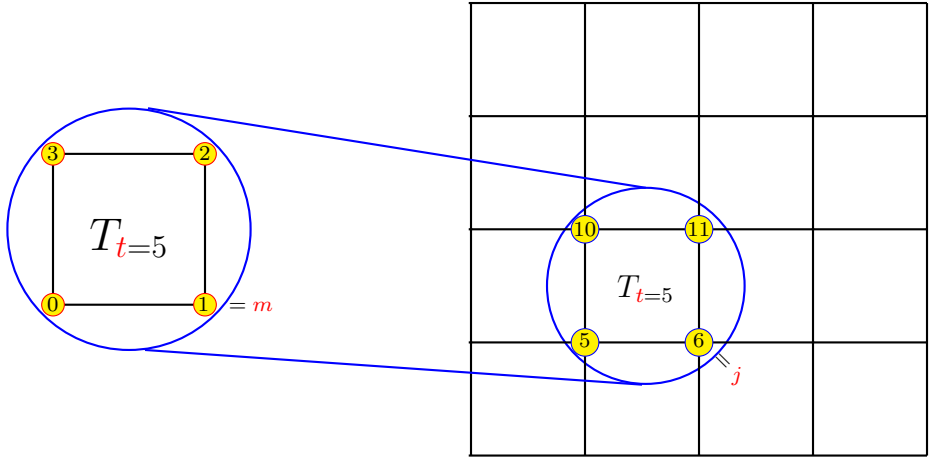
\begin{figure}
  \centering
\begin{tikzpicture}[scale=1]
\draw[thick,black] (2,2) -- (2,4) -- (4,4) -- (4,2) -- (2,2);
\draw[thick,blue] (3,3) circle (1.6);
\draw[thick,blue] (9.8,2.25) circle (1.3);
%grid
\draw[thick] (7.47,0) grid [step=1.5] (13.5,6);
\draw[thick,blue] (3.25,4.605) -- (9.8,3.55);
\draw[thick,blue] (3.21,1.4) -- (9.75,0.94);
\draw (3,2.5) node[above] {\LARGE$T_{{\color{red}t}=5}$};
\path[fill=yellow,draw=red] (2,2) circle (0.9ex);
\draw (2,1.78) node[above] {\scriptsize$0$};
\path[fill=yellow,draw=red] (4,2) circle (0.9ex);
\draw (4.4,1.78) node[above] {\scriptsize$1\hspace{5pt}={\color{red}m}$};
\path[fill=yellow,draw=red] (4,4) circle (0.9ex);
\draw (4,3.79) node[above] {\scriptsize$2$};
\path[fill=yellow,draw=red] (2,4) circle (0.9ex);
\draw (2,3.79) node[above] {\scriptsize$3$};
\path[fill=yellow,draw=blue] (9,1.5) circle (1.1ex);
\draw (9,1.30) node[above] {\scriptsize$5$};
\path[fill=yellow,draw=blue] (10.5,1.5) circle (1.1ex);
\draw (10.5,1.30) node[above] {\scriptsize{$6$}};
\draw (10.8,1.00) node[above] {\scriptsize{\rotatebox[origin=c]{-45}{$=$}}};
\draw (11.0,0.80) node[above] {\scriptsize{\color{red}$j$}};
\path[fill=yellow,draw=blue] (9,3) circle (1.1ex);
\draw (9,2.78) node[above] {\scriptsize$10$};
\path[fill=yellow,draw=blue] (10.5,3) circle (1.1ex);
\draw (10.5,2.78) node[above] {\scriptsize$11$};
\draw (9.8,1.9) node[above] {\large$T_{{\color{red}t}=5}$};
\end{tikzpicture}
\caption{Example: the fine grid $\T_h$ consists of $16$ quadratic grid elements. We start counting with $0$. For the 6'th element $T_5$ (i.e. $t=5$), we find 4 global nodes with global indices $5$, $6$, $10$ and $11$. Each of these global indices $j$, can be mapped to a local index $m$ and vice versa. This is done by $\sigma$.
For instance, for $t=5$, we have $\sigma(t,0)=5$, $\sigma(t,1)=6$, $\sigma(t,2)=11$ and $\sigma(t,3)=10$. For $m=1$ and $j=6$ as in the graphic, we have $\bigsigma_{\!t}[m][j]=1$, because $\sigma(t,m)=j$.}
\label{figure-sigma}
\end{figure}
\begin{definition}[Element stiffness and mass matrices]
\label{element-matrices}For $T_t \in \T_h$ with index $t$ we define the (fine grid) element stiffness matrix $\mathbf{A}_t \in \R^{c_d \times c_d}$ by
\begin{align*}
\mathbf{A}_t[m][n] := \int_{T_t} \kappa \nabla \phi_{\sigma(t,n)} \cdot \nabla \phi_{\sigma(t,m)} \qquad \mbox{for all } 0\le n,m <c_d
\end{align*}
and the (fine grid) element mass matrix
$\mathbf{M}_t \in \R^{c_d \times c_d}$ by
\begin{align*}
\mathbf{M}_t [m][n] := \int_{T_t}  \phi_{\sigma(t,n)} \hspace{2pt} \phi_{\sigma(t,m)} \qquad \mbox{for all } 0\le n,m <c_d.
\end{align*}
The corresponding global block matrices that store all these element contributions (i.e. stiffness and mass matrix in a discontinuous Galerkin discretization on the fine grid) shall be denoted by $\mathbf{A}_{dc}$ and $\mathbf{M}_{dc}$ respectively.
\end{definition}
\begin{definition}[Global stiffness and mass matrices]
  \label{def:global-matrices}
By $\mathbf{A}_h \in \R^{N_h \times N_h}$ we denote the global fine stiffness matrix with entries $\mathbf{A}_h[i][j]=(\kappa \nabla \phi_j,\nabla \phi_i)_{L^2(\Omega)}$ and by $\mathbf{M}_h \in \R^{N_h \times N_h}$ the global fine mass matrix with entries $\mathbf{M}_h[i][j]=( \phi_j,\phi_i)_{L^2(\Omega)}$. Analogously we denote by $\mathbf{A}_H \in \R^{N_H \times N_H}$ the coarse stiffness matrix and by $\mathbf{M}_H \in \R^{N_H \times N_H}$ the coarse mass matrix. Note that the above matrices are with respect to {\it all} nodes including the whole boundary.
\end{definition}
Assuming that the local contributions $\mathbf{A}_t,\mathbf{M}_t \in \R^{c_d \times c_d}$ are computed for all $t \in \T_h$ and assuming that $\bigsigma_{\!t}$ is available, we can assemble the global (fine) stiffness and mass matrix by
\begin{align}
\label{global-mass-and-stiffness-matrix}\mathbf{A}_h = \sum_{t \in \T_h} \bigsigma_{\!t} \mathbf{A}_t  \bigsigma_{\!t}^{\top} \qquad \mbox{and} \qquad
\mathbf{M}_h = \sum_{t \in \T_h} \bigsigma_{\!t} \mathbf{M}_t  \bigsigma_{\!t}^{\top}.
\end{align}
Since $V_h$ and $V_H$ do still incorporate all boundary nodes on $\partial \Omega$, we require matrices
that erases the unnecessary rows and columns in the system matrices (and load vectors) that are associated with the DOFs on the Dirichlet boundary part $\Gamma_D$.
\begin{definition}[Boundary correction matrices]
  \label{def:bc-corr}
We define the boundary correction (or restriction) matrices by
$\mathbf{B}^h \in \R^{N_h \times N_h}$ by
\begin{align*}
\mathbf{B}^h[i][j] =
\begin{cases}
1 \qquad &\mbox{if } i=j \mbox{ and } z_i \in \mathcal{N}_h \setminus \Gamma_D\\
0 \qquad &\mbox{else.}
\end{cases}
\end{align*}
and analogously
$\mathbf{B}^H \in \R^{N_H \times N_H}$ by
\begin{align*}
\mathbf{B}^H[i][j] =
\begin{cases}
1 \qquad &\mbox{if } i=j \mbox{ and } Z_i \in \mathcal{N}_H \setminus \Gamma_D\\
0 \qquad &\mbox{else.}
\end{cases}
\end{align*}
\end{definition}
As the coarse boundary condition vector is in general not available,
we describe an easy way to compute it. We define a vertex map that receives the index of a coarse node and the index of a fine node. If the coordinates of the coarse node are identical to the coordinates of the fine node, the vertex map is $1$ (true). In any other case, the vertex map is $0$ (false). In algebraic form, we describe the vertex map by the matrix $\mathbf{V}^h\in \R^{N_H \times N_h}$ that is given by
\begin{equation*}
  \mathbf{V}^h[i][j] :=
  \begin{cases}
    1 \qquad& \text{if }\mbox{\rm coord}(Z_i)=\mbox{\rm coord}(z_j),
    \quad\text{with } Z_i \in \mathcal{N}_H, z_j \in \mathcal{N}_h\\
    0 \qquad& \text{else.}
  \end{cases}
\end{equation*}
For instance, the matrix $\mathbf{V}^h$ can be easily computed by using an interpolation matrix $\mathbf{P}_h$
as defined in \eqref{projection-matrix} below. If $\mathbf{P}_h$ denotes such a matrix (expressing a function on the coarse grid in terms of fine degrees of freedom) we can define
  \begin{equation*}
    \mathbf{V}^h[i][j] := (\mathbf{P}_h[i][j] \equiv 1.0).
  \end{equation*}
%\comment{CE}{add more details on constraints handling\ldots{}}
\subsection{Local restriction matrices}
In order to localize computations to a patch $\mathcal{U}_{\ell}$ with $K_{\ell} \in \mathcal{T}_H$ we require a restriction operator $\mathcal{R}_{\ell} : V_h \rightarrow \Vhell$. For $v_h \in V_h$ we define the nodal interpolation $\mathcal{R}_{\ell}( v_h) \in \Vhell$ by
%{\Huge\textcolor{red}{notation overloading $\mathcal{R}$!}}
\begin{align*}
\mathcal{R}_{\ell}( v_h)( z ) =
\begin{cases}
 v_h( z ) \quad &\mbox{if } z \in \mathcal{N}_{\ell,h}, \\
 0 \quad &\mbox{else}.
\end{cases}
\end{align*}
The algebraic version of the restriction operator is denoted by $\mathbf{R}_{\ell}^h \in \R^{N_{\ell,h} \times N_{h}}$ and defined by the entries
\begin{align}
\label{restriction-op}\mathbf{R}_{\ell}^h [i][j] =
\begin{cases}
1 \qquad &\mbox{if } z_{\ell,i} = z_{j}\\%\in \mathcal{N}_{\ell,h}\\
0 \qquad &\mbox{else}
\end{cases}
%\hspace{50pt}\quad \mbox{for } 0 \le j < N_{\ell,h}
\end{align}
and where $i$ is the index that corresponds to the fine node $z_i \in \mathcal{N}_h$. Hence, we get the local mass matrix $\mathbf{M}_{\ell}$ (respectively local stiffness matrix $\mathbf{A}_{\ell}$) from the global mass matrix $\mathbf{M}_h$ (respectively stiffness matrix $\mathbf{A}_h$) by matrix multiplication, i.e.
\begin{align*}
 \mathbf{M}_{\ell} &= \mathbf{R}_{\ell}^h \mathbf{M}_h {\mathbf{R}_{\ell}^h}^{\!\top} \hspace{50pt} \mathbf{M}_{\ell} \in  \R^{N_{\ell,h} \times N_{\ell,h}} \mbox{ is local mass matrix,}\\
 \mathbf{A}_{\ell} &= \mathbf{R}_{\ell}^h \mathbf{A}_h {\mathbf{R}_{\ell}^h}^{\!\top} \hspace{52pt} \mathbf{A}_{\ell} \in  \R^{N_{\ell,h} \times N_{\ell,h}} \hspace{2pt} \mbox{ is local stiffness matrix.}
\end{align*}
Recall that the entries of $\mathbf{A}_{\ell}$ are given by
\begin{align*}
%  \mathbf{A}_{\ell} \in \R^{N_{\ell,h} \times N_{\ell,h}} \quad \mbox{has the entries } \quad
\mathbf{A}_{\ell}[i][j] = \mathcal{A}( \phi_{\ell,j}, \phi_{\ell,i} ) \quad \mbox{for} \enspace 0 \le i,j < N_{\ell,h}.
\end{align*}
Besides the restriction $\mathbf{R}_{\ell}^h$ to fine grid nodes in $\mathcal{U}_{\ell}$ we also require a restriction $\mathbf{R}_{\ell}^H$ to active coarse grid nodes in $\mathcal{U}_{\ell}$. We can define $\mathbf{R}_{\ell}^H \in \R^{N_{\ell,H} \times N_{H}}$
 analogously by the entries
\begin{align}
\label{restriction-op-coarse}\mathbf{R}_{\ell}^H [i][j] :=
\begin{cases}
1 \qquad &\mbox{if } Z_{\ell,i} = Z_j \\% \in \mathcal{N}_{\ell,H}\\
0 \qquad &\mbox{else}
\end{cases}
%\hspace{50pt}\quad \mbox{for } 0 \le j < N_{\ell,H}
\end{align}
and where $j$ is the global index that corresponds with the coarse node $Z_j \in \mathcal{N}_H$.
The local restriction matrices $\mathbf{R}_{\ell}^H$ and $\mathbf{R}_{\ell}^h$ need to be stored only temporary. Both matrices (as well as $\mathbf{A}_{\ell}$) can be deleted as soon as the $\ell$'th corrector matrix is computed (cf. Section~\ref{subsubsection-assembling-solving-loc-problems}).
\subsection{An algebraic characterization of the space $W_h$}
\label{subsubsection-algebraic-L2-projection}
Recall the notation
\begin{align*}
  V_H = \mbox{span}\{ \Phi_i| \hspace{2pt} 0\le i \le N_H - 1\} \quad \mbox{and} \quad V_h = \mbox{span}\{ \phi_j| \hspace{2pt} 0 \le j \le N_h - 1\}.
\end{align*}
This subsection describes how we can characterize the kernel of the projection $\Ic$. In order to illustrate our method, we restrict our considerations to the choice that $\Ic$ defines the $L^2$-projection. It is obvious that a function $v_h\in \VhGD$ is in the kernel of the global $L^2$-projection (i.e. $\Ic(v_h)=0$) if it holds
\begin{align}
\label{projection-contraint-new}(v_h,\Phi_i)_{L^2(\Omega)}=0 \quad &\mbox{for all } Z_i \in \mathcal{N}_H \setminus \Gamma_D.
\end{align}
In order to handle this constraint, we first observe that any coarse basis function $\Phi_i$ can be easily expressed in terms of fine basis functions by
$$%\begin{align}
\Phi_i = \sum_{j=0}^{N_h -1} \Phi_i(z_j) \phi_j.
$$%\end{align}
Hence we have
$$(\Phi_i , \Phi_j )_{L^2(\Omega)}=\sum_{k,\ell=0}^{N_h-1} \Phi_i(z_k) \hspace{2pt} (\phi_k , \phi_\ell )_{L^2(\Omega)} \hspace{2pt} \Phi_j(z_\ell).$$
Consequently, we can define the projection matrix $\mathbf{P}_h \in \R^{N_H \times N_h}$ from the coarse-mesh Lagrange space to the fine-mesh Lagrange space by
\begin{align}\label{projection-matrix}\mathbf{P}_h := \left(
      \begin{matrix}
        \Phi_0(z_0) & \cdots & \Phi_0(z_{N_h-1}) \\
        \vdots & \ddots & \vdots \\
        \Phi_{N_H-1}(z_0) & \cdots & \Phi_{N_H-1}(z_{N_h-1})
      \end{matrix}
    \right)
\end{align}
and relate the coarse and the fine mass matrix via
\begin{align*}
\mathbf{M}_H = \mathbf{P}_h \mathbf{M}_h \mathbf{P}_h^{\top}.
\end{align*}
With that we can see that the analytical constraint \eqref{projection-contraint-new}, i.e. $\Ic(v)=0$, can be equivalently expressed through the algebraic constraint
$(\mathbf{B}^H \mathbf{P}_h \mathbf{M}_h) \mathbf{v} = 0$. Motivated by these considerations we define the global constraint matrix
$\mathbf{C}_h \in \R^{N_H \times N_h}$ by
$$\mathbf{C}_h :=
\mathbf{P}_h \mathbf{M}_h.$$
\begin{remark}
The definition of $\mathbf{C}_h= \mathbf{P}_h \mathbf{M}_h$ might be surprising since it does not account for the distinction between Dirichlet-nodes and Neumann-nodes. In fact, the natural way is to define
$\mathbf{C}_h :=
\mathbf{B}^H \mathbf{P}_h \mathbf{M}_h$, where the restriction (or boundary correction) matrix $\mathbf{B}^H$ is used to remove the coarse basis functions associated with nodes on $\Gamma_D$. However, the boundary matrix $\mathbf{B}^H$ causes that $\mathbf{B}^H \mathbf{P}_h \mathbf{M}_h$ has not a maximal rank and hence the arising saddle point problem would be singular. This would cause numerical issues for the method with patches $\mathcal{U}_{\ell}=\Omega$. For that reason, we define $\mathbf{C}_h = \mathbf{P}_h \mathbf{M}_h$ and note that the arising constraint would be stronger than condition \eqref{projection-contraint-new}. However, the smaller error that we make in the definition of  $\mathbf{C}_h$ is already corrected in the next step by using local restrictions $\mathbf{R}_{\ell}^H$ and $\mathbf{R}_{\ell}^h$.
\end{remark}
Recalling the definition of the local restriction matrix $\mathbf{R}_{\ell}^h$ given by (\ref{restriction-op}), we can define the localization of
%the quasi-interpolation
$\mathbf{C}_h$ to the patch $\mathcal{U}_{\ell}$ (and to the correct boundary nodes) by
\begin{align*}
 \mathbf{C}_{\ell} &= \mathbf{R}_{\ell}^H \mathbf{C}_h {\mathbf{R}_{\ell}^h}^{\!\top} \hspace{50pt} \mbox{where } \mathbf{C}_{\ell} \in  \R^{N_{\ell,H} \times N_{\ell,h}} \mbox{ is the local constraint matrix.}
\end{align*}
Observe that if $\mathbf{C}_{\ell} \mathbf{v} = 0$ for some $\mathbf{v} \in \R^{N_{\ell,h}}$, then the represented function $v = \sum_{j=0}^{N_{\ell,h}-1} \mathbf{v}_j \phi_{\ell,j} \in \Vhell$ has the property
\begin{equation*}
0 = \int_{\mathcal{U}_{\ell}} v \Phi_{\ell,i} = \int_{\Omega} v \Phi_{\ell,i}  \qquad \mbox{for all } 0 \le i \le N_{\ell,H}-1,
\end{equation*}
i.e. $L^2$-orthongonality for all coarse basis functions $\Phi_{j}$ that have a support intersecting $\mathcal{U}_{\ell} = \mbox{supp}\hspace{2pt}v$. Consequently, the property holds for the whole space, $v \perp_{L^2} \VHGD = 0$, which means that $v$ is in the kernel of $L^2$-projection $I_H$ as desired. The converse conclusion follows analogously.
\begin{remark}\label{remark-on-size-of-C-ell}
The matrix $\mathbf{C}_{\ell}$ fully represents the local constraints and maps a function from the fine scale finite element space $\Vhell$ onto the coarse finite element space $V_H$ (restricted to the local subdomain $\mathcal{U}_{\ell}$). Since there are only constraints for coarse vertices, the dimension of the first component of $\mathbf{C}_{\ell}$ is small. If we use a localization parameter $k$ with $k \simeq C |\log(H)|$ (as it will be suggested by Theorem
\ref{theorem-apriori-linear} below), we obtain that $N_{\ell,h} = \mathcal{O}((H |\log(H)|) /h
)^{-d})$ (which is the dimension of $\Vhell$) and that $N_{H,\ell}$ (the number of coarse nodes in $\mathcal{U}_{\ell}$) grows proportional to $|\log(H)|^d$.
\end{remark}
\subsection{Assembling of a local load vector}
\label{subsubsection-assembling-load-vector}
Let us again fix some coarse element $K_{\ell} \in \T_H$. Beside assembling the local stiffness matrices and the local constraints matrix, we also need to compute the load vector that corresponds to the right hand side in (\ref{local-corrector-problem}), i.e. the term
$$- \int_{K_{\ell}} \kappa \nabla \Phi_i \cdot \nabla \phi_j$$
for every coarse and fine basis function $\Phi_i\in\VHGD$ and $\phi_j \in\VhGD$ with support on $K_{\ell}$.
We start with defining a suitable (algebraic) restriction operator for coarse grid nodes (in $K_{\ell}$).
\begin{definition}[Coarse-node-in-coarse-element restriction]
Recall that each coarse element $K_{\ell} \in \T_H$ contains $c_d$ coarse nodes. Let the global indices of these nodes be denoted by $p_{0}(\ell)< p_{1}(\ell) < \cdots <  p_{c_d-1}(\ell)$. Then the coarse element restriction matrix
$\mathbf{T}_{\ell}^H \in \R^{c_d \times N_H}$ is given for $0 \le i < c_d$ and $0 \le j < N_H$ by
\begin{align}
\label{restriction-op-coarse-element} \mathbf{T}_{\ell}^H [i][j] =
\begin{cases}
1 \qquad &\mbox{if } j=p_i(\ell) \\
0 \qquad &\mbox{else.}
\end{cases}
\end{align}
\end{definition}
As for the global stiffness matrix $\mathbf{A}_h$ in \eqref{global-mass-and-stiffness-matrix}, we can obtain the local stiffness matrix on $K_{\ell}$ (i.e. with entries $(\kappa \nabla \phi_i, \nabla \phi_j)_{L^2(K_{\ell})}$) from the element stiffness matrices $\mathbf{A}_t$. After that, we can restrict the resulting matrix to the fine basis functions that belong to fine nodes $z_i \in \mathcal{U}_{\ell}$ by using $\mathbf{R}_{\ell}^h$ (see \eqref{restriction-op}). We obtain for $0 \le j < N_{h}$ and $0 \le i < N_{\ell,h}$
\begin{align*}
\Big( ( \underset{t\subset K_{\ell}}{\sum_{t \in \T_h}} \bigsigma_{\!t} \mathbf{A}_t \bigsigma_{\!t}^{\top} ) {\mathbf{R}_{\ell}^h}^{\!\top} \Big) [j][i]
= (\kappa  \nabla \phi_{\ell,i}, \nabla \phi_j)_{L^2(K_{\ell})}.
\end{align*}
Consequently we can define the matrix ${\mathbf{r}}_\ell \in \R^{c_d \times N_{\ell,h}}$ that stores the $c_d$ load vectors as its rows
\begin{align}
\label{formula-local-rhs}{\mathbf{r}}_\ell
:= - \mathbf{T}_{\ell}^H \mathbf{B}^H \mathbf{P}_h ( \underset{t\subset K_{\ell}}{\sum_{t \in \T_h}} \bigsigma_{\!t} \mathbf{A}_t \bigsigma_{\!t}^{\top} ) {\mathbf{R}_{\ell}^h}^{\!\top},
\end{align}
i.e. ${\mathbf{r}}_\ell$ as defined above fulfills either ${\mathbf{r}}_\ell[i]=0$ for $Z_{p_i(\ell)} \in \Gamma_D$ or else its rows
transposed are given by
\begin{align*}
{\mathbf{r}}_\ell[i] = -
\left( \begin{matrix}
(\kappa \nabla \Phi_{p_i(\ell)}, \nabla \phi_{\ell,0} )_{L^2(K_{\ell})},\\
\vdots\\
  \hspace{14pt}( \kappa \nabla \Phi_{p_i(\ell)}, \nabla \phi_{\ell,M_{\ell,h}-1} )_{L^2(K_{\ell})} \hspace{8pt}
\end{matrix} \right),
\end{align*}
where $\Phi_{p_i(\ell)}$ is the $p_i(\ell)$'th coarse basis function (i.e. the global index is $p_i(\ell)$ and the local index in $K_{\ell}$ is $i$).
If ${\mathbf{r}}_\ell[i]=0$, no local problem has to be solved and the local corrector is zero.
\subsection{Assembly and solution of a local problem}
\label{subsubsection-assembling-solving-loc-problems}
Observe that (\ref{local-corrector-problem}) must be solved for every $K _{\ell}\in
\T_H$ and every coarse basis function $\Phi_m$ that has a support on
$K_{\ell}$ (i.e. for $\Phi_{p_i(\ell)}$ with $0\le i < c_d$, except the ones that belong to nodes on $\Gamma_D$) and
recall that the correct boundary condition on $\partial \mathcal{U}_{\ell}$ is already included in the local stiffness matrix $\mathbf{A}_{\ell}$.
Let us fix $\ell \in \{ 0 , \ldots, |\T_H| -1 \}$ and %consequently also $K_{\ell}$ and $\mathcal{U}_{\ell}$ with $k\in \mathbb{N}$.
%an associated basis function $\boldsymbol{\Phi}_{\ell}:=\Phi_i$ with support on $\mathcal{U}_{\ell}$.
 a coarse basis function $\Phi_{p_i(\ell)}$. In the light of the discussion in Section~\ref{subsubsection-algebraic-L2-projection}, we can formulate the local problem (\ref{local-corrector-problem}) in the following way.
\begin{definition}[Continuous formulation of a local problem]
Let $0\le i <c_d$ and let us denote $w_{\ell,i}:=\Qh^{K_{\ell}}(\Phi_{p_i(\ell)})$. Then $w_{\ell,i} \in V_h(\mathcal{U}_{\ell})$ is characterized by the property $\Ic(w_{\ell,i}) =0$ and the property that it solves
  \begin{align}
    \label{local-corrector-problem-new}
    \int_{\mathcal{U}_{\ell}} \kappa \nabla w_{\ell,i}\cdot \nabla w_h = - \int_{K_{\ell}} \kappa
    \nabla \Phi_{p_i(\ell)} \cdot \nabla w_h
\end{align}
for all $w_h \in V_h(\mathcal{U}_{\ell})$ with $\Ic(w_h) =0$.
\end{definition}
Problem (\ref{local-corrector-problem-new}) can be obviously interpreted as a saddle point problem. Hence, we obtain the following algebraic formulation using the notation from the previous subsections.
\begin{definition}[Algebraic formulation of a local problem]
\label{def-algebraic-formulation-og-local-problems}
Let $0\le i <c_d$. The algebraic version of problem (\ref{local-corrector-problem-new}) is the following saddle point problem. Find the tuple $(\mathbf{w}_{\ell}[i],\boldsymbol{\lambda}_{\ell}[i]) \in \R^{N_{\ell,h}} \times \R^{N_{\ell,H}}$; with
\begin{align}
\label{local-corrector-problem-new-algebraic}  \mathbf{A}_{\ell} \hspace{2pt} \mathbf{w}_{\ell}[i] + \mathbf{C}_\ell^{\top} \boldsymbol{\lambda}_{\ell}[i] &= {\mathbf{r}}_\ell[i]\\
\nonumber  \mathbf{C}_\ell \hspace{2pt} \mathbf{w}_{\ell}[i] &= 0.
\end{align}
Here, $\mathbf{w}_{\ell}[i]$ is the coefficient vector for the solution $w_{\ell,i} = \Qh^{K_{\ell}}(\Phi_{p_i(\ell)})$ of (\ref{local-corrector-problem-new}), i.e.
$$w_{\ell,i} = \sum_{j=0}^{N_{\ell,h}-1} \mathbf{w}_{\ell}[i][j] \hspace{2pt} \phi_{\ell,j},$$
and $\boldsymbol{\lambda}_{\ell}[i]$ is the corresponding Lagrange multiplier.
\end{definition}
We can state this problem also in Schur complement formulation.
\begin{remark}[Schur complement]
The Schur complement matrix $\mathbf{S}_{\ell}$ associated with problem \eqref{local-corrector-problem-new-algebraic} is given by
\begin{align}
\label{schurcomplement}\mathbf{S}_{\ell} := (\mathbf{C}_\ell \hspace{2pt} \mathbf{A}^{\hspace{-2pt}-1}_\ell \hspace{2pt} \mathbf{C}_\ell^{\top}).
\end{align}
Hence the solution $\mathbf{w}_{\ell}[i]$ of \eqref{local-corrector-problem-new-algebraic} can be written as
\begin{align}
\label{schurcomplement-2}
\mathbf{w}_{\ell}[i] = \mathbf{A}^{\hspace{-2pt}-1}_\ell {\mathbf{r}}_\ell[i] - (\mathbf{A}^{\hspace{-2pt}-1}_\ell \mathbf{C}_\ell^{\top}) \boldsymbol{\lambda}_{\ell}[i],
\end{align}
where
$\boldsymbol{\lambda}_{\ell}[i] \in \mathbb{R}^{N_{\ell,H}}$ solves
\begin{align}
\label{lagrangemultiplier}\mathbf{S}_{\ell} \boldsymbol{\lambda}_{\ell}[i] = (\mathbf{C}_\ell \hspace{2pt} \mathbf{A}^{\hspace{-2pt}-1}_\ell) \mathbf{r}_\ell[i].
\end{align}
\end{remark}
The common approach is to solve systems such as (\ref{schurcomplement})-(\ref{lagrangemultiplier}) iteratively with an approximate Schur complement matrix. As the system (\ref{lagrangemultiplier}) is only of size
$N_{\ell,H} \times N_{\ell,H}$ and must be solved $c_d$ times (for different righthand sides corresponding to each coarse basis function with support in $\mathcal{K}_\ell)$) it is faster to compute the whole Schur-complement matrix, solve it directly and apply back-substitution for each right-hand-side vector.
Solving the local problem \eqref{local-corrector-problem-new-algebraic} for all (transposed) rows of the matrix ${\mathbf{r}}_\ell \in \R^{c_d \times N_{\ell,h}}$ can be hence obtained in the following way. It involves a pre-processing step that is independent of ${\mathbf{r}}_\ell$ and a post-processing step that must be performed for each row of ${\mathbf{r}}_\ell$.
$\\$
{\bf Pre-processing steps.}
\begin{enumerate}
\item Compute the matrix $\mathbf{Y}_{\ell} := \mathbf{A}^{\hspace{-2pt}-1}_\ell \hspace{2pt} \mathbf{C}_\ell^{\top}$. This involves to solve $N_{\ell,H}$ problems of size $N_{\ell,h} \times N_{\ell,h}$, i.e. for all $0 \le m<N_{\ell,H}$ we need to solve for $\mathbf{Y}_{\ell}[m] \in \R^{N_{\ell,h}}$ with
$$\mathbf{A}_{\ell} \hspace{3pt} \mathbf{Y}_{\ell}[m] = \mathbf{C}_\ell^{\top}[m].$$
The $(1 \times N_{\ell,h})$-matrix $(\mathbf{Y}_{\ell}[m])^{\top}$ forms the $m$'th row of $\mathbf{Y}_{\ell}$.
\item Assemble the Schur complement $\mathbf{S}_{\ell} = \mathbf{C}_\ell \mathbf{Y}_{\ell}$ by matrix multiplication and compute $\mathbf{S}_{\ell}^{-1}$. Since $\mathbf{S}_{\ell}$ is only a $(N_{\ell,H} \times N_{\ell,H})$-matrix its inversion is cheap.
\end{enumerate}
$\\$
{\bf Post-processing steps for all $0 \le i < c_d$.}
\begin{enumerate}
\item If $\mathbf{r}_{\ell}[i]\neq0$, solve for $\mathbf{q}_{\ell}[i] \in \R^{N_{\ell,h}}$ with $\mathbf{A}_{\ell} \mathbf{q}_{\ell}[i] = \mathbf{r}_{\ell}[i]$.
\item Since
%$\mathbf{Y}_{\ell}= \mathbf{A}^{\hspace{-2pt}-1}_\ell \hspace{2pt} \mathbf{C}_\ell^{\top}$ and
$\mathbf{S}_{\ell}^{-1}$ is precomputed, we obtain $\boldsymbol{\lambda}_{\ell}[i]$ from $\mathbf{q}_{\ell}[i]$ via equation \eqref{lagrangemultiplier}, i.e. set
$$\boldsymbol{\lambda}_{\ell}[i] %= (\mathbf{S}_{\ell}^{-1} \mathbf{C}_\ell \hspace{2pt} \mathbf{A}^{\hspace{-2pt}-1}_\ell) \mathbf{r}_\ell[i]
= \mathbf{S}_{\ell}^{-1} \mathbf{C}_\ell \mathbf{q}_{\ell}[i].$$
\item Using the precomputed matrix $\mathbf{Y}_{\ell}= \mathbf{A}^{\hspace{-2pt}-1}_\ell \hspace{2pt} \mathbf{C}_\ell^{\top}$ and inserting $\boldsymbol{\lambda}_{\ell}[i]$ in \eqref{schurcomplement-2} we obtain
$$\mathbf{w}_{\ell}[i] = \mathbf{q}_{\ell}[i] - \mathbf{Y}_{\ell} \boldsymbol{\lambda}_{\ell}[i].$$
\end{enumerate}
When the post-processing step is concluded for the $\ell$'th local problem, a term of the form $(\mathbf{T}_{\ell}^H)^{\top} \mathbf{w}_{\ell} \mathbf{R}_{\ell}^h$ needs to be stored in a global corrector matrix $\mathbf{Q}_h$. Once this is done, all local matrices involved in the pre- and post-processing steps above are no longer required and can be deleted.
\begin{remark}[Cost]
Recall that $N_{\ell,H}\approx |\log(H)|^d$ is typically a small number (see also Remark \ref{remark-on-size-of-C-ell} above). The pre-processing step requires to solve $N_{\ell,H}$ equations of size $N_{\ell,h} \times N_{\ell,h}$ and to invert one matrix of size $N_{\ell,H} \times N_{\ell,H}$ (cost $\mathcal{O}(N_{\ell,H}^3)$). And in the post-processing step, for each $i\in \{0,\cdots,c_d-1\}$, we only need to solve one additional problem. In total, for one patch $\mathcal{U}_{\ell}$, the procedure involves to solve $(c_d + N_{\ell,H})$ equations of dimension $N_{\ell,h} \times N_{\ell,h}$ and $N_{\ell,H}$ equations of dimension $N_{\ell,H} \times N_{\ell,H}$. This also justifies why we solve the saddle point problem (\ref{local-corrector-problem-new-algebraic}) with a direct inversion of the Schur complement instead of using an iterative solver like the Uzawa solver. Roughly speaking, if the average number of iterations of an
iterative solver is larger than $(c_d+N_{\ell,H})/c_d$, then the direct inversion above is the cheaper
approach. This is in most cases fulfilled.
\end{remark}
\begin{remark}
Note that we can practically use the fact that the Lagrange basis
functions of $V_H$ have a partition of unity property, which implies
that it is only required to solve the local corrector problem
\eqref{local-corrector-problem} $d \cdot |\T_H|$ times in the case of
a triangulation and $(d+1) \cdot |\T_H|$ times in the case of a
quadrilation. We do not consider this in the algorithms. However, a corresponding modification ist straightforward.
\end{remark}
\begin{algorithm}
  \DontPrintSemicolon
  \KwData{$N_h$, $N_H$, $N_{\T_H}$, \tcp*{fine space / coarse space / coarse grid size}}
  \AddData{$\mathbf{A}_{\text{dc}}$, $\mathbf{M}_{\text{dc}}$,
    \tcp*{element stiffness / mass matrix, cf.\hspace*{.6ex}Def.\hspace*{.6ex}\ref{element-matrices}}}
  % \AddData{$\mathbf{b}^H$,
  %   \tcp*{coarse Dirichlet constraints vector}}
  \AddData{$\mathbf{B}^H$,
    \tcp*{$\tt N_H \times N_H$~~boundary correction matrix}}
  \AddData{$\mathbf{P}_h$
%      = \left(
%      \begin{matrix}
%        \Phi_0(z_0) & \cdots & \Phi_0(z_{N_h-1}) \\
%        \vdots & \ddots & \vdots \\
%        \Phi_{N_H-1}(z_0) & \cdots & \Phi_{N_H-1}(z_{N_h-1})
%      \end{matrix}
%    \right)$
    \tcp*{$\tt N_H \times N_h$~~projection matrix \hspace*{.6ex}in\hspace*{.6ex}~~(\ref{projection-matrix})}}
  \AddData{$\mathbf{R}_{\ell}^h$, $\mathbf{R}_{\ell}^H$, $\mathbf{T}_{\ell}^H$
    \tcp*{local restrictions in~~(\ref{restriction-op}),(\ref{restriction-op-coarse}),(\ref{restriction-op-coarse-element})}}
  \AddData{$\bigsigma$,
    \tcp*{per elem.\hspace*{1ex}to conforming map \hspace*{1ex}in~~\eqref{loc-to-glob-map}}}
  % \tcp*{DG to conforming mapping in \eqref{loc-to-glob-map}}}
  % \hfill\parbox{4.5cm}{\tcp*[f]{$\tt N_H \times N_h$ interp.}\par\tcp*[f]{~~~~~~ matrix~}}}
  \Fn{computeCorrections}{
    compute $\mathbf{A}_h = \bigsigma \mathbf{A}_{dc}  \bigsigma^{\top}$\tcp*{$\tt N_h \times N_h$ Stiffness matrix}
    compute $\mathbf{M}_h = \bigsigma \mathbf{M}_{dc}  \bigsigma^{\top}$\tcp*{$\tt N_h \times N_h$ Mass ~~~~~matrix}
%    compute $\mathbf{M}_H := \mathbf{P}_h \mathbf{M}_h \mathbf{P}_h^{\top}$\tcp*{$\tt N_H \times N_H$ Mass ~~~~~matrix}
    compute $\mathbf{C}_h := \mathbf{P}_h \mathbf{M}_h $\tcp*{$\tt N_H \times N_h$ constraints ~~matrix}
    $\mathbf{Q}_h := $ Matrix$(N_{H}, N_{h})$ \tcp*{$\tt N_{H} \times
      N_{h}$ Corrector matrix}
    \tcc{-- foreach $K_\ell \in \mathcal{T}_H$ --}
    \lnl{alg:main:subdomains}\ForEach{$0 \le \ell < N_{\T_H}$}
    {
      $N_{\ell,h} := \rows(\mathbf{R}_{\ell}^h)$
      \tcp*{local fine ~~space size}
      % \KwData{local fine space size $N_{\ell,h}$\tcp*{from \hspace{2pt}$\mathbf{R}_{\ell}^h$}}
      $N_{\ell,H} := \rows(\mathbf{R}_{\ell}^H)$
      \tcp*{local coarse space size}
      %\AddData{{\rm local coarse space size} $N_{\ell,H}$\tcp*{from $\mathbf{R}_{\ell}^H$}}
      %define $R_{\ell} : V_h \rightarrow $\;
      % \nlset{REM}
      $\mathbf{A}_{\ell} := \mathbf{R}_{\ell}^h \mathbf{A}_h {\mathbf{R}_{\ell}^h}^{\!\top}$\tcp*{$\tt N_{\ttell,h} \times N_{\ttell,h}$ Stiffness matrix}
      $\mathbf{C}_\ell := \mathbf{R}_{\ell}^H \mathbf{C}_h {\mathbf{R}_{\ell}^h}^{\!\top}$\tcp*{$\tt N_{\ttell,H} \times N_{\ttell,h}$ constraints ~~matrix}
      ${\mathbf{r}}_\ell := - \mathbf{T}_{\ell}^H \mathbf{B}^H
      \mathbf{P}_h (\bigsigma_{\!\ell} \mathbf{A}_t \bigsigma_{\!\ell}^{\top} {\mathbf{R}_{\ell}^h}^{\!\top})$
      \tcp*{$\tt \!~c_{d} \times N_{\ttell,h}$  load vec.\hspace{-5pt} matrix}
      \tcc{-- compute inverse operator --}
      $A_{\ell}^{\text{inv}} := \mathbf{A}_{\ell}^{-1}$\tcp*{e.g. using sparse ~~LU}
      \tcc{-- precomputations related to the operator --}
      $\mathbf{Y}_{\ell} := $ Matrix$(N_{\ell,H}, N_{\ell,h})$ \tcp*{$\tt N_{\ell,H} \times N_{\ell,h}$ matrix}
%      \ForEach{$\Phi_i \in \basis(V_H)$,
%        with $\supp(\Phi_i) \cap U(K_\ell) = \emptyset$
%        $(0\le i<N_{\ell,H})$}
      \ForEach{$0 \le m < N_{\ell,H}$}{
        % $c_m := \mathbf{C}_\ell^{\top}[m]$\;
        % \Solve $\mathbf{A}_{\ell}\, y_{m} = c_m$\;
        % $\mathbf{Y}_{\ell}[m] := y_m$\;
        $\mathbf{Y}_{\ell}[m] := A_{\ell}^{\text{inv}}( \mathbf{C}_\ell^{\top}[m])$\;
      }
      \tcc{-- compute inverse Schur complement --}
      $S^{\text{inv}}_{\ell} := (\mathbf{C}_\ell \mathbf{Y}_{\ell})^{-1}$ \tcp*{$\tt N_{\ell,H} \times N_{\ell,H}$ matrix}
      \tcc{-- compute correction for each coarse space function --}
      \tcc{-- which has a support on $K_{\ell}$ ~~~~~~~~~~~~~~~~~~~~~~~~--}
%      \ForEach{$\Phi_j \in \basis(V_H)$,
%        with $\supp(\Phi_j) \cap U(K_\ell) = \emptyset$
%        $(0\le j<N_{\ell,H})$}
      $c_d := \rows(T^H_{\ell})$\;
      $\mathbf{w}_{\ell} := $ Matrix$(c_d, N_{\ell,h})$ \tcp*{$\tt c_d \times N_{\ell,h}$ matrix}
      \lnl{alg:main:vertices}\ForEach{$0\le i < c_d$}
      {
        \tcc{compute $\tt \mathbf{w}_{\ell}[i] = \mathbf{A}^{\hspace{-2pt}-1}_\ell {\mathbf{r}}_\ell[i] - (\mathbf{A}^{\hspace{-2pt}-1}_\ell \mathbf{C}_\ell^{\top}) \mathbf{S}_{\ell}^{-1} (\mathbf{C}_\ell \hspace{2pt} \mathbf{A}^{\hspace{-2pt}-1}_\ell) \mathbf{r}_\ell[i]$}
        compute $q_i = A_{\ell}^{\text{inv}}(\mathbf{r}_\ell[i])$\tcp*{$\tt ~~ q_i := \hspace{-2pt}~~~~~~ \mathbf{A}_{\ell}^{-1} r_\ell[i]$}
        % \Solve $(C_i B^T) v = C_i q$\;
        % \Solve $v = S_{\text{inv}} C_i q$\;
        % \Solve $S_i g = v * C_i$\;
        % $w := q - g$\;
        %%% \Solve ??? $\mathbf{A}_{\ell} w := q - \mathbf{C}_\ell^{\top} (S_{\text{inv}} (\mathbf{C}_\ell q))$\;
        compute $\lambda_i = S^{\text{inv}}_{\ell}(\mathbf{C}_\ell q_i)$\tcp*{$\tt \boldsymbol{\lambda}_{\ell}[i] := \hspace{-2pt}~~~~\mathbf{S}_{\ell}^{-1} \mathbf{C}_\ell \mathbf{q}_{\ell}[i]$}
        compute $\mathbf{w}_{\ell}[i] = q_i - \mathbf{Y}_{\ell} \lambda_i$\tcp*{$\tt \mathbf{w}_{\ell}[i] := \mathbf{q}_{\ell}[i] - \mathbf{Y}_{\ell} \boldsymbol{\lambda}_{\ell}[i]$}
      }
      \tcc{update correction}
      $\mathbf{Q}_h := \mathbf{Q}_h + (\mathbf{T}_{\ell}^H)^{\top} \mathbf{w}_{\ell} \mathbf{R}_{\ell}^h$
    }
    \Return $\mathbf{Q}_h$\;
  }
  \caption{\label{alg:main}
    Computation of global corrector matrix $\mathbf{Q}_h$.}
\end{algorithm}
\begin{remark}
The algorithm can also be formulated with the boundary matrix $\mathbf{B}^H$. In this case, there is a small overhead in terms of the number of local problems to be solved, i.e. we solve problems for right hand sides ${\mathbf{r}}_\ell[i]$ that correspond to inactive coarse basis functions (basis functions belonging to Dirichlet-nodes).
\end{remark}
\subsection{The global corrector matrix}
To store the information that we obtained from the solutions of the local problems, we
introduce the global corrector matrix $\mathbf{Q}_h \in \R^{N_H \times N_h}$.
\begin{definition}[Global corrector matrix $\mathbf{Q}_h$]
\label{def-global-correction-matrix}
Recall the matrix $\mathbf{w}_{\ell} \in \R^{c_d \times N_{\ell,h}}$ introduced in Definition \ref{def-algebraic-formulation-og-local-problems} and recall that it is related to the correctors $\Qh^{K_{\ell}}(\Phi_{p_i(\ell)})$ (solving equation (\ref{local-corrector-problem-new})) via
$$\Qh^{K_{\ell}}(\Phi_{p_i(\ell)}) = \sum_{j=0}^{N_{\ell,h}-1} \mathbf{w}_{\ell}[i][j] \hspace{2pt} \phi_{\ell,j}.$$
With \eqref{relation-loc-cor-glob-cor} and the previously defined local restriction matrices $\mathbf{R}_{\ell}^h \in \R^{N_{\ell,h} \times N_{h}}$ and
$\mathbf{T}_{\ell}^H \in \R^{c_d \times N_{H}}$ we get the corrector matrix $\mathbf{Q}_h \in \R^{N_H \times N_h}$ via
$$\mathbf{Q}_h := \sum_{K_{\ell} \in \T_H} (\mathbf{T}_{\ell}^H)^{\top} \mathbf{w}_{\ell} \mathbf{R}_{\ell}^h.$$
\end{definition}
Hence, for any coarse function $\Phi_H \in \VHGD$, we can compute $Q_h(\Phi_H) \in \VhGD$ easily from $\mathbf{Q}_h$. For instance, let $\boldsymbol{\Phi} \in \R^{N_H}$ be the vector with entries $\boldsymbol{\Phi}[i]=\Phi_H(Z_i)$. Then we have
$$Q_h(\Phi_H) = \sum_{i=1}^{N_H-1} \left( \mathbf{Q}_h^{\top} \boldsymbol{\Phi} \right)\hspace{-2pt}[i] \hspace{2pt}\Phi_i.$$
The complete assembly of the global corrector matrix $\mathbf{Q}_h$ is summarized in Algorithm {\bf 1}.
%\textcolor{blue}
{Note that the global corrector matrix $\mathbf{Q}_h$ and the global stiffness matrix $\mathbf{A}_h$ are the only relevant matrices that need to be stored at this point. All other matrices are no longer required (or can be recomputed cheaply).}
\section{The LOD for linear elliptic problems}
\label{section-lod-for-elliptic-problems}
We are prepared to state the first full example for an application of the LOD.
Given $f\in L^{2}( \Omega) $, we seek the weak solution of
\begin{equation*}\label{eq:model}
  \begin{aligned}
    -\nabla \cdot \kappa \nabla u &= f\quad \text{in }\Omega, \\
    u & = 0 \quad \text{on }\Gamma_D, \\
    \kappa \nabla u \cdot n &= 0 \quad \text{on }\Gamma_N,
  \end{aligned}
\end{equation*}
i.e., we seek $u\in H^1_{\Gamma_D}(\Omega)$ that satisfies
\begin{equation}\label{e:modelproblem}
  \mathcal{A}\left(  u,v\right)=\int_{\Omega}\kappa \nabla u\cdot \nabla
  v =\int_{\Omega}fv=:\mathcal{F}( v) \quad\text{for all } v\in H^1_{\Gamma_D}(\Omega).
\end{equation}
\subsection{Method and convergence results}
\label{subsection-method-convergence}
With the definitions from Section~\ref{section-lod-initial-defs}, we can state the Local Orthogonal Decomposition method (LOD) for model problem \eqref{eq:model}.
%to account for the homogenous boundary condition on $\partial \Omega$ (i.e. we need to change the solution spaces). For that purpose, we define $W_h(U_k(K)) := W_{h}(U_k(K)) \cap H^1_0(\Omega)$.
% a homogenous Dirichlet boundary condition everywhere, hence we have $\Gamma=\partial \Omega$.
\begin{definition}[LOD approximation for problem \eqref{e:modelproblem}]
  \label{definition-loc-lod-approx}
  Recall Definition \ref{def-loc-lod} for a given localization parameter $k\in \mathbb{N}$.
  If $u_H \in \VHGD$ solves
  \begin{eqnarray}
    \label{lod-problem-eq}
    \int_{\Omega} \kappa \nabla (u_H+\Qh(u_H)) \cdot \nabla (\Phi_H+\Qh(\Phi_H))&=&
    \int_{\Omega} f (\Phi_H+\Qh(\Phi_H)) \quad \mbox{for all } \Phi_H \in \VHGD,
  \end{eqnarray}
  the final LOD approximation is given by $u_{\LOD}=u_H+\Qh(u_H)$.
\end{definition}
The Galerkin solution $u_h\in \VhGD$ which satisfies
\begin{equation}\label{e:modelref}
  \mathcal{A}(u_h,v) = F(v) \quad\text{ for all }v\in \VhGD
\end{equation}
can be considered as a reference solution in the sense that
$u_{\LOD}=u_H+\Qh(u_H)$ is constructed to approximate $u_h$ with a desired
accuracy of at least O$(H)$. This approximation quality can be
quantified:
\begin{theorem}[A priori error estimate]
  \label{theorem-apriori-linear}Assume that the localization parameter fulfills $k \gtrsim m |\log(H)|$ for some $m \in \mathbb{N}$. Then, there exists a
  positive constant $C$ that depends on the space dimension $d$, on
  $\Omega$, $\gmin$, $\gmax$ and interior angles of the partitions,
  but not on the mesh sizes $H$ and $h$, such that
  \begin{align*}
    \| u_h - (u_H+\Qh(u_H))\|_{L^2(\Omega)} &\le C (H + H^{r m})^2 \quad \mbox{and}\\
    \| u_h - (u_H+\Qh(u_H))\|_{H^1(\Omega)} +  \| u_h - u_H\|_{L^2(\Omega)} &\le C (H + H^{r m}),
  \end{align*}
  for some constant $r>0$ that depends linearly on the square root of
  the contrast.
\end{theorem}
The theorem was proved in \cite{MaP14,HeP13,HeM14}. Practically, numerical
experiments indicate that the choice $m\in\{1,2,3\}$ typically yields
good results even for high contrast cases \cite{2016arXiv160106549P}. We refer to the numerical experiments in
\cite{HeM14,MaP14}.
\subsection{Assembly and solution of the global problem}
\label{subsection-assembling-solving-global-problem}
\subsubsection{Formal description}
Assume that all local problems are solved (i.e. solved for every $K _{\ell}\in
\T_H$ and every coarse basis function $\Phi_{p_i(\ell)}$ with support on
$K_{\ell}$, where $0\le i < c_d$) so that
$\Qh^{K_{\ell}}(\Phi_{p_i(\ell)})$
is available and we can write
\begin{align*}
\Qh(\Phi_H)=\sum_{K_{\ell} \in \T_H} \sum_{i=0}^{c_d-1} \Phi_H(Z_{p_i(\ell)}) \Qh^{K_{\ell}}(\Phi_{p_i(\ell)})
\end{align*}
for any $\Phi_H \in \VHGD$. Consequently we can assemble the (global) LOD stiffness matrix $\mathbf{A}_H^{\LOD}\in \R^{N_H \times N_H}$ that is given by the entries
\begin{align*}
  \mathbf{A}_H^{\LOD}[m][n] :=
  \begin{cases}
  \mathcal{A}( \Phi_{n}+\Qh(\Phi_n), \Phi_{m}+\Qh(\Phi_m) ) \qquad &\mbox{for } \enspace Z_m,Z_n \in \mathcal{N}_H \setminus \Gamma_D,\\
  \hspace{75pt}0 &\mbox{else}
  \end{cases}
\end{align*}
and LOD load vector $\mathbf{f}_H \in \R^{N_H}$ given by
\begin{align*}
  \mathbf{f}_H[m] :=
  \begin{cases}
  ( f, \Phi_{m}+\Qh(\Phi_m) )_{L^2(\Omega)} \qquad &\mbox{for } \enspace Z_m \in \mathcal{N}_H \setminus \Gamma_D,\\
  \hspace{35pt}0 &\mbox{else.}
  \end{cases}
\end{align*}
%With that, we update $\mathbf{A}_H^{\LOD}$ and $\mathbf{f}_H$ via
%$$\mathbf{A}_H^{\LOD} := \mathbf{B}^H \mathbf{A}_H^{\LOD} \mathbf{B}^H \qquad \mbox{and} \qquad \mathbf{f}_H := \mathbf{B}^H \mathbf{f}_H.$$
With that, the algebraic version of \eqref{lod-problem-eq} hence reads: find ${\mathbf{u}}_H^{\LOD} \in \R^{N_H}$ with
$$\mathbf{A}_H^{\LOD} {\mathbf{u}}_H^{\LOD} = \mathbf{f}_H.$$
Once this is solved, the final LOD approximation is given by
\begin{align*}
 u_{\LOD} = \sum_{m=0}^{N_H-1}  {\mathbf{u}}_H^{\LOD}[m] (\Phi_{m}+\Qh(\Phi_m)) \in \VhGD.
\end{align*}
%\textcolor{blue}
{The condition number of the LOD system matrix $\mathbf{A}_H^{\LOD}$
  is of order $1/H^2$, i.e. of the same order as the condition number
  for standard finite elements on the coarse scale. A corresponding
  estimate for the condition number of $\mathbf{A}_H^{\LOD}$ is for
  instance given in \cite[Lemma 1]{Peterseim.Schedensack:2016}.}
Depending on the coarse space size one can employ a direct solver or a
standard iterative solver, like the multigrid method.
\subsubsection{Algorithmic realization}
%So far we only presented a formal algebraic description
%
  Next we describe an efficient algorithmic realization of how to assemble and solve the global problem.
When all local problems are solved, the global corrector matrix $\mathbf{Q}_h \in \R^{N_H \times N_h}$ is available (cf. Definition \ref{def-global-correction-matrix}). From $\mathbf{Q}_h$,
%$\mathbf{B}^H$,
$\mathbf{B}^H$ and the projection matrix $\mathbf{P}_h$ we hence obtain the global LOD system matrix by matrix multiplication
$$\mathbf{A}_H^{\LOD} = \mathbf{B}^H (\mathbf{P}_h+\mathbf{Q}_h) \mathbf{A}_h (\mathbf{P}_h+\mathbf{Q}_h)^{\top} \mathbf{B}^H.$$
Similarly, we get the load vector $\mathbf{f}_H$ by
$$\mathbf{f}_H = \mathbf{B}^H (\mathbf{P}_h+\mathbf{Q}_h) \mathbf{f}_h,$$
where $\mathbf{f}_h \in \R^{N_h}$ denotes the classical FEM load vector with entries $\mathbf{f}_h[i]=(f,\phi_i)_{L^2(\Omega)}$ for $0\le i < N_h$. Now we can solve for ${\mathbf{u}}_H^{\LOD} \in \R^{N_H}$ with
$$\mathbf{A}_H^{\LOD} {\mathbf{u}}_H^{\LOD} = \mathbf{f}_H$$
and obtain the final coefficient vector ${\mathbf{u}}_h^{\LOD} \in \R^{N_h}$ of our LOD approximation by
$${\mathbf{u}}_h^{\LOD} := (\mathbf{P}_h+\mathbf{Q}_h)^{\top} {\mathbf{u}}_H^{\LOD}.$$
The procedure is summarized in Algorithm {\bf 2}.
\begin{algorithm}
  \DontPrintSemicolon
  \KwData{$\mathcal{T}_h$, $\mathcal{T}_H$,\tcp*{fine / coarse mesh}}
  \AddData{$N_h$, $N_H$, \tcp*{fine / coarse space size}}
  \AddData{$\mathbf{P}_h$
%      = \left(
%      \begin{matrix}
%        \Phi_0(z_0) & \cdots & \Phi_0(z_{N_h-1}) \\
%        \vdots & \ddots & \vdots \\
%        \Phi_{N_H-1}(z_0) & \cdots & \Phi_{N_H-1}(z_{N_h-1})
%      \end{matrix}
%    \right)$
    \tcp*{$\tt N_H \times N_h$ interp.\hspace*{1.5ex}matrix in \eqref{projection-matrix}}}
  \AddData{$\mathbf{A}_h$\tcp*{$\tt N_h \times N_h$ stiffness matrix}}
  \AddData{$\mathbf{f}_h$\tcp*{$\tt N_h \times 1$ (fine) load vector}}
  \AddData{$\mathbf{B}^H$\tcp*{$\tt N_H \times N_H$ boundary correction matrix}}
  \AddData{$\mathbf{Q}_h$\tcp*{$\tt N_H \times N_h$ global corrector matrix}}
  % \hfill\parbox{4.5cm}{\tcp*[f]{$\tt N_H \times N_h$ interp.}\par\tcp*[f]{~~~~~~ matrix~}}}
  \Fn{solveLODSystem}{
    compute $\mathbf{A}_H^{\LOD} = \mathbf{B}^H (\mathbf{P}_h+\mathbf{Q}_h) \mathbf{A}_h (\mathbf{P}_h+\mathbf{Q}_h)^{\top} \mathbf{B}^H$
    \tcp*{$\tt N_H \times N_H$ LOD sys.mat.}
    compute $\mathbf{f}_H = \mathbf{B}^H (\mathbf{P}_h+\mathbf{Q}_h) \mathbf{f}_h$\tcp*{$\tt N_H \times 1$\hspace*{4pt} LOD rhs.vec.}
    \Solve $\mathbf{A}_H^{\LOD} {\mathbf{u}}_H^{\LOD} = \mathbf{f}_H$\tcp*{solve \hspace*{5pt}LOD system\hspace*{11.5pt}}
    compute ${\mathbf{u}}_h^{\LOD} := (\mathbf{P}_h+\mathbf{Q}_h)^{\top} {\mathbf{u}}_H^{\LOD}$\tcp*{$\tt N_h \times 1$\hspace*{4pt} LOD solution}
    \Return ${\mathbf{u}}_h^{\LOD}$\;
  }
  \caption{Computation of the final LOD approximation ${\mathbf{u}}_h^{\LOD}$.}
\end{algorithm}
\begin{remark}
Typically it is possible to replace the right hand side in \eqref{lod-problem-eq} by $(f,\Phi_H)_{L^2(\Omega)}$ without a significant loss in accuracy. In particular if $f$ is a slow variable. Consequently, we do not require the correctors any longer to compute the load vector in the global LOD system. This can turn out to be an immense computational advantage if the LOD system has to be solved for several source terms $f$. In this case, we can fully reuse $\mathbf{A}_H^{\LOD}$ and quickly assemble $\mathbf{f}_H$ with entries $\mathbf{f}_H[i]=(f,\Phi_i)_{L^2(\Omega)}$ (only involving coarse basis functions).
\end{remark}
\begin{remark}
In the case that the coefficient $\kappa$ has certain structural properties (such as periodicity), it might be possible to only assemble some of the local correctors $\Qh^K$ and reuse them on the different location in $\Omega$. This is possible if a corrector can be expressed as a rotation and translation of another corrector. With that the computational complexity can be decreased significantly. This has been exploited in the context of acoustic scattering in \cite{GaP15}.
\end{remark}
\subsection{Petrov-Galerkin version of the method}
In some cases it can happen that the fine space $V_h$ is so large that the storing of the full system matrix $\mathbf{A}_h$ (respectively the storing of the corrector matrix $\mathbf{Q}_h$) becomes too memory demanding. In such cases we cannot afford the multiplication of $N_h \times N_h$ matrices as frequently done in Algorithm 1 and 2. To overcome the issue that the size of $V_h$ exceeds the computational resources, a Petrov-Galerkin (PG) formulation of the LOD can be used. This method allows an on-the-fly assembling of the LOD system matrix $\mathbf{A}_H^{\LOD}$ on the expense that we lose the symmetry. Let us start with describing the Petrov-Galerkin LOD from the analytical point of view.
\subsubsection{Description and properties of the PG-LOD}
We use the notation introduced earlier in this section.
\begin{definition}[PG LOD approximation]
  \label{definition-loc-pg-lod-approx}
Let $k\in \mathbb{N}$ be fixed and let
$\Qh: \VHGD \rightarrow \VhGD$ denote the
corresponding correction operator as in Definition \ref{def-loc-lod}.
%Again, we denote $\Rh(\Phi_H):=\Phi_H+\Qh(\Phi_H)$ for $\Phi_H \in \VHGD$.
If $u_H^{\PG} \in \VHGD$ solves
  \begin{eqnarray}
    \label{pg-lod-problem-eq}
    \int_{\Omega} \kappa \nabla ( u_H^{\PG} + \Qh(u_H^{\PG})) \cdot \nabla \Phi_H&=&
    \int_{\Omega} f  \Phi_H \quad \mbox{for all } \Phi_H \in \VHGD,
  \end{eqnarray}
the final Petrov-Galerkin LOD approximation is given by $u^{\PG-\LOD}:=u_H^{\PG}+\Qh(u_H^{\PG})$.
\end{definition}
Obviously, the standard formulation of the LOD only differs from the PG formulation by the choice of test functions in (\ref{pg-lod-problem-eq}). In particular, the solving of the local corrector problems is identical for both methods. From the analytical point of view, the change of test functions in (\ref{lod-problem-eq}) does not have a crucial influence. We still have well-posedness of the PG-LOD solution and the obtained convergence rates are the same as for the original method. We summarize the corresponding main result in the following theorem, which is proved in \cite{EGH15}.
\begin{theorem}[A priori error estimate for the PG-LOD]
  \label{theorem-apriori-linear-pg}
  %Assume that $\kappa \in
  %L^\infty(\Omega,\R^{d\times d}_{\mbox{\tiny sym}})$, $f \in
  %L^2(\Omega)$ and
  Assume that $k \gtrsim m |\log(H)|$ for some $m \in \mathbb{N}$. Furthermore, let the positive constants $C$ and $r$ be as in Theorem \ref{theorem-apriori-linear}. Then the left side of (\ref{pg-lod-problem-eq}) represents a coercive bilinear form on $\VHGD$, i.e.
  \begin{align}
\label{coercivity-pg-lod}  \int_{\Omega} \kappa \nabla (\Phi_H + \Qh(\Phi_H)) \cdot \nabla \Phi_H \ge C (\alpha - C H^{rm} ) \| \Phi_H \|_{H^1(\Omega)}^2.
  \end{align}
 Consequently, problem (\ref{pg-lod-problem-eq}) is well-posed and the PG-LOD approximation fulfills the same error estimates as the standard LOD approximation, i.e. we have
  \begin{align*}
    \| u_h - u_H^{\PG} - \Qh(u_H^{\PG})\|_{L^2(\Omega)} &\le C (H + H^{r m})^2 \quad \mbox{and}\\
    \| u_h - u_H^{\PG} - \Qh(u_H^{\PG}) \|_{H^1(\Omega)} +  \| u_h - u_H^{\PG}\|_{L^2(\Omega)} &\le C (H + H^{r m}).
  \end{align*}
\end{theorem}
For even sharper results in $L^2$ we refer to \cite{2017arXiv170208858G}.
\begin{remark}[Relevance of $L^2$-approximations]
\label{pg-lod-remark-l2-approx}Theorem \ref{theorem-apriori-linear-pg} contains $L^2$- and $H^1$ error estimate for the full RG-LOD approximation $u_H^{\PG}+\Qh(u_H^{\PG})$. However, to compute it from the coarse part $u_H^{\PG}$, we need to know the operator $\Qh$. But recall that the algebraic version of $\Qh$ is represented by the corrector matrix $\mathbf{Q}_h\in \R^{N_h \times N_h}$, which is of the same size (and even less sparse) than the global stiffness matrix $\mathbf{A}_h \in \R^{N_h \times N_h}$. So if we do not have the capacities to store $\mathbf{A}_h$, neither do we have the capacities to store $\mathbf{Q}_h$. Consequently, even though $u_H^{\PG}$ might be available, $\Qh(u_H^{\PG})$ is typically not. Hence, our final approximation $u_H^{\PG}$ is only an $L^2$-approximation, instead of a full $H^1$-approximation (as e.g. $u_H+\Qh(u_H)$). Hence, the relevant estimate that remains from Theorem \ref{theorem-apriori-linear-pg} is the $L^2$-error estimate
$$\| u_h - u_H^{\PG}\|_{L^2(\Omega)} \le C (H + H^{r m}).$$
However, note that given the RG-LOD solution $u_H$ we could go easily back to the local problems and recompute with known right hand side to form the full fine scale solution without having to store the correctors $\Qh(u_H)$. This strategy allows to obtain $H^1$-approximations through local post-processing.
\end{remark}
\subsubsection{Computational advantages and disadvantages}
Let us now describe the computational advantages and disadvantages of the Petrov-Galerkin formulation. We start with the advantages to demonstrate how the PG formulation overcomes the capacity issues.
$\\$
{\bf Advantages.} The basic advantage of the PG-LOD is that matrices of size $N_h \times N_h$ have to be handled at no point. Operations either involve $N_{\ell,h} \times N_{\ell,h}$-matrices ($N_{\ell,h}$ is the number of fine nodes in the patch $U_k(K_{\ell})$) or they involve $N_{H} \times N_{H}$-matrices (where $N_H$ denotes the size of $V_H$). The reason why this is possible is that no corrector-to-corrector communication is required for the PG-LOD. For instance, in order to assemble the system matrix that is associated with standard LOD (cf. (\ref{lod-problem-eq})), we need to compute entries such as
\begin{align*}
\int_{\Omega} \kappa \nabla \left( \Phi_i + \Qh(\Phi_i) \right) \cdot  \nabla \left( \Phi_j + \Qh(\Phi_j) \right)
\end{align*}
for two coarse basis functions $\Phi_i$ and $\Phi_j$. It is impossible to compute this entry without knowing both $\Qh(\Phi_i)$ and $\Qh(\Phi_j)$ at the same time. Consequently correctors must be stored so that they can communicate with each other. For the PG-LOD, system matrix entries are always of the structure
\begin{align*}
\int_{\Omega} \kappa \nabla \left( \Phi_i + \Qh(\Phi_i) \right) \cdot  \nabla \Phi_j
=
\underset{K \subset \mbox{\rm supp}(\Phi_i)}{\sum_{K \in \T_H}} \int_{U_{k}(K)} \kappa \nabla \left( \Phi_i + \Qh^{K}(\Phi_i) \right) \cdot  \nabla \Phi_j,
\end{align*}
which can be assembled (respectively updated) after a corrector
$\Qh^{K}(\Phi_i)$ is computed. If desired, $\Qh^{K}(\Phi_i)$ can be
immediately deleted after this.
%\textcolor{blue}
{Global matrices such as the corrector matrix $\mathbf{Q}_h$ or the stiffness matrix $\mathbf{A}_h$ do neither have to be stored nor explicitly computed. Only a sparse global system matrix $\mathbf{A}_H^{\PG-\LOD}$ of size $N_H \times N_H$ is required.} Consequently, the storage requirements are significantly lower for the Petrov-Galerkin version.
$\\$
{\bf Trade-offs.} In comparison to (\ref{lod-problem-eq}), we observe that the PG-LOD system given by (\ref{pg-lod-problem-eq}) can no more be represented by a symmetric matrix, which formally excludes the usage of certain efficient algebraic solvers that rely on symmetry. However, having a closer look, we see that the method only suffers from a mild loss of symmetry in the sense that the PG-LOD is still symmetric if there is no localization and that the lack of symmetry can be hence quantified by the exponential decay property. A symmetric approximation can be for instance obtained by using $\frac{1}{2}\mathbf{A}_H^{\PG-\LOD}+\frac{1}{2}(\mathbf{A}_H^{\PG-\LOD})^{\top}$ for the system matrix.
%The computational cost for solving (\ref{pg-lod-problem-eq}) might hence slightly increase (compared to \eqref{lod-problem-eq}).
The second trade-off is rather subtle. Theorem \ref{theorem-apriori-linear-pg} predicts a coercivity constant that can be disturbed by a term of order $\mathcal{O}(H^{rm})$. Even though this seems to be mostly unproblematic for small $H$, there is formally no guarantee that $(\alpha - C H^{rm} )$ is always positive. This can only be guaranteed by a numerical investigation of the eigenvalues. However, we also note that non-positivity has never been observed in numerical experiments. So far it seems that the result (\ref{coercivity-pg-lod}) is not yet optimal and the coercivity appears to be always fulfilled in practical applications.
Another trade-off was already mentioned in Remark \ref{pg-lod-remark-l2-approx}.
%(it is not a real disadvantage since the original LOD cannot overcome this problem either).
If we are in a scenario where the PG-LOD is used to decrease the memory demand, then the correctors $\Qh$ will not be stored. Hence we will not be able to compute $\Qh(u_H^{\PG})$ from $u_H^{\PG}$ and have to be content with an $L^2$-approximation of $u_h$. In many applications this is enough. Especially when considering that $u_H^{\PG}$ can be stored with significantly lower costs than the full fine scale approximation $u_H^{\PG}+\Qh(u_H^{\PG})$. If local fine scale information is required afterwards by a user, it is possible to perform a local "real time" post-processing where the missing fine-scale information is only computed in the relevant region.
\subsubsection{Realization}
Even though the (algebraic) realization differs only slightly from the realization of the classical LOD, these differences are essential. The local corrector problems are computed in the same way as before, however, instead of storing their solutions in a global corrector matrix, their contributions are directly added to the global PG-LOD system matrix and can be immediately deleted afterwards.
$\\$
To summarize the basic procedure, let us fix a coarse element $K_{\ell} \in \T_H$ and a corresponding coarse layer patch $\mathcal{U}_{\ell}=U_{k}(K_{\ell})$. The following step has to be repeated for every $K_{\ell}$.
For every coarse basis function $\Phi_{p_i(\ell)}$ (with $0\le i < c_d$ and only if $Z_{p_i(\ell)} \not\in \Gamma_D$) we solve for the local corrector
for $\Qh^{K_{\ell}}(\Phi_{p_i(\ell)})$ according to (\ref{local-corrector-problem-new}).
%After the $c_d$ correctors are computed we update the system matrix $\mathbf{A}_H^{\PG-\LOD}$.
After $\Qh^{K_{\ell}}(\Phi_{p_i(\ell)})$ is computed, we update the system matrix $\mathbf{A}_H^{\PG-\LOD}$.
To do that we visit every coarse basis function $\Phi_m \in \VHGD$ with $\mbox{\rm supp}(\Phi_m)\cap\mathcal{U}_{\ell}\neq \emptyset$. For each such $\Phi_m$ we can make the update
\begin{align*}
\mathbf{A}_H^{\PG-\LOD}[m][p_i(\ell)] := \mathbf{A}_H^{\PG-\LOD}[m][p_i(\ell)] + \int_{\mathcal{U}_{\ell}} \kappa \nabla \left( \Phi_{p_i(\ell)} + \Qh^{K_{\ell}}(\Phi_{p_i(\ell)}) \right)
\cdot \nabla \Phi_m.
\end{align*}
When all loops have terminated we can incorporate the homogenous boundary condition by multiplying $\mathbf{A}_H^{\PG-\LOD}$ at both sides with the boundary correction matrix $\mathbf{B}^H$, i.e. $\mathbf{A}_H^{\PG-\LOD}:=\mathbf{B}^H \mathbf{A}_H^{\PG-\LOD} \mathbf{B}^H$. In total, we computed the correct PG-LOD system matrix with entries
\begin{align*}
\mathbf{A}_H^{\PG-\LOD}[m][n]
&= \int_{\Omega} \kappa \nabla \left( \Phi_n + \Qh(\Phi_n) \right) \cdot \nabla \Phi_m\\
&= \underset{K_{\ell} \subset \mbox{\rm supp}(\Phi_n)}{\sum_{K_{\ell} \in \T_H}}
%\int_{\mbox{\rm supp}(\Phi_m)}
\int_{\mathcal{U}_{\ell}}
 \kappa \nabla \left( \Phi_n + \Qh^{K_{\ell}}(\Phi_n) \right) \cdot \nabla \Phi_m.
\end{align*}
Observe that we could generate a new local fine mesh for every $\mathcal{U}_{\ell}$. Basically, there is no need for a global fine mesh $\mathcal{T}_h$.
The remaining procedure is straightforward. Since the right hand side of (\ref{pg-lod-problem-eq}) only involves standard coarse functions we can set $\mathbf{f}_H^{\PG-\LOD}:=\mathbf{B}^H \mathbf{f}_H$, where $\mathbf{f}_H \in \mathbb{R}^{N_H}$ denotes the standard coarse load vector with entries $\mathbf{f}_H[m]= \int_{\Omega} f \Phi_m$. Consequently it only remains to solve for ${\mathbf{u}}_H^{\PG-\LOD} \in \R^{N_H}$ with
$$\mathbf{A}_H^{\PG-\LOD} {\mathbf{u}}_H^{\PG-\LOD} =\mathbf{f}_H^{\PG-\LOD}.$$
The final approximation is given by
$$
u_H^{\PG} = \sum_{m=0}^{N_H-1} {\mathbf{u}}_H^{\PG-\LOD}[m] \Phi_m.
$$
Accordingly modified algorithms can be formulated analogously to the algorithms presented for the standard (symmetric) LOD.
\section{The treatment of rough boundary data and sources}
\label{section5:LOD-including-source-and-boundary-terms}

In this section we discuss how we can incorporate nonhomogeneous boundary conditions in the LOD approximations and how we can treat regions in which the source term $f$ becomes close to singular.
\subsection{Model problem and discretization}
We consider the following problem with mixed Dirichlet and Neumann boundary conditions. Find $u$ with
\begin{equation*}\label{eq:model2}
  \begin{aligned}
    -\nabla \cdot \kappa \nabla u &= f\quad \hspace{6pt}\text{in }\Omega, \\
    u & = g \quad \hspace{7pt}\text{on }\Gamma_D, \\
    \kappa \nabla u \cdot n &= q \quad \hspace{7pt}\text{on }\Gamma_N.
  \end{aligned}
\end{equation*}
In addition to the previous assumptions on $\Omega$ and $\kappa$ we assume that the Dirichlet boundary values
fulfill $g \in H^{\frac{1}{2}}(\Gamma_D)$ and that the Neumann boundary values fulfill $q \in L^2(\Gamma_N)$.
The weak formulation of problem (\ref{eq:model2}) reads: find $u \in H^1(\Omega)$, with $T_D(u)=g$, such that
\begin{align*}
\int_{\Omega} \kappa \nabla u \cdot \nabla v = \int_{\Omega} f v + \int_{\Gamma_N} q v \quad \mbox{for all } v\in H^1_{\Gamma_D}(\Omega).
\end{align*}
Assume that $g$ is sufficiently regular so that point evaluations are
possible. Then we can define $g_H \in V_H$ as the function that is
uniquely determined by the nodal values $g_H(Z)=g(Z)$ for all $Z \in
\mathcal{N}_H \cap \Gamma_D$ and $g_H(Z)=0$ for all $Z \in
\mathcal{N}_H \setminus \Gamma_D$. Using this, we define the (fine
scale) Dirichlet extension $g_h \in V_h$ uniquely by the nodal values
$g_h(z)=g(z)$ for all $z \in \mathcal{N}_h \cap \Gamma_D$ and
$g_h(z)=g_H(z)$ for all $z \in \mathcal{N}_h \setminus \Gamma_D$. With
this, we avoid degeneracy of $g_h$ for $h$ tending to zero.
%\textcolor{blue}
{Note that the extension $g_h$ needs to be explicitly constructed.} The reference problem reads: find $U_h \in \VhGD$ with
\begin{align}
\label{reference-problem-eq} \mathcal{A}( U_h , v_h) = \int_{\Omega} f v_h - \int_{\Omega} \kappa \nabla g_h \cdot \nabla v_h +\int_{\Gamma_N} q v_h \quad \mbox{for all } v_h \in \VhGD.
\end{align}

The final fine scale approximation is then given by $u_h:=U_h+g_h \in V_h$. Now observe that problem (\ref{reference-problem-eq}) is basically of the same structure as the homogenous problem \eqref{e:modelref}. This suggest to apply the same methodology as before.
Unfortunately, the correctors introduced in Definition \ref{def-loc-lod} might not be sufficient to construct accurate approximations, if e.g. the Dirichlet boundary condition is highly oscillatory. The slight difference that the right hand side is no longer purely represented by an $L^2$-function $f$ (but by a less regular functional
%$\mathcal{F}$
which is only in the dual space of $\VhGD$) makes it necessary to introduce additional correctors to preserve the previous convergence rates. We call these new correctors {\it source correctors}.
\subsection{Source correctors}
In this section we introduce {\it source term correctors}. They are defined analogously to the correctors $\Qh$. Their purpose is to captured oscillatory effects that are produced by a general source. For that purpose, we split the right hand side of (\ref{reference-problem-eq}) into two parts. One part (we shall denote by $\mathcal{F}$) that has basically a coarse scale structure and that can be considered as harmless if ignored by the fine grid, and a second part (we shall denote by $\mathcal{F}^{\text{s}}$) which might have a considerable influence on the oscillations of $u_h$. Hence we let $\mathcal{F} : H^1(\Omega) \rightarrow \R$ and $\mathcal{F}^{\text{s}} : H^1(\Omega) \rightarrow \R$ be source functionals such that
\begin{align*}
\mathcal{F} (v_h) + \mathcal{F}^{\text{s}} (v_h)  = \int_{\Omega} f v_h - \int_{\Omega} \kappa \nabla g_h \cdot \nabla v_h +\int_{\Gamma_N} q v_h \quad \mbox{for all } v_h \in \VhGD.
\end{align*}
We only wish to introduce additional correctors for the $\mathcal{F}^{\text{s}}$-contribution. It can incorporate source terms and boundary conditions and we assume that it is of the structure
\begin{align*}
\mathcal{F}^{\text{s}}(v) = \int_{\Omega} \eta_1 v + \kappa \nabla \eta_2 \cdot \nabla v + \int_{\partial \Omega} \eta_3 v
\end{align*}
with some given $\eta_1 \in L^2(\Omega)$, $\eta_2 \in H^1(\Omega)$ and $\eta \in H^{1/2}(\partial \Omega)$.
Typical choices would be
\begin{align*}
\mathcal{F}^{\text{s}}(v) &= \int_{\Omega} \kappa \nabla g_h \cdot \nabla v + \int_{\Gamma_N} q v \quad \mbox{(boundary source) or},\\
\mathcal{F}^{\text{s}}(v) &= \int_{\Omega} fv + \kappa \nabla g_h \cdot \nabla v + \int_{\Gamma_N} q v \quad \mbox{(total source).}
\end{align*}
We define the localization of $\mathcal{F}^{\text{s}}$ to a coarse element $K\in \T_H$ by
\begin{align*}
\mathcal{F}^{\text{s}}_K(v) := \int_{K} \eta_1 v + \kappa \nabla \eta_2 \cdot \nabla v + \int_{\partial \Omega \cap K} \eta_3 v.
\end{align*}
With that, we can define local source correctors.
\begin{definition}[Source term correctors]
  \label{definition-lod-source-term-correctors}
  Let $\mathcal{F}^{\text{s}}$ be fixed according to the previous discussion. For a given positive $k\in \mathbb{N}$ and $K\in \T_H$ we define the local source corrector
  $\QFh^{K} \in W_h(U_k(K))$ as the solution of
\begin{align}
\label{local-source-corrector-problem}\int_{U_k(K)} \kappa \nabla \QFh^{K} \cdot \nabla w_h = - \mathcal{F}^{\text{s}}_K(w_h) \qquad \mbox{for all } w_h \in  W_h(U_k(K)).
\end{align}
The corresponding global corrector is given by
\begin{align*}
\QFh:=\sum_{K\in \T_H} \QFh^{K}.
\end{align*}
\end{definition}
Note that it is desirable that $\mathcal{F}^{\text{s}}$ only contains locally supported sources, i.e. that $\mathcal{F}^{\text{s}}_K(w_h)=0$ for most of the coarse elements $K$. The more elements with $\mathcal{F}^{\text{s}}_K(w_h)\neq0$, the more local problems to solve. Therefore $\mathcal{F}^{\text{s}}$ typically only contains boundary terms or parts of $f$ in a small region where $f$ might become close to singular.
\subsection{Formulation of the method and error estimates}
Using Definition \ref{def-loc-lod} and \ref{definition-lod-source-term-correctors} we propose the following LOD approximation.
\begin{definition}[LOD approximation for boundary value problems]$\\$
\label{definition-loc-lod-approx-bc}For fixed $\mathcal{F}^{\text{s}}$ and fixed $k\in \mathbb{N}$
the LOD approximation to problem (\ref{reference-problem-eq}) is given by $u_{\LOD}:=U_H+\Qh(U_H)-B_{h}$, where $U_H \in \VHGD$ solves:
\begin{eqnarray}
\label{lod-problem-eq-mixed-bc} \lefteqn{\int_{\Omega} \kappa \nabla (\Id+\Qh)(U_H) \cdot \nabla (\Id+\Qh)(\Phi_H)}\\
\nonumber&=& (\mathcal{F}+\mathcal{F}^{\text{s}})((\Id+\Qh)(\Phi_H)) + \int_{\Omega} \kappa \nabla \QFh \cdot \nabla (\Id+\Qh)(\Phi_H)
\end{eqnarray}
for all $\Phi_H \in \VHGD$.
\end{definition}
The following a priori error estimate is proved in \cite{HeM14}.
\begin{theorem}
  \label{theorem-apriori-linear-mixed-bc}Assume that $k \gtrsim m |\log(H)|$ for some $m \in \mathbb{N}$. Furthermore, let
\begin{align*}
\mathcal{F}^{\text{s}}(v) = \int_{\Omega} \kappa \nabla g_h \cdot \nabla v + \int_{\Gamma_N} q v \qquad \mbox{and} \qquad
\mathcal{F}(v) &= \int_{\Omega} fv.
\end{align*}
Recall that we compute the source corrector $\QFh$ only with respect to $\mathcal{F}^{\text{s}}$. Then the LOD approximation $u_{\LOD}$ introduced in Definition \ref{definition-loc-lod-approx-bc} fulfills the estimates
  \begin{align*}
    \| u_h - u_{\LOD}\|_{L^2(\Omega)} &\le C (H + H^{r m})^2 \quad \mbox{and}\quad
    \| u_h - u_{\LOD}\|_{H^1(\Omega)} \le C (H + H^{r m}),
  \end{align*}
 where $C$ and $r>0$ are as in Theorem \ref{theorem-apriori-linear}.
\end{theorem}
\subsection{Algebraic realization}
The algebraic realization is straightforward following the ideas presented in Section~\ref{section-lod-for-elliptic-problems}. We only need to solve one additional linear elliptic problem more for each coarse element $K_{\ell}\in\T_H$ with $\mathcal{F}^{\text{s}}_{K_{\ell}}\neq 0$
%, we require a modification of the boundary correction matrix $\mathbf{B}^H$
and we need to assemble an additional vector that stores the entries $\int_{\Omega} \kappa \nabla \QFh \cdot \nabla (\Id+\Qh)(\Phi_j)$.
\subsubsection{Solving of the additional local problem}
Let us fix a coarse element $K _{\ell}\in \T_H$. First recall that for every coarse basis function $\Phi_{p_i(\ell)}$ ($0\le i < c_d$, see Section~\ref{subsubsection-assembling-load-vector}) we need to solve for $w_{\ell,i}:=\Qh^{K_{\ell}}(\Phi_{p_i(\ell)})\in \Vhell$ with $\Ic(w_{\ell,i}) =0$ and
  \begin{align*}
    \int_{\mathcal{U}_{\ell}} \kappa \nabla w_{\ell,i}\cdot \nabla w_h = - \int_{K_{\ell}} \kappa
    \nabla \Phi_{p_i(\ell)} \cdot \nabla w_h
\end{align*}
for all $w_h \in \Vhell$ with $\Ic(w_h) =0$. Now, if $\mathcal{F}^{\text{s}}_{K_{\ell}}\neq 0$ we also need to solve
for $\hat{w}_{\ell} := \QFh^{K_{\ell}}\in \Vhell$ with $\Ic(\hat{w}_{\ell}) =0$ and
\begin{align}
 \label{loc-source-corrector-problem-new} \int_{\mathcal{U}_{\ell}} \kappa \nabla \hat{w}_{\ell} \cdot \nabla w_h = - \mathcal{F}^{\text{s}}_K(w_h)
\end{align}
for all $w_h \in \Vhell$ with $\Ic(w_h) =0$. However, both problems only differ in their source terms, where the inverse of the Schur-complement matrix $\mathbf{S}_{\ell}^{-1}$ is already precomputed for solving the original corrector problems. Let us introduce a notation for the algebraic version of (\ref{loc-source-corrector-problem-new}).
\begin{definition}[Algebraic formulation of (\ref{loc-source-corrector-problem-new})]
\label{def-algebraic-formulation-source-corrector-problem-new}
Let the load vector $\hat{\mathbf{r}}_{\ell} \in \R^{N_{\ell,h}}$ be given by the entries
$
 \hat{\mathbf{r}}_{\ell}[j] := - \mathcal{F}^{\text{s}}_K( \phi_{\ell,j} ) \enspace \mbox{for } 0\le j < N_{\ell,h}.
$
The algebraic version of problem (\ref{loc-source-corrector-problem-new}) is the following saddle point problem. Find the tuple $(\hat{\mathbf{w}}_{\ell},\hat{\boldsymbol{\lambda}}_{\ell}) \in \R^{N_{\ell,h}} \times \R^{N_{\ell,H}}$ with
\begin{align}
\label{local-corrector-problem-source-new-algebraic}  \mathbf{A}_{\ell} \hspace{2pt} \hat{\mathbf{w}}_{\ell} + \mathbf{C}_\ell^{\top} \hat{\boldsymbol{\lambda}}_{\ell} &= \hat{{\mathbf{r}}}_\ell\\
\nonumber  \mathbf{C}_\ell \hspace{2pt} \hat{\mathbf{w}}_{\ell} &= 0.
\end{align}
Here, $\hat{\mathbf{w}}_{\ell}$ is the coefficient vector for the solution $\hat{w}_{\ell} = \QFh^{K_{\ell}}$.
\end{definition}
Exploiting the notation from Section~\ref{subsubsection-assembling-solving-loc-problems}, we obtain $\hat{\mathbf{w}}_{\ell}$ in three steps:
\begin{enumerate}
\item Solve for $\hat{\mathbf{q}}_{\ell} \in \R^{N_{\ell,h}}$ with $\mathbf{A}_{\ell} \hat{\mathbf{q}}_{\ell} = \hat{\mathbf{r}}_{\ell}$.
\item Since
$\mathbf{S}_{\ell}^{-1}$ is precomputed, we obtain $\hat{\boldsymbol{\lambda}}_{\ell}$ from $\hat{\mathbf{q}}_{\ell}$ via
$\hat{\boldsymbol{\lambda}}_{\ell}
= \mathbf{S}_{\ell}^{-1} \mathbf{C}_\ell \hat{\mathbf{q}}_{\ell}.$
\item Using the precomputed matrix $\mathbf{Y}_{\ell}= \mathbf{A}^{\hspace{-2pt}-1}_\ell \hspace{2pt} \mathbf{C}_\ell^{\top}$ we obtain
$\hat{\mathbf{w}}_{\ell} = \hat{\mathbf{q}}_{\ell} - \mathbf{Y}_{\ell} \hat{\boldsymbol{\lambda}}_{\ell}.$
\end{enumerate}
Observe that this procedure only involves one single (low dimensional) system of equations to solve.
\subsubsection{Assembly and solution of the global problem}
The procedure is basically analogous to the case of a homogenous boundary condition.
The global corrector matrix $\mathbf{Q}_h$ is assembled identically as before. The same holds for the interpolation matrix $\mathbf{P}_h$. With that, we obtain the LOD stiffness matrix $\mathbf{A}_H^{\LOD} \in \R^{N_H \times N_H}$ (associated with the left hand side of (\ref{lod-problem-eq-mixed-bc})) by
$$\mathbf{A}_H^{\LOD} = \mathbf{B}^H (\mathbf{P}_h+\mathbf{Q}_h) \mathbf{A}_h (\mathbf{P}_h+\mathbf{Q}_h)^{\top} \mathbf{B}^H.$$
In order to assemble the LOD load vector $\mathbf{f}_H \in \R^{N_H}$, we first need
to assemble the vector $\hat{\mathbf{f}}_h \in \R^{N_h}$ that stores the information gained from source correctors. For a given fine basis function $\phi_j$ the corresponding entry of $\hat{\mathbf{f}}_h$ is given by
\begin{align*}
\hat{\mathbf{f}}_h[j] &:= \int_{\Omega} \kappa \nabla \QFh \cdot \nabla \phi_j = \sum_{\ell = 0}^{N_H-1}
\int_{\mathcal{U}_{\ell}} \kappa \nabla \QFh^{K_{\ell}} \cdot \nabla \phi_j\\
&= \sum_{\ell = 0}^{N_H-1}
 \sum_{i = 0}^{N_{\ell,h}-1}
\hat{\mathbf{w}}_{\ell}[i] \int_{\mathcal{U}_{\ell}} \kappa \nabla \phi_{\ell,i} \cdot \nabla \phi_j.
\end{align*}
Consequently we obtain $\hat{\mathbf{f}}_h$ by matrix multiplication and summation as
\begin{align*}
\hat{\mathbf{f}}_h = - \mathbf{A}_h \hat{\mathbf{w}}_h,
%= - \sum_{\ell = 0}^{N_H-1} {\mathbf{R}_{\ell}^h}^{\!\top} \mathbf{A}_{\ell} \hat{\mathbf{w}}_{\ell}.
    \qquad
    \mbox{where} \enspace
\hat{\mathbf{w}}_h := - \sum_{\ell =
    0}^{N_H-1} {\mathbf{R}_{\ell}^h}^{\!\top} \hat{\mathbf{w}}_{\ell}.
\end{align*}
The standard load vector $\mathbf{f}_h \in \R^{N_h}$ associated with a classical fine element method on the fine grid $\T_h$ is given by
$$\mathbf{f}_h[i]:= (\mathcal{F}+\mathcal{F}^{\text{s}})(\phi_i).$$
In total, we obtain the LOD load vector as $\mathbf{f}_H = \mathbf{B}^H (\mathbf{P}_h + \mathbf{Q}_h ) (\mathbf{f}_h + \hat{\mathbf{f}}_h) \in \R^{N_H}$ (i.e. the vector associated with the right hand side of (\ref{lod-problem-eq-mixed-bc})).
%by \begin{align*}\mathbf{f}_H := \mathbf{B}^H (\mathbf{P}_h + \mathbf{Q}_h ) (\mathbf{f}_h + \hat{\mathbf{f}}_h).\end{align*}
Using this, we can solve for ${\mathbf{U}}_H^{\LOD} \in \R^{N_H}$ with
$$\mathbf{A}_H^{\LOD} {\mathbf{U}}_H^{\LOD} =
\mathbf{B}^H (\mathbf{P}_h + \mathbf{Q}_h ) (\mathbf{f}_h - \mathbf{A}_h \hat{\mathbf{w}}_h)$$
%\mathbf{f}_H$$
and obtain the final solution vector ${\mathbf{u}}_h^{\LOD} \in \R^{N_h}$ of our LOD approximation by
$${\mathbf{u}}_h^{\LOD} := (\mathbf{P}_h+\mathbf{Q}_h)^{\top} {\mathbf{u}}_H^{\LOD}
+ \hat{\mathbf{w}}_h,$$
%%- \sum_{\ell = 0}^{N_H-1} {\mathbf{R}_{\ell}^h}^{\!\top} \hat{\mathbf{w}}_{\ell},$$
respectively $u_{\LOD}:= \sum_{i=0}^{N_h-1} {\mathbf{u}}_h^{\LOD}[i] \phi_i$.
\iffalse
\todo[inline]{
Alternative formulation:
$$\mathbf{A}_H^{\LOD} {\mathbf{U}}_H^{\LOD} = \mathbf{B}^H (\mathbf{P}_h + \mathbf{Q}_h ) (\mathbf{f}_h + \mathbf{A}_h \hat{\mathbf{w}}_h)$$
$${\mathbf{u}}_h^{\LOD} := (\mathbf{P}_h+\mathbf{Q}_h)^{\top} {\mathbf{u}}_H^{\LOD} - \hat{\mathbf{w}}_h$$
furthermore I'd rather change the definition of $\hat{\mathbf{w}}_h$
to have opposite sign, the it would read:
$$\mathbf{A}_H^{\LOD} {\mathbf{U}}_H^{\LOD} = \mathbf{B}^H (\mathbf{P}_h + \mathbf{Q}_h ) (\mathbf{f}_h - \mathbf{A}_h \hat{\mathbf{w}}_h)$$
$${\mathbf{u}}_h^{\LOD} := (\mathbf{P}_h+\mathbf{Q}_h)^{\top} {\mathbf{u}}_H^{\LOD} + \hat{\mathbf{w}}_h$$
this is then consistent with the usual formulation for inhomogeneous
dirichlet boundary conditions, where we subtract the known BC matching function.
}
\fi
With these changes, Algorithm {\bf{1}} and {\bf{2}} can be modified in an obvious way.
\section{Linear elliptic eigenvalue problems}
\label{section:LOD-for-EV}
In this section we describe how the LOD can be applied to solve eigenvalue problems. Eigenvalue problems have a distinct status since the LOD is particularly efficient for tackling them, even if the diffusion coefficient has no multiscale character.
We consider the following linear eigenvalue problem with a homogenous Dirichlet boundary condition, i.e. $\Gamma_D=\partial \Omega$, $g=0$ and $\Gamma_N=\emptyset$. We seek tuples $(\lambda^{(n)},u^{(n)})\in \R \times H^1_0(\Omega)$ with
\begin{equation}\label{e:modelproblem:ev}
  \mathcal{A}(  u^{(n)},v )
  %=\int_{\Omega}\kappa \nabla u^{(n)}\cdot \nabla v
 =\lambda^{(n)} \int_{\Omega} u^{(n)} v \quad\text{for all } v\in H^1_0(\Omega).
\end{equation}
We assume that the eigenvalues are ordered and that the eigenvalues are $L^2$-normalized, i.e. we have $\lambda^{(n)}\le\lambda^{(n+1)}$ and $\| u^{(n)} \|_{L^2(\Omega)}=1$ for all $n\in\mathbb{N}$. Note that we always have $\lambda^{(0)}>0$.
The corresponding fine-scale reference solutions are given by the tuples $(\lambda_h^{(n)},u_h^{(n)})\in \R \times \VhGD$ with $\| u_h^{(n)} \|_{L^2(\Omega)}=1$ and
\begin{equation}\label{e:modelproblem:ev:ref}
  \mathcal{A}(  u^{(n)}_h,v_h )
  %=\int_{\Omega}\kappa \nabla u^{(n)}\cdot \nabla v
 =\lambda^{(n)}_h \int_{\Omega} u^{(n)}_h v_h \quad\text{for all } v_h \in \VhGD
\end{equation}
and where $0<\lambda^{(n)}_h\le\lambda^{(n+1)}_h$.
\subsection{LOD approximation of a linear eigenvalue problem}
 The LOD method for approximating the eigenpair $(\lambda^{(n)},u^{(n)})$ of (\ref{e:modelproblem:ev}) can be implemented in a straightforward way. After solving all local problems as described in Section~\ref{subsection-assembling-solving-local-problems} we can assemble the LOD system matrix $\mathbf{A}_H^{\LOD}$ as described in Section~\ref{subsection-assembling-solving-global-problem}. This has to be done only once and can be seen as a one-time preprocessing step. The local solutions and the LOD system matrix are stored and can be reused for all eigenvalues that we are interested in. We define the LOD approximation as follows.
\begin{definition}[Eigenvalue LOD approximation]
  \label{definition-lod-approx-ev-prob}
Let $k\in \mathbb{N}$ be fixed and let
$\Qh: \VHGD \rightarrow \VhGD$ denote the
corresponding correction operator as in Definition \ref{def-loc-lod}.
%We denote $(\Id+\Qh)(\Phi_H):=\Phi_H+\Qh(\Phi_H)$ for $\Phi_H \in \VHGD$.
We seek tuples $(\lambda^{(n)}_H,u^{(n)}_H)\in \R \times \VHGD$ with $\| u^{(n)}_H +\Qh( u^{(n)}_H ) \|_{L^2(\Omega)}=1$ such that for all $\Phi_H \in \VHGD$
  \begin{align}
    \label{pg-lod-problem-eq-ev-prob}
    \mathcal{A} ( u^{(n)}_H +\Qh( u^{(n)}_H ) , \Phi_H +\Qh( \Phi_H) ) &=&
    \lambda^{(n)}_H \int_{\Omega} ( u^{(n)}_H +\Qh( u^{(n)}_H ) ) \hspace{3pt} ( \Phi_H +\Qh( \Phi_H) )
  \end{align}
 and $0<\lambda^{(n)}_H \le \lambda^{(n+1)}_H$ for $0 \le n < N_H$.
We denote the arising LOD approximations of the eigenvectors by
$u_{\LOD}^{(n)}:=(\Id+\Qh)(u_H^{(n)})$.
\end{definition}
The eigenvalue problem (\ref{pg-lod-problem-eq-ev-prob}) can be solved with any favorite solver. The main cost (for computing $\mathbf{A}_H^{\LOD}$ and $\mathbf{Q}_h$) arise only once. Hence the computational advantage of the method becomes bigger the more eigenvalues we want to compute. Furthermore, the obtained convergence rates in $H$ surpass the classical convergence rates. The following result was proved in \cite{MaP15}.
\begin{theorem}
\label{theorem-lod-eigenvalue-error}
Recall the reference problem (\ref{e:modelproblem:ev:ref}). By $C$ we denote a generic constant as in Theorem \ref{definition-lod-approx-ev-prob}. Let $U_k(K)=\Omega$ for all $K \in \T_H$ and assume that $H \le \left( \sqrt{n} \lambda_h^{(n)} \right)^{-1/2}$. Then it holds
\begin{align}
\label{lod-eigenvalue-error}\frac{\lambda_H^{(n)}-\lambda_h^{(n)}}{\lambda_h^{(n)}} \le C n (\lambda_h^{(n)})^2 H^4
\end{align}
for all $0\le n <N_H$.
%
%Let the multiplicity of the eigenvalue $\lambda_h^{(n)}$ be denoted by $m$, i.e. we have $\lambda_h^{(n)}=\ldots=\lambda_h^{(n+m-1)}$. The corresponding eigenspace is given by the orthonormal basis $\{ u_h^{(n)},\ldots,u_h^{(n+m-1)}\}$. Let $u_{\LOD}^{(n)},\ldots,u_{\LOD}^{(n+m-1)}$ denote the Rayleigh-Ritz LOD approximations of \eqref{pg-lod-problem-eq-ev-prob}. Furthermore denote $\rho_{\lambda}:=\max_{j \not\in \{ n,n+1,\ldots,n+m-1\}} \frac{\lambda_h^{(n)}}{|\lambda_h^{(n)}-\lambda_H^{(j)}|}$. If $H$ is such that $H \le C n^{-1/3} (1+\rho_{\lambda})^{-1/3} (\lambda_h^{(n)})^{-1/2}$ then there exists an orthonormal basis $\{ \tilde{u}_h^{(n)},\ldots,\tilde{u}_h^{(n+m-1)}\}$ of $\{ u_h^{(n)},\ldots,u_h^{(n+m-1)}\}$ such that for all $j\in \{ 0,\ldots, m-1\}$
%\begin{align*}
%\| \tilde{u}_h^{(n+j)} - u_{\LOD}^{(n+j)} \|_{H^1(\Omega)} &\le C \left( \sqrt{n} (\lambda_h^{(n)})^{3/2} + H n(1 + \rho_{\lambda})  (\lambda_h^{(n)})^{2} \right) H^2,\\
%\| \tilde{u}_h^{(n+j)} - u_{\LOD}^{(n+j)} \|_{L^2(\Omega)} &\le C \left( n(1 + \rho_{\lambda})  (\lambda_h^{(n)})^{3/2} \right) H^3.
%\end{align*}
\end{theorem}
Besides the fourth order convergence rates for the eigenvalues, it can be also shown that the corresponding eigenfunctions converge with higher rates (cf. \cite{MaP15}). Concerning the $H^1$-error, the convergence is of quadratic order (in $H$) and concerning the $L^2$-error of cubic order.
These rates are higher than the rates that we obtained for the LOD for standard linear elliptic problems. Observe that these high convergence rates allow for much coarser grids and hence for a reduced computational complexity.
%Basically, Theorem \ref{theorem-lod-eigenvalue-error} predicts a quadratic order convergence (in $H$) for the $H^1$-error and even a cubic order convergence for the $L^2$-error. These rates are higher than the rates that we obtained for the LOD for standard linear elliptic problems. Observe that these high convergence rates allow for much coarser grids and hence for a reduced computational complexity.
\begin{remark}[Truncation]
If the localization parameter $k$ is chosen such that $k \ge m |\log(H)|$, the truncation error will be of order $\mathcal{O}(H^{rm})$ (with $r$ as in Theorem \ref{theorem-apriori-linear}). Consequently, for properly chosen $m$, the convergence rates in Theorem \ref{theorem-lod-eigenvalue-error} remain valid even for the localized method.
\end{remark}
\subsection{Two-grid post-processing}
\label{sec:two-grid-postprocessing}
The high convergence rates depicted in Theorem \ref{theorem-lod-eigenvalue-error} can be even improved by applying a two-grid post-processing technique as initially suggested in \cite{XuZh01}. The post-processing technique can be applied, if it is affordable to solve global (linear elliptic, non-eigenvalue) problems in the full fine space $\VhGD$. We define the post-processed LOD approximation as follows.
\begin{definition}[Post-processed eigenvalue LOD approximation]
\label{definition-lod-approx-ev-prob-postprocess}
Let $(\lambda^{(n)}_H,u_{\LOD}^{(n)})\in \R \times \VhGD $ (with $0 \le n < N_H$) denote the eigenpair approximations obtained with the LOD as stated in Definition \ref{definition-lod-approx-ev-prob}. We call $(\lambda^{(n)}_{H,\mbox{\tiny{post}}},u_{\mbox{\tiny{LOD,post}}}^{(n)})\in \R \times \VhGD$ the corresponding post-processed approximations if $u_{\mbox{\tiny{LOD,post}}} \in \VhGD$ solves
\begin{align*}
\mathcal{A} ( u_{\mbox{\tiny{LOD,post}}}^{(n)} ,  v_h ) &=
    \lambda^{(n)}_H \int_{\Omega} u_{\mbox{\tiny{LOD}}}^{(n)} \hspace{2pt} v_h \qquad \mbox{for all } v_h \in \VhGD
\end{align*}
and where we define
$$
\lambda^{(n)}_{H,\mbox{\tiny{post}}} := \frac{\mathcal{A} ( u_{\mbox{\tiny{LOD,post}}}^{(n)} ,  u_{\mbox{\tiny{LOD,post}}}^{(n)} )}{\| u_{\mbox{\tiny{LOD,post}}}^{(n)} \|_{L^2(\Omega)}^2}.
$$
\end{definition}
Observe that the post-processing step involves to solve an additional linear elliptic problem in the full fine space $\VhGD$. Before discussing the feasibility of this step, we present a corresponding a priori error estimate.
\begin{theorem}
\label{theorem-on-post-processed-LOD-appr}
Let $(\lambda^{(n)}_H,u_{\LOD}^{(n)})\in \R \times \VhGD $ (with $0 \le n < N_H$) denote the LOD approximations as in Definition \ref{definition-lod-approx-ev-prob} and let $(\lambda^{(n)}_{H,\mbox{\tiny{post}}},u_{\mbox{\tiny{LOD,post}}}^{(n)})\in \R \times \VhGD$ denote the post-processed LOD approximations as in Definition \ref{definition-lod-approx-ev-prob-postprocess}. Then it holds
\begin{align*}
| \lambda^{(n)}_{H,\mbox{\tiny{post}}} - \lambda_h^{(n)} | &\le C (\lambda^{(n)}_{H} -  \lambda_h^{(n)})^2
+ C (\lambda_h^{(n)})^2 \| u_{\LOD}^{(n)} - u_h^{(n)}\|_{L^2(\Omega)}^2 \qquad \mbox{and}\\
\| u_{\mbox{\tiny{LOD,post}}}^{(n)} - u_h^{(n)} \|_{H^1(\Omega)} &\le C |\lambda^{(n)}_{H} -  \lambda_h^{(n)}|
+ C \lambda_h^{(n)} \| u_{\LOD}^{(n)} - u_h^{(n)}\|_{L^2(\Omega)}.
\end{align*}
Consequently, we obtain that the eigenvalue $\lambda^{(n)}_{H,\mbox{\tiny{post}}}$ converges at least with order $\mathcal{O}(H^8)$ to $\lambda_h^{(n)}$ and that the $H^1$-error between $u_{\mbox{\tiny{LOD,post}}}^{(n)}$ and $u_h^{(n)}$ converges to zero at least with order $\mathcal{O}(H^4)$.
\end{theorem}
In view of Theorem \ref{theorem-on-post-processed-LOD-appr} we can see that the LOD can be a powerful tool to tackle linear eigenvalue problems even if $\kappa$ has no multiscale character. If the cost for solving a full linear system on the fine scale is still feasible, it can be highly efficient to apply the LOD with pre- and post-processing. The extremely high convergence rates in $H$ (at least $\mathcal{O}(H^8)$ if we are interested in the eigenvalues) allow to choose a very coarse grid $\T_H$. Depending on how coarse we choose $\T_H$, the truncation might be even skipped completely. First, we solve the corrector problems in a preprocessing step. This involves a number of linear elliptic fine scale problems that can be solved in parallel. After the correctors are available, we assemble the global LOD stiffness matrix $\mathbf{A}_H^{\LOD}$, which has a very low dimension. All operations of the chosen algebraic eigenvalue solver, only involve $\mathbf{A}_H^{\LOD}$ and can be hence performed quickly and sequential. Once the eigenpairs are computed, we apply the postprocessing step, which only involves one global fine scale problem per eigenpair. Again, we can do this in parallel. With this approach, we can decrease the computational complexity.
The implementation of pre- and post-processing is obvious, following the ideas presented in Section~\ref{section-lod-for-elliptic-problems}.
\section{Numerical Examples}
The algorithms present in this paper are available as a prototype
implementation \footnote{\url{https://gitlab.dune-project.org/christi/dune-py-lod/}}.
In this section we use this
implementation to solve an elliptic boundary value problem and an eigenvalue problem.
As the proposed algorithmic approach to the LOD is purely algebraic,
the prototype is implemented in \texttt{python} 2.7 using the modules
\texttt{numpy} \cite{walt2011numpy} and \texttt{scipy}
\cite{millman2011scipy}. The numerical example is implemented in
\texttt{src/eigenvalues.py} and the actual LOD algorithm is
implemented in the \texttt{compute\_correction} function of the
\texttt{lod} module.
As input it requires different matrices, most of them directly related
to the fine scale model:
\begin{itemize}[noitemsep]
\item \texttt{Adc}: Sparse matrix with element wise contributions to
  the stiffness matrix, i.e. each diagonal block contains a particular
  local stiffness matrix $\mathbf{A}_t$, see definition
  \ref{element-matrices}. Basically this is the conforming fine scale
  operator assembled with discontinuous Galerkin shape functions.
\item \texttt{fdc}: Vector with element wise contributions to the
  right-hand-side. Similarly defined as $\mathbf{A}_t$, it is
  necessary to compute the source corrector. In this implementation we
  do not distinguish between different right-hand-side contributions,
  but just compute a full corrector. Note, that one may safe a bit,
  if the correct only needs to be computed for the boundary and not
  for the source term.
\item \texttt{BH}: Coarse-Mesh boundary correction $\mathbf{B}^H$, see
  definition \ref{def:bc-corr}. This matrix has
  exactly one 1 in each row, every other entry is 0. It removes all
  constraint unknowns.
\item \texttt{Mh}: Fine-mesh Mass-matrix, see definition
  \ref{def:global-matrices}. In many cases it can be convinient to
  also compute the mass matrix element wise, see
  \ref{element-matrices}, and then compute the global matrix as given
  in \eqref{global-mass-and-stiffness-matrix}.
\item \texttt{Ph}:
  Projection matrix $\mathbf{P}_h$ from the
  coarse-mesh Lagrange space to the fine-mesh Lagrange space,
  see equation~\eqref{projection-matrix}.
\item \texttt{Sigma}: Index-mapping from local to global dofs of the
  fine scale discretization, a block matrix containing all
  $\bigsigma_{\!t}$, see definition \ref{local-to-global-mapping}.
\item \texttt{SubInfo} provides a list with different restriction
  operators for each
  sub-domain (actually our implementation provides this information via
  a generator). Each entry has to provide the following details:
  \begin{itemize}[nosep]
  \item \texttt{SubInfo.R} Restriction operator mapping from the fine
    space restriction to the patch, see equation~\eqref{restriction-op}.
  \item \texttt{SubInfo.RH} Restriction operator mapping from the
    coarse space restriction to the patch, see equation~\eqref{restriction-op-coarse}.
  \item \texttt{SubInfo.TH} Restriction operator from the coarse
    space to a coarse cell, see equation~\eqref{restriction-op-coarse-element}.
  \end{itemize}
\end{itemize}
In our examples this fine scale model is
implemented using the DUNE \cite{dune08:1,dune08:2} framework, a modern \texttt{C++} library for
grid based methods. It requires the DUNE core modules in version 2.4
\cite{blatt2016distributed} and DUNE-PDELab \cite{bastian2010generic} in the 2.4
compatible version.
For the sub-domain information we provide some
additional infrastructure, so that it can be computed in
\texttt{python}, eventhough it will usually be faster to assemble this
information also in the framework.

The examples presented in the following are run on an AMD
Epic 7501 server. Due to the global interpreter lock, python does not
immediately allow for thread parallelization. Surely there are a
range of different approaches for python to work around the
limitations in multi-threading, but this would exceed the scope of an
illustrative prototype implementations. We therefor decided to enforce
all computations to be run sequentially on a single core. In
particular for the eigenvalue solves, ARPACK's OpenMP parallelization
was limited to a single thread to make timings comparable.

Furthermore, all our numerical experiments are performed for bilinear finite elements (cf. \ref{notation-quadrilaterals}) on uniform quadrilateral meshes with square elements.
\subsection{Elliptic Problem}
\label{sec:ell_problem}
\begin{figure}
  \centering
  \includegraphics[height=0.325\linewidth]{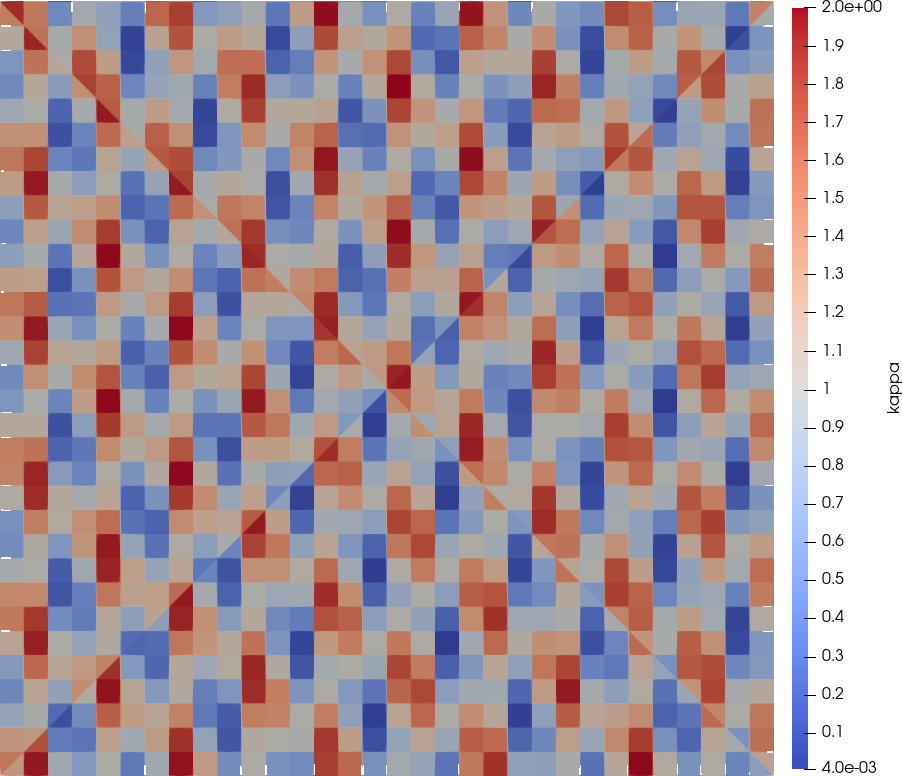}\qquad\qquad
  \includegraphics[height=0.325\linewidth]{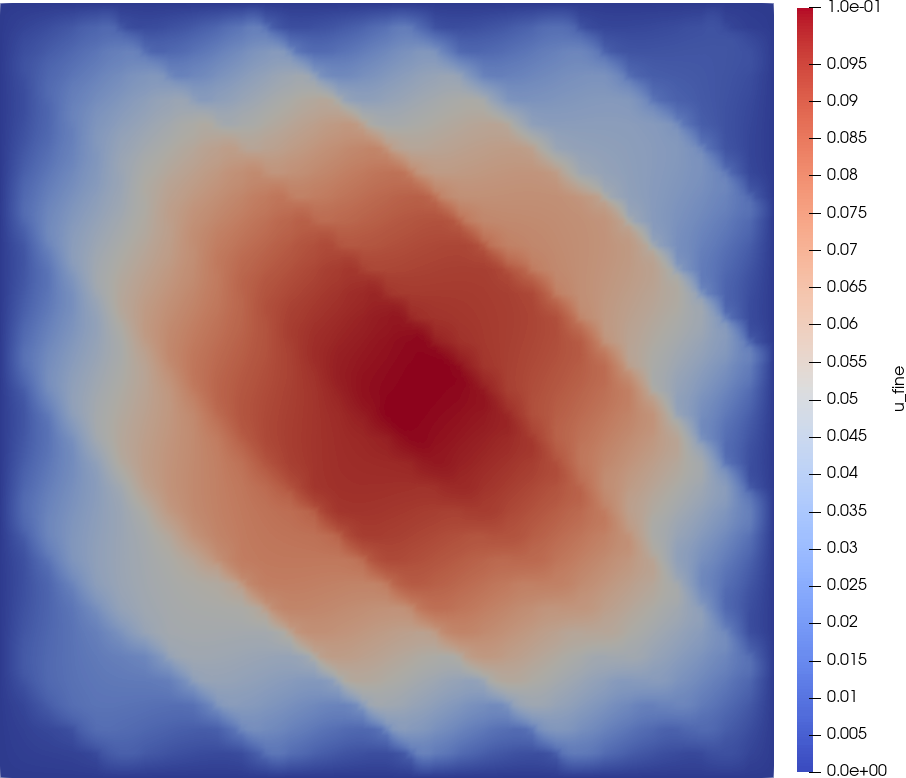}
  \caption{Model Problem \ref{sec:ell_problem}. Left: Multiscale coefficient $\kappa$ for $\varepsilon=2^{-5}$. Right: Reference solution $u$ obtained on a fine grid with resolution $h=2^{-7}$. As a remark, the LOD solution for $H=2^{-3}$ is visually indistinguishable from this reference solution.}
   \label{fig:elliptic_kappa_solution}
\end{figure}
\begin{figure}
  \includegraphics[height=0.325\linewidth]{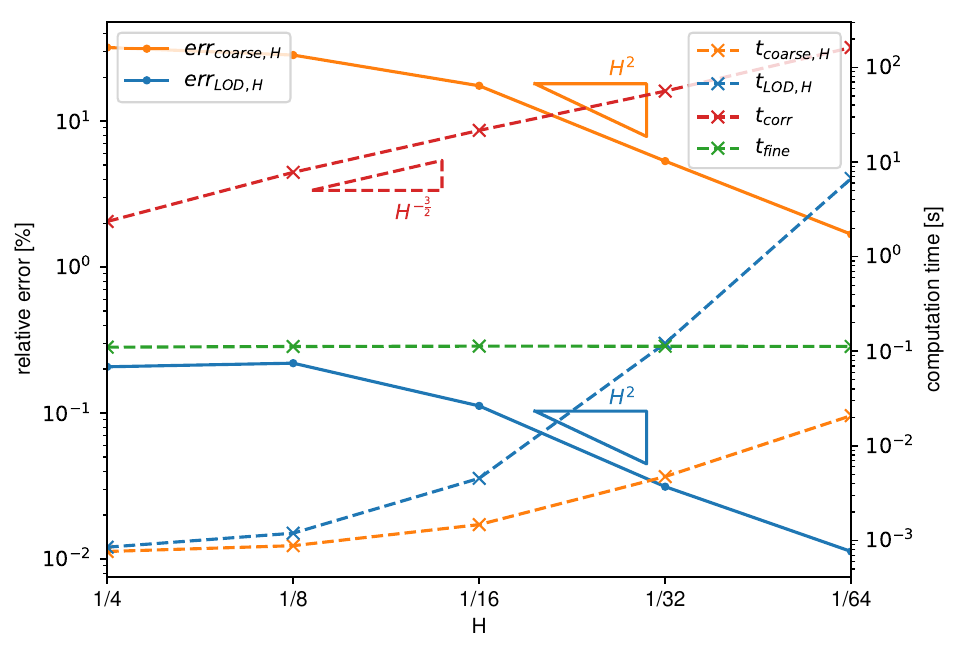}\qquad
  \includegraphics[height=0.325\linewidth]{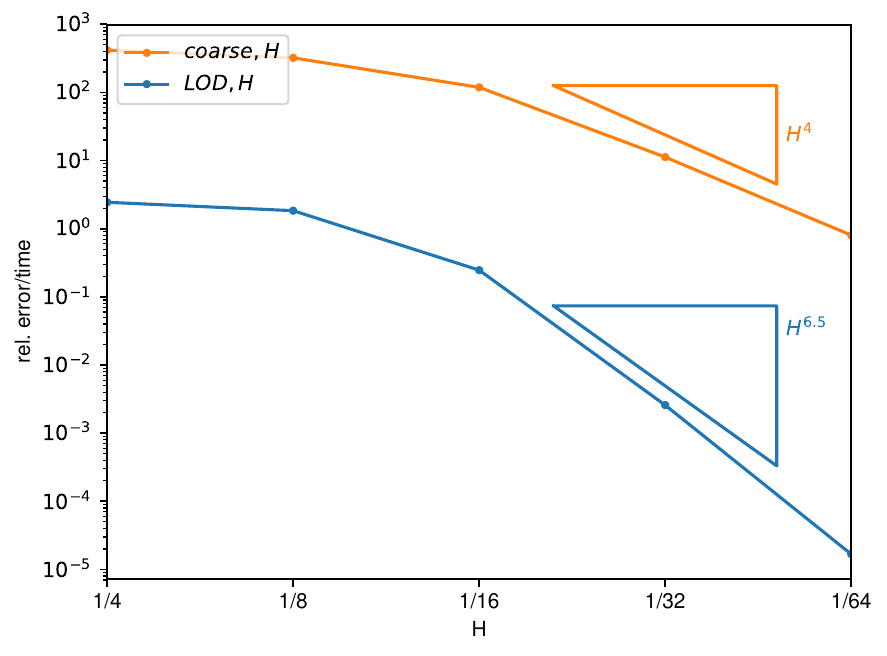}\\
  \caption{Model Problem \ref{sec:ell_problem}. Left: $L^2$-errors for different coarse resolutions $H$. We compare the error $\operatorname{err}_{\coarse,H}$ for the standard approach with the LOD-error $\operatorname{err}_{\LOD,H}$ on the same mesh with overlap parameter $k \approx - 0.9 \log_2(H)$. In addition, we list the times for solving the coarse problem $t_{\coarse,H}$, the macroscopic LOD problem $t_{\LOD,H}$, computing the correctors $t_{\corr}$  and for solving the full fine grid problem $t_{\fine}$. Right: Performance comparison for solving the standard coarse problem compared to the macroscopic LOD-problem.}
   \label{fig:elliptic_performance}
\end{figure}
We first want to exemplify the algorithm and the basic
performance using a standard elliptic model problem on the unit square
$\Omega=(0,1)^2$, with a heterogenous coefficient $\kappa$ and a constant
source term $f = 1$:
\begin{equation*}\label{eq:model}
  \begin{aligned}
    -\nabla \cdot \kappa \nabla u &= 1\quad \text{ in }\Omega, \\
    u & = 0 \quad \text{ on }\partial\Omega. \\
  \end{aligned}
\end{equation*}
The diffusion coefficient 
\begin{equation*}
\kappa(x) = 1 + 10^{-8}
  + \frac{1}{2} \sin( \left\lfloor x_1 + x_2 \right\rfloor + \left\lfloor \frac{x_1}\varepsilon \right\rfloor + \left\lfloor \frac{x_2}\varepsilon \right\rfloor)
  + \frac{1}{2} \cos( \left\lfloor x_2 - x_1 \right\rfloor + \left\lfloor \frac{x_1}\varepsilon \right\rfloor + \left\lfloor \frac{x_2}\varepsilon \right\rfloor)
 \end{equation*} 
   exhibits a multiscale structure.
  The fine scale reference solution is computed on level 7 (of quadrilateral mesh), i.e. $h = 2^{-7}$ and $\varepsilon$ is chosen to be $\varepsilon = 2^{-5} = 4h$. Figure~\ref{fig:elliptic_kappa_solution} shows the coefficient and the aforementioned fine scale solution used as a reference for the computation of errors. The same fine mesh is also used in all LOD computations. As for the LOD algorithm, we follow the descriptions presented in Section~\ref{section5:LOD-including-source-and-boundary-terms}, which means that we included a source corrector $\mathcal{Q}_{1,h}$ for constant inputs in order to improve the accuracy. For optimal convergence orders, the overlap parameter $k$ needs to be proportional to $\log(H)$, which is why we chose $k \approx - 0.9 \log_2(H)$ in this numerical experiment.
 Timings for the different phases are presented in Figure~\ref{fig:elliptic_performance}. It should be noted that the implementation is to illustrate some central properties of the LOD in terms of performance, which is why it is not tuned for efficiency. In practice there are a couple of possible improvements
as we will discuss later. 

In Figure~\ref{fig:elliptic_performance}, $t_{\coarse,H}$ denotes the time for solving the elliptic test problem using standard $Q_1$-FEM on the coarse mesh $\T_H$; $t_{\LOD,H}$ denotes the time for solving the macroscopic LOD problem \eqref{lod-problem-eq-mixed-bc} on $\T_H$; $t_{\corr}$ denotes the time for computing all the local correctors and for assembling the LOD-stiffness matrix associated with \eqref{lod-problem-eq-mixed-bc}; finally, $t_{\fine}$ denotes the time for solving the test problem on the full fine mesh $\T_h$ with resolution $h=2^{-7}$. The relative errors are computed as
 \begin{align*}
  \operatorname{err}_{\coarse,H} = \frac{\| u - u_H \|_{L^2(\Omega)}}{\| u \|_{L^2(\Omega)}}
\qquad 
\mbox{and}\qquad 
\operatorname{err}_{\LOD,H} = \frac{\| u - u_{\LOD} \|_{L^2(\Omega)}}{\| u \|_{L^2(\Omega)}},
\end{align*}
where $u$ is the fine scale reference solution, $u_H$ the standard finite element approximation on the coarse mesh and $u_{\LOD}$ the LOD approximation for a given coarse mesh $\T_H$, fine mesh $\T_h$ and overlap parameter $k \approx - 0.9 \log_2(H)$. From the graphs in Figure~\ref{fig:elliptic_performance} we can see that the main cost of the LOD account for the assembly of the correctors which roughly grows (in a sequential implementation) with the rate $H^{-3/2}$, whereas the cost for solving the global LOD system are negligible in a one-shot simulation. We can also clearly see that the method becomes inefficient if the fine mesh $\T_h$ is only slightly finer than the coarse mesh $\T_H$. In practical situations, the fine mesh $\T_h$ needs to be significantly finer than $\T_H$, where the regime $h\simeq H^2$ is often reasonable. We can also see that the costs for solving the global LOD system grow faster compared to the costs for solving the standard coarse system of the same dimension. The reasons is that the number of nonzeros in the system matrix grows by the factor $k^d$, where $k$ is the localization parameter. This pollution of the sparsity structure causes an increased computational complexity. However, comparing these increased costs with the improvement of the accuracy in the right graph of Figure~\ref{fig:elliptic_performance}, we can see that the LOD is significantly more efficient if we just look at the macroscopic solve. This is an important observation since the one-time costs for computing the correctors become negligible if there is either a high capacity for parallelization or if the calculations have to be repeated for many source terms. In the next model problem we will face such a situation where the overhead caused by the computation of the correctors is compensated by its repeated usage and where the LOD outperforms the standard approach on the fine grid.

\color{black}
\begin{figure}
\centering
\begin{tabular}{cc@{\,\,}c@{\,\,}c@{\,\,}c}
  &
  \includegraphics[width=0.2\linewidth]{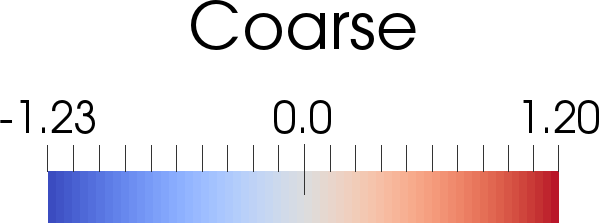} &
  \includegraphics[width=0.2\linewidth]{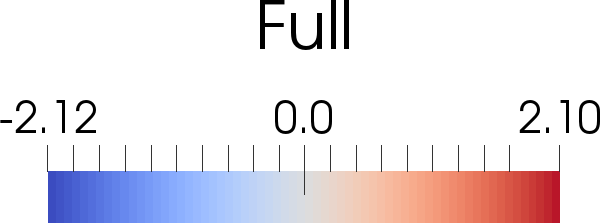} &
  \includegraphics[width=0.2\linewidth]{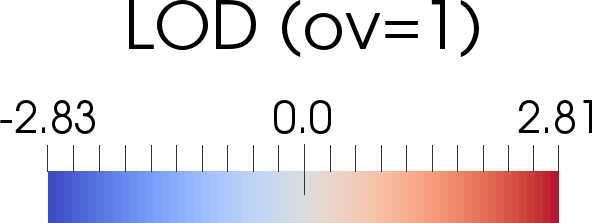} &
  \includegraphics[width=0.2\linewidth]{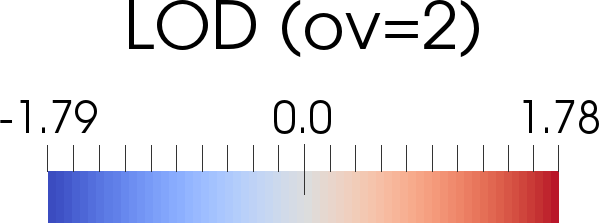} \\
  \raisebox{1.7cm}{\rotatebox{90}{\textsf{mode 1}}} &
  \includegraphics[width=0.2\linewidth]{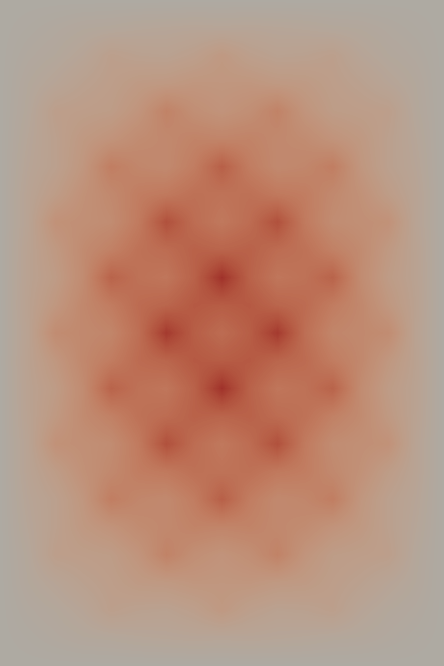} &
  \includegraphics[width=0.2\linewidth]{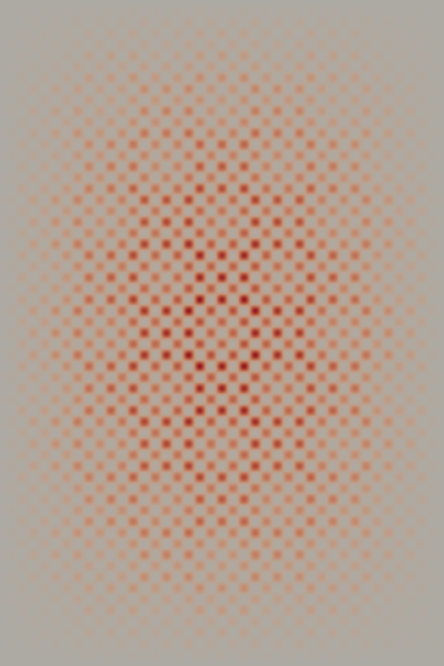} &
  \includegraphics[width=0.2\linewidth]{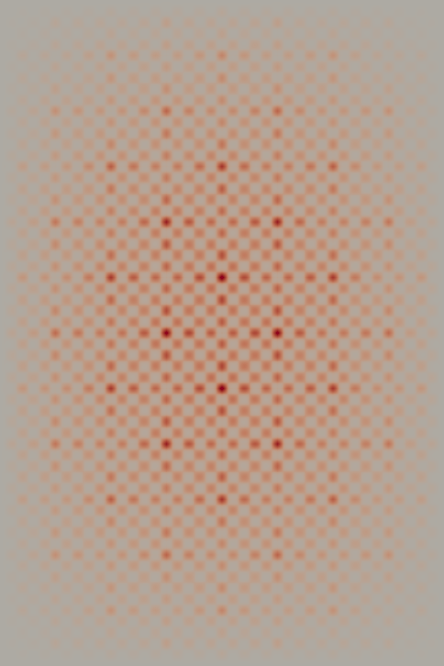} &
  \includegraphics[width=0.2\linewidth]{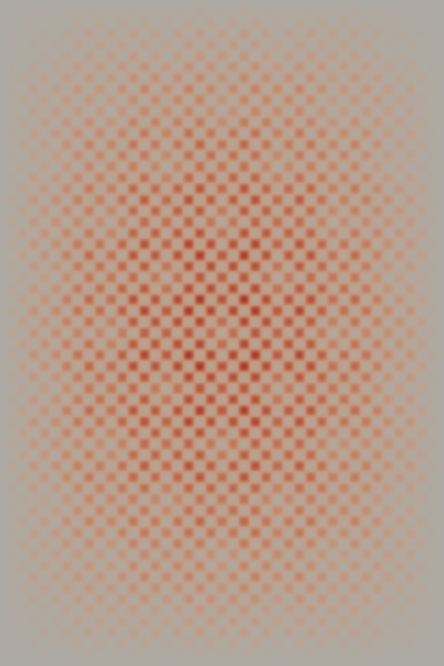} \\
  \raisebox{1.7cm}{\rotatebox{90}{\textsf{mode 6}}} &
  \includegraphics[width=0.2\linewidth]{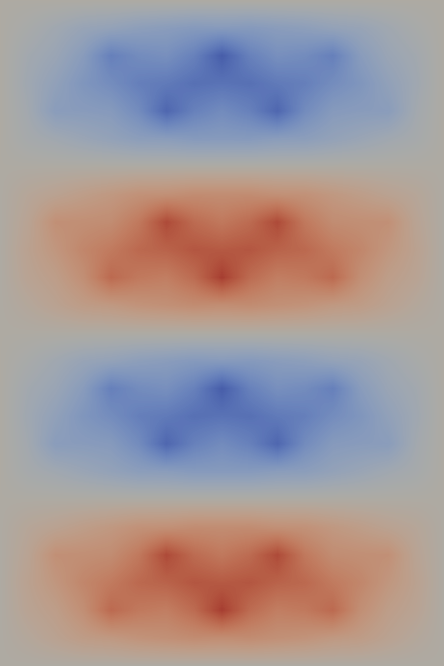} &
  \includegraphics[width=0.2\linewidth]{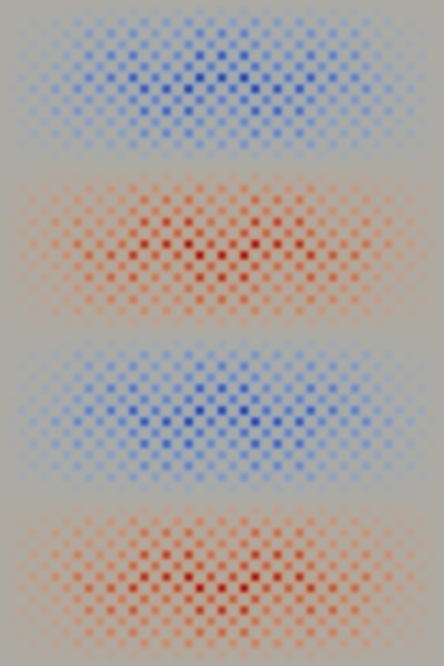} &
  \includegraphics[width=0.2\linewidth]{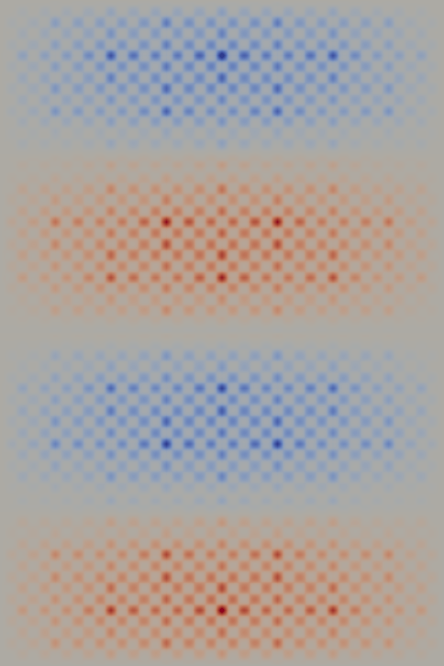} &
  \includegraphics[width=0.2\linewidth]{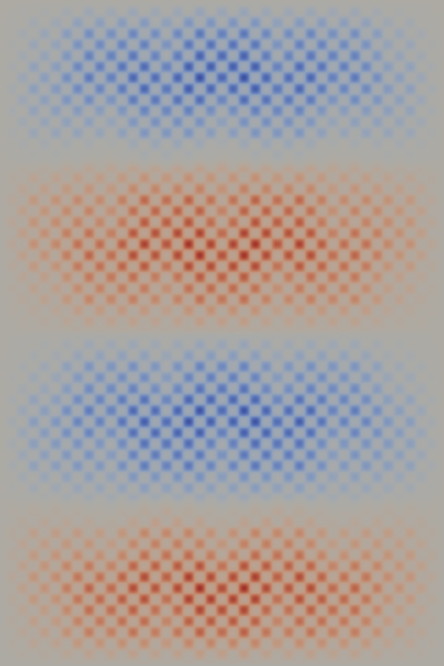} \\
  \raisebox{1.7cm}{\rotatebox{90}{\textsf{mode 11}}} &
  \includegraphics[width=0.2\linewidth]{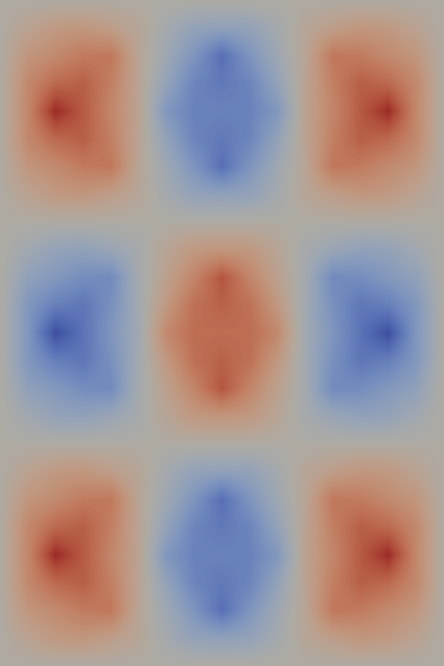} &
  \includegraphics[width=0.2\linewidth]{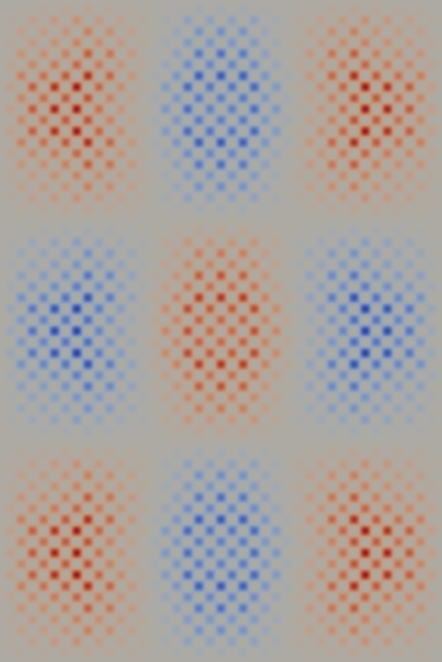} &
  \includegraphics[width=0.2\linewidth]{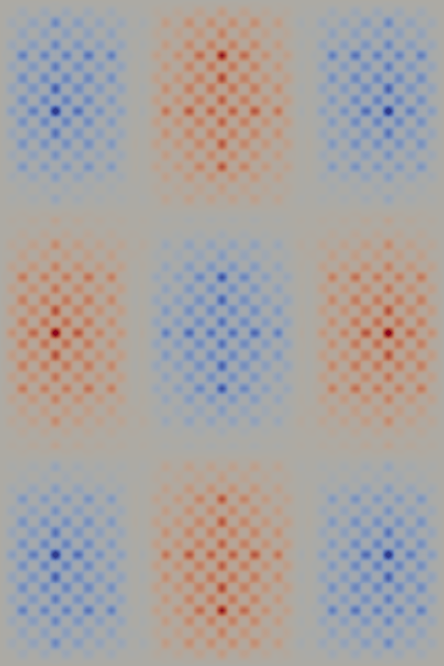} &
  \includegraphics[width=0.2\linewidth]{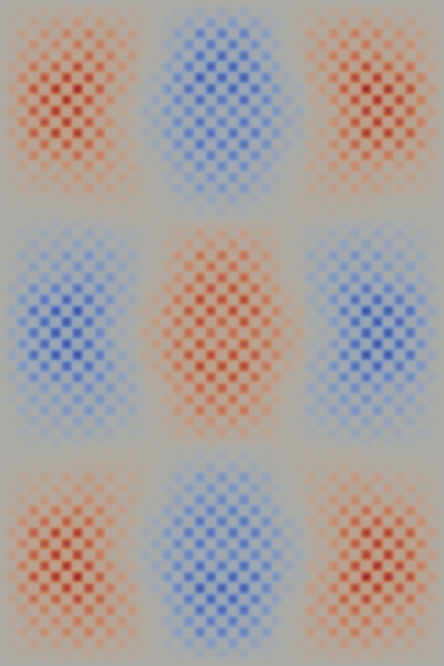} \\
  \raisebox{1.7cm}{\rotatebox{90}{\textsf{mode 20}}} &
  \includegraphics[width=0.2\linewidth]{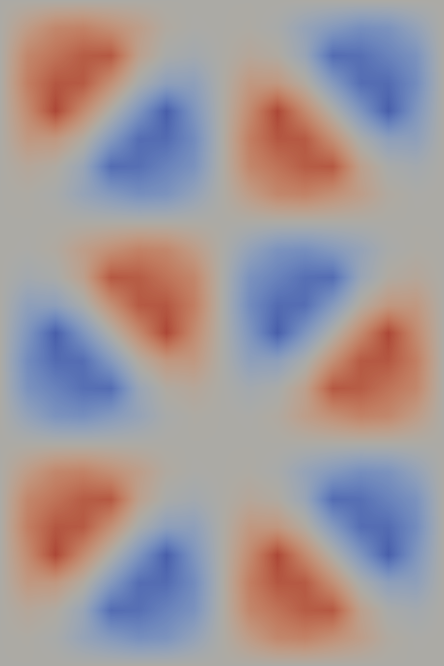} &
  \includegraphics[width=0.2\linewidth]{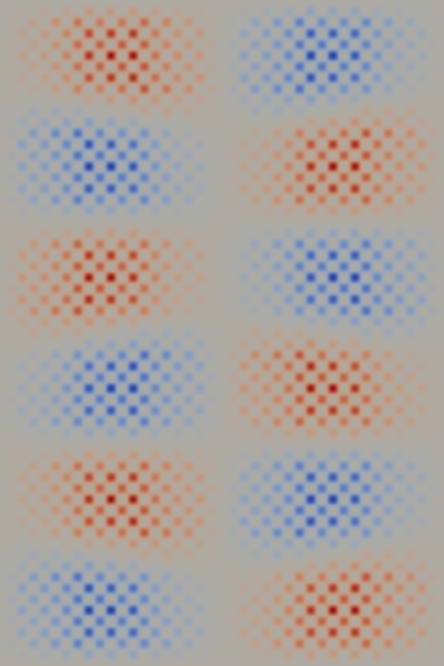} &
  \includegraphics[width=0.2\linewidth]{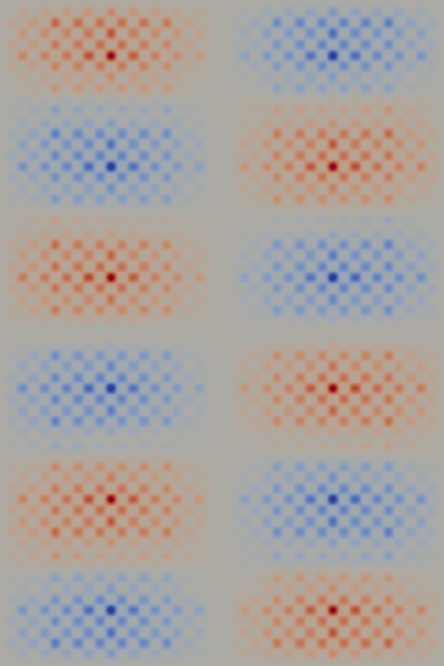} &
  \includegraphics[width=0.2\linewidth]{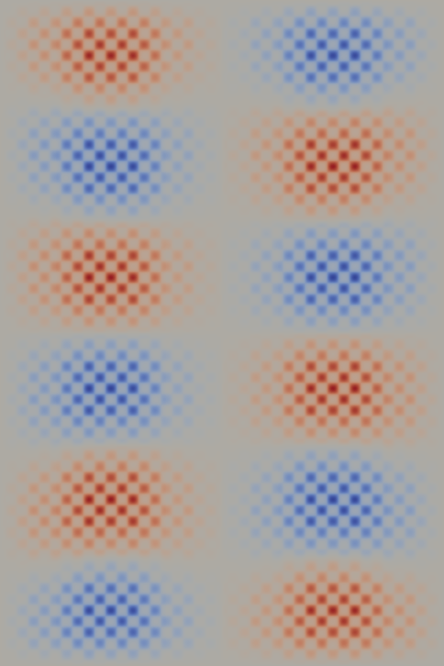}
\end{tabular}
  \caption{Model Problem \ref{sec:ev_problem} with $\nu=20$ and $\gamma=2\cdot10^4$. Comparison of the 1st, 6th, 11th, and 20th eigenmode computed for problem
    \eqref{eq:ev-problem} using the coarse ($H=2^{-3}$), the full
    ($h=2^{-8}$) and the LOD discretization ($H=2^{-3}$, overlap $H$
    and $2H$, i.e. $k=1$ and $k=2$) .}
   \label{comparison}
\end{figure}
\begin{figure}
  \includegraphics[height=0.325\linewidth]{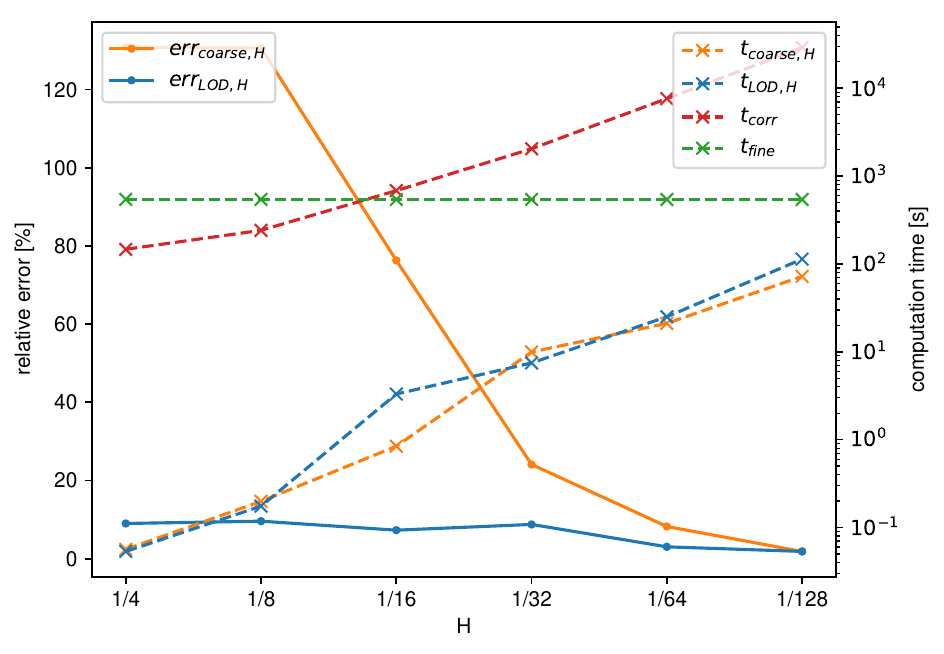}\qquad
  \includegraphics[height=0.325\linewidth]{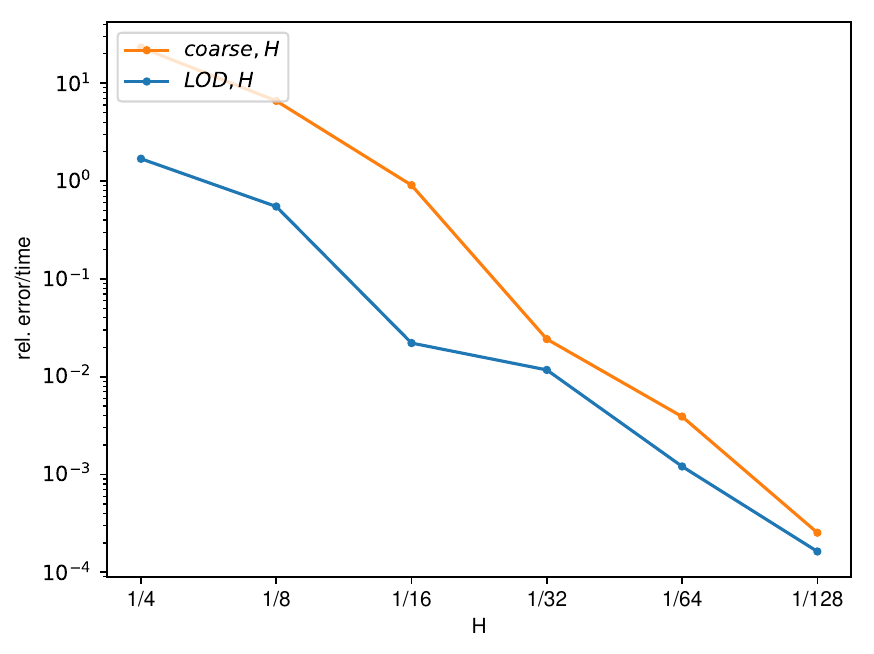}
  \caption{Model Problem \ref{sec:ev_problem} with $\nu=20$ and $\gamma=2\cdot10^4$. Timing and accuracy for different resolution $H$ of the coarse mesh (left) and performance as error/time (right).}
   \label{fig:performance_H}
\end{figure}
\begin{figure}
  \includegraphics[height=0.325\linewidth]{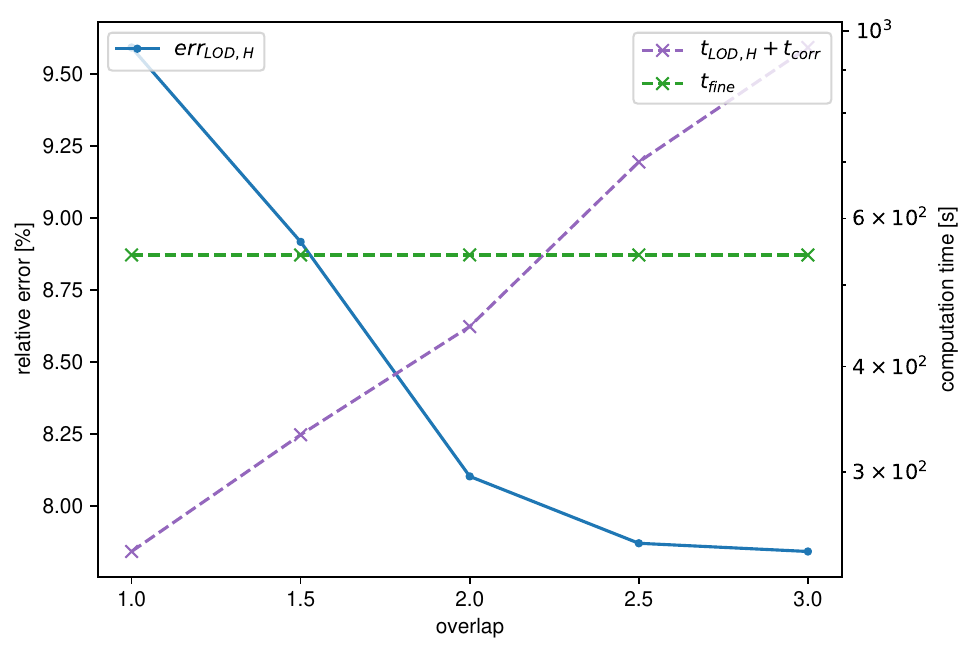}\qquad
  \includegraphics[height=0.325\linewidth]{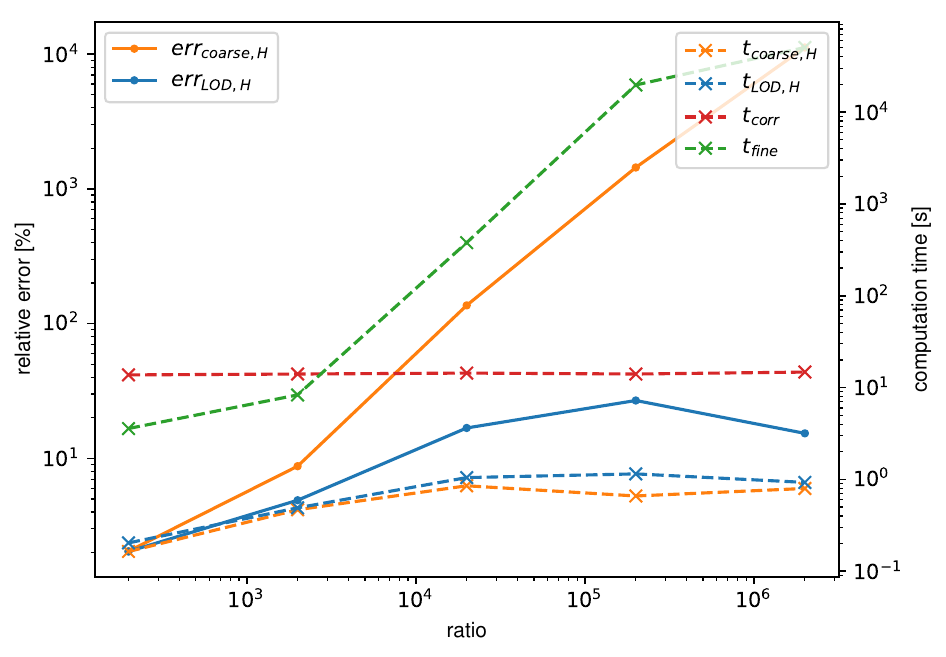}
  \caption{Model Problem \ref{sec:ev_problem}. Timing and accuracy of the LOD for (left) different overlap
    sizes ($H=2^{-3}$, $h=2^{-8}$, $\gamma=2\cdot 10^4$) and (right)
    different ratios $\gamma$ ($H=2^{-3}$, $h=2^{-8}$, overlap
    $H$).}
   \label{fig:performance_and_accuracy}
\end{figure}

\subsection{Eigenvalue Problem}
\label{sec:ev_problem}
In the previous model we demonstrated the basic performance of the algorithm for one fixed source term. However, the big advantage of the LOD is getting more pronounced, if computations have to be repeated with different source terms. To pick a natural setting where this is the case, our next model is an eigenvalue problem. Here we follow the LOD algorithm described in Section~\ref{section:LOD-for-EV}. We consider computing the 20 smallest eigenvalues
for the stationary linear Schr\"odinger equation with a discontinuous potential $V$. The problem is inspired by the Kronig-Penney model (cf. \cite{Varga:2011:CNA:2000306}) and exhibits fine scale heterogeneities and high contrast of order $\gamma$. Here we seek eigenfunctions $u$ for the 20 smallest eigenvalues $\lambda$ such that
\begin{subequations}
\begin{align}
  \label{eq:ev-problem}
  -\Delta u + V(x)u &= \lambda u\phantom{0} \quad\text{in }\Omega,  \\
  u &= 0\phantom{\lambda u} \quad\text{on }\partial\Omega,
  \intertext{with}
           V(x) &= \gamma \Big\lceil \cos( \pi \nu (x_1+0.1) )
                  \cos( \pi \nu x_2 ) \Big\rceil\,.
\end{align}
\end{subequations}
As the \texttt{python} implementation is not parallelized, we enforced
sequential computation of the eigenvalue problem. Different
resolutions for the coarse space and different overlap sizes are
considered. We solve for the 20 smallest eigenvalues and eigenmodes
and compare the obtained eigenvalues for the coarse discretization,
the full discretization on a fine mesh with resolution $h=2^{-8}$ and the LOD. Overall the obtained timings are
convincing and the results of the LOD yield good accuracy. All
eigenvalue problems are solved using the \texttt{ARPACK} library
\cite{arpack1998} with the shift-inverse method. For the inner solve (in the inverse power iteration)
we used a CG Krylov solver with the \texttt{pyamg}\cite{BeOlSc2008}
algebraic multigrid as a
preconditioner, which then leads to reproducible robust results.
In our numerical experiments we used the following parameters, if not
indicated differently for a particular experiment.
The domain $\Omega=(0,2)\times(0,3)$, the wave-number $\nu=20$, and the
ratio $\gamma=2\cdot10^4$. Figure~\ref{comparison} shows four selected
eigenmodes for a fine scale resolution of $h=2^{-8}$ and a coarse
resolution of $H=2^{-3}$. We compare the coarse simulation (which is too
coarse to actually pick up the fine scale structure of the solution),
the full simulation, using the fine scale discretization and two
LOD simulations, with an overlap of $H$ and $2H$.
We measured the computational time and the error of different
setup. Figure~\ref{fig:performance_H} shows the relative error and the costs (\text{relative error / time}) for different coarse mesh sizes and a fixed overlap of $H$. From the left graph in Figure~\ref{fig:performance_H} we can see that the LOD error $\operatorname{err}_{\LOD,H}$ remains small for all coarse mesh sizes, where the total LOD computing time $t_{\LOD,H}+t_{\corr}$ is smaller than the computing time $t_{\fine}$ for the full fine scale discretization, which is a remarkable observation. We can also see that the ``offline costs'' $t_{\corr}$ (i.e. the time for computing the corrections) still amounts for the major costs of the LOD, at least for reasonably coarse mesh sizes $H$.  

Even more interesting results are depicted in Figure~\ref{fig:performance_and_accuracy}, which shows the relative error for overlap sizes and different coefficient ratios $\gamma$. The error was measured with respect to the eigenvalues of the
fine scale solution. Given the vectors $\Lambda_{\coarse,H}$,
$\Lambda_{LOD,H}$ and $\Lambda_{\fine}$ of the first 20 eigenvalues, the
relative errors are computed as
\begin{align*}
  \operatorname{err}_{\coarse,H} &= \| \Lambda_{\coarse,H} -
                                \Lambda_{\fine}\|_{\infty}\,, &  \operatorname{err}_{\LOD,H} &= \| \Lambda_{\LOD,H} - \Lambda_{\fine}\|_{\infty}\,.
\end{align*}
For the computation time, we consider the time $t_{\fine}$ for computing
the 20 smallest eigenvectors and their eigenmodes using the fine scale
discretization, the time $t_{\coarse,H}$ for the eigenvalue solve of the
coarse system
and $t_{\LOD,H}$ for the eigenvalue solve of the LOD system. In addition
the LOD requires the usual preprocessing, in order to compute
the correction of coarse space basis. The time necessary for this
setup phase is denoted by $t_{\corr}$ as before.
The left graph of Figure~\ref{fig:performance_and_accuracy} shows that even for a high contrast parameter $\gamma$, an overlap of two coarse elements is sufficient so that truncation error becomes negligible. The right figure shows the efficiency of the LOD compared to a regular fine-scale computation (using AMG), in particular for large $\gamma$. For $\gamma \approx 5\cdot 10^4$ both approaches show an equal performance. For $\gamma \ge 5\cdot 10^4$ the efficiency of the LOD compared to the standard implementation is continuously increasing. In this regime, we observe that even for a sequential implementation of the LOD method, the total run time $t_{\LOD,H}+t_{\corr}$ is significantly below the run times $t_{\coarse}/t_{\fine}$ of the multigrid alternative. Here the LOD was up to 1000 times faster than the multigrid alternative, which shows its enormous potential for even larger problems.
   
\subsection{Efficiency}
These examples just consider
the actual solve and we didn't compute the additional post-processing
described in section \ref{sec:two-grid-postprocessing} and the timings exclude the
assembly of the different operators. All computation done in DUNE are
very fast and negligible compared to the overall solving time. In
practice, one would use several improvements, in particular, the
different sub-domain problems (algorithm \ref{alg:main}, line
\ref{alg:main:subdomains}) are completely independent and can be
solved in parallel, which yields a perfect speedup and is easy to
implement for modern many-core systems. Additionally the
necessary solves for different right-hand-sides in algorithm \ref{alg:main}, line
\ref{alg:main:vertices} allows for a slight reformulation of the
algorithm and using vector instructions (e.g. SSE, AVX, Neon) to
compute all updates in a single run.

\section{Conclusions}
In this contribution we presented an efficient implementation of the
Localized Orthogonal Decomposition (LOD), including several
applications and variations of the methodology. The efficiency of the
algorithms is verified in numerical experiments, where we demonstrated
that the approach can be even very powerful in its sequential
implementation. This aspect is specifically stressed by Figure
\ref{fig:performance_and_accuracy} where we compare the CPU times for
the LOD with the CPU times of an efficient algebraic multigrid solver
(AMG) on the fine scale depending on $\gamma$. It can be observed that
the computational complexity of the LOD is independent of the contrast
parameter $\gamma$, whereas the reference solver is not robust. For
$\gamma=10^6$, the run times using algebraic multigrid were of order
$10^3$ times higher than for the LOD. Furthermore, we observed that
the AMG implementation suffered from \quotes{false eigenvalues} caused
by numerical rounding errors. Ratios of order
$\gamma>\mathcal{O}(10^6)$ could no longer be handled by the AMG,
whereas the LOD was still performing well and with the same run times
as for small values of $\gamma$.

\section*{Acknowledgements} This work was supported by:
  the German Research Foundation (DFG) through the Priority Programme
  `Software for Exascale Computing (SPP 1648)' (grant EN-1042/2-1 and
  EN-1042/2-2) [supporting C. Engwer],
  the Swedish Research Council (grant 2016-03339) [supporting
  P. Henning] and (grant 2015-04964) [supporting A. M\r{a}lqvist],
  and the \emph{Hausdorff Institute for Mathematics} in Bonn
  for the kind hospitality of during the trimester program on multiscale
  problems
  [supporting P. Henning, A. M\r{a}lqvist and D. Peterseim].

\def\cprime{$'$}


\begin{thebibliography}{10}
\expandafter\ifx\csname url\endcsname\relax
  \def\url#1{\texttt{#1}}\fi
\expandafter\ifx\csname urlprefix\endcsname\relax\def\urlprefix{URL }\fi
\expandafter\ifx\csname href\endcsname\relax
  \def\href#1#2{#2} \def\path#1{#1}\fi

\bibitem{MaP14}
A.~M{\aa}lqvist, D.~Peterseim, Localization of elliptic multiscale problems,
  Math. Comp. 83~(290) (2014) 2583--2603.
\newblock \href {http://dx.doi.org/10.1090/S0025-5718-2014-02868-8}
  {\path{doi:10.1090/S0025-5718-2014-02868-8}}.

\bibitem{HeP13}
P.~Henning, D.~Peterseim, Oversampling for the {M}ultiscale {F}inite {E}lement
  {M}ethod, SIAM Multiscale Model. Simul. 11~(4) (2013) 1149--1175.
\newblock \href {http://dx.doi.org/10.1137/120900332}
  {\path{doi:10.1137/120900332}}.

\bibitem{HeM14}
P.~Henning, A.~M{\aa}lqvist, Localized orthogonal decomposition techniques for
  boundary value problems, SIAM J. Sci. Comput. 36~(4) (2014) A1609--A1634.
\newblock \href {http://dx.doi.org/10.1137/130933198}
  {\path{doi:10.1137/130933198}}.

\bibitem{EGM13}
D.~Elfverson, E.~H. Georgoulis, A.~M{\aa}lqvist, D.~Peterseim, Convergence of a
  discontinuous {G}alerkin multiscale method, SIAM J. Numer. Anal. 51~(6)
  (2013) 3351--3372.
\newblock \href {http://dx.doi.org/10.1137/120900113}
  {\path{doi:10.1137/120900113}}.

\bibitem{EGM13b}
D.~Elfverson, E.~H. Georgoulis, A.~M{\aa}lqvist, An adaptive discontinuous
  {G}alerkin multiscale method for elliptic problems, Multiscale Model. Simul.
  11~(3) (2013) 747--765.
\newblock \href {http://dx.doi.org/10.1137/120863162}
  {\path{doi:10.1137/120863162}}.

\bibitem{Elf14}
D.~Elfverson, A discontinuous {G}alerkin multiscale method for
  convection-diffusion problems, uppsala University Preprint (2014).

\bibitem{Mal11}
A.~M{\aa}lqvist, Multiscale methods for elliptic problems, Multiscale Model.
  Simul. 9~(3) (2011) 1064--1086.
\newblock \href {http://dx.doi.org/10.1137/090775592}
  {\path{doi:10.1137/090775592}}.

\bibitem{HHM15}
F.~Hellman, P.~Henning, A.~M{\aa}lqvist, Multiscale mixed finite elements,
  Discrete Contin. Dyn. Syst. Ser. S 9~(5) (2016) 1269--1298.
\newblock \href {http://dx.doi.org/10.3934/dcdss.2016051}
  {\path{doi:10.3934/dcdss.2016051}}.

\bibitem{HMP13b}
P.~Henning, P.~Morgenstern, D.~Peterseim, Multiscale partition of unity, in:
  M.~Griebel, M.~A. Schweitzer (Eds.), Meshfree Methods for Partial
  Differential Equations VII, Vol. 100 of Lecture Notes in Computational
  Science and Engineering, Springer International Publishing, 2015, pp.
  185--204.
\newblock \href {http://dx.doi.org/10.1007/978-3-319-06898-5_10}
  {\path{doi:10.1007/978-3-319-06898-5_10}}.

\bibitem{AbH15}
A.~Abdulle, P.~Henning, A reduced basis localized orthogonal decomposition, J.
  Comput. Phys. 295 (2015) 379--401.
\newblock \href {http://dx.doi.org/10.1016/j.jcp.2015.04.016}
  {\path{doi:10.1016/j.jcp.2015.04.016}}.

\bibitem{MaP15}
A.~M{\aa}lqvist, D.~Peterseim, Computation of eigenvalues by numerical
  upscaling, Numer. Math. 130~(2) (2015) 337--361.
\newblock \href {http://dx.doi.org/10.1007/s00211-014-0665-6}
  {\path{doi:10.1007/s00211-014-0665-6}}.

\bibitem{2015arXiv151005792M}
A.~M{\aa}lqvist, D.~Peterseim, Generalized finite element methods for quadratic
  eigenvalue problems, ESAIM Math. Model. Numer. Anal. 51~(1) (2017) 147--163.
\newblock \href {http://dx.doi.org/10.1051/m2an/2016019}
  {\path{doi:10.1051/m2an/2016019}}.

\bibitem{BrP14}
D.~L. Brown, D.~Peterseim, A multiscale method for porous microstructures,
  Multiscale Model. Simul. 14~(3) (2016) 1123--1152.
\newblock \href {http://dx.doi.org/10.1137/140995210}
  {\path{doi:10.1137/140995210}}.

\bibitem{2016arXiv160106549P}
D.~Peterseim, R.~Scheichl, Robust numerical upscaling of elliptic multiscale
  problems at high contrast, Comput. Methods Appl. Math. 16~(4) (2016)
  579--603.
\newblock \href {http://dx.doi.org/10.1515/cmam-2016-0022}
  {\path{doi:10.1515/cmam-2016-0022}}.

\bibitem{HeM17}
F.~Hellman, A.~M{\aa}lqvist, Contrast independent localization of multiscale
  problems, Multiscale Model. Simul. 15~(4) (2017) 1325--1355.
\newblock \href {http://dx.doi.org/10.1137/16M1100460}
  {\path{doi:10.1137/16M1100460}}.

\bibitem{2017arXiv170208858G}
D.~{Gallistl}, D.~{Peterseim}, {Numerical stochastic homogenization by
  quasilocal effective diffusion tensors}, ArXiv e-prints\href
  {http://arxiv.org/abs/1702.08858} {\path{arXiv:1702.08858}}.

\bibitem{2018arXiv180701741F}
M.~{Feischl}, D.~{Peterseim}, {Sparse Compression of Expected Solution
  Operators}, arXiv e-prints (2018) arXiv:1807.01741\href
  {http://arxiv.org/abs/1807.01741} {\path{arXiv:1807.01741}}.

\bibitem{HMP14}
P.~Henning, A.~M{\aa}lqvist, D.~Peterseim, A localized orthogonal decomposition
  method for semi-linear elliptic problems, M2AN Math. Model. Numer. Anal.
  48~(5) (2014) 1331--1349.
\newblock \href {http://dx.doi.org/10.1051/m2an/2013141}
  {\path{doi:10.1051/m2an/2013141}}.

\bibitem{AbH14c}
A.~Abdulle, P.~Henning, Localized orthogonal decomposition method for the wave
  equation with a continuum of scales, Math. Comp. 86~(304) (2017) 549--587.
\newblock \href {http://dx.doi.org/10.1090/mcom/3114}
  {\path{doi:10.1090/mcom/3114}}.

\bibitem{Peterseim.Schedensack:2016}
D.~Peterseim, M.~Schedensack, Relaxing the {CFL} condition for the wave
  equation on adaptive meshes, J. Sci. Comput. 72~(3) (2017) 1196--1213.

\bibitem{Maier2018}
R.~Maier, D.~Peterseim,
  \href{https://doi.org/10.1007/s10543-018-0735-8}{Explicit computational wave
  propagation in micro-heterogeneous media}, BIT Numerical Mathematics\href
  {http://dx.doi.org/10.1007/s10543-018-0735-8}
  {\path{doi:10.1007/s10543-018-0735-8}}.
\newline\urlprefix\url{https://doi.org/10.1007/s10543-018-0735-8}

\bibitem{MaPr15}
A.~M{\aa}lqvist, A.~Persson, Multiscale techniques for parabolic equations,
  Numer. Math. 138~(1) (2018) 191--217.
\newblock \href {http://dx.doi.org/10.1007/s00211-017-0905-7}
  {\path{doi:10.1007/s00211-017-0905-7}}.

\bibitem{refId0}
{M\aa{}lqvist, Axel}, {Persson, Anna},
  \href{https://doi.org/10.1051/m2an/2016054}{A generalized finite element
  method for linear thermoelasticity}, ESAIM: M2AN 51~(4) (2017) 1145--1171.
\newblock \href {http://dx.doi.org/10.1051/m2an/2016054}
  {\path{doi:10.1051/m2an/2016054}}.
\newline\urlprefix\url{https://doi.org/10.1051/m2an/2016054}

\bibitem{2018arXiv180100615A}
R.~{Altmann}, E.~{Chung}, R.~{Maier}, D.~{Peterseim}, S.-M. {Pun},
  {Computational multiscale methods for linear heterogeneous poroelasticity},
  arXiv e-prints (2018) arXiv:1801.00615\href {http://arxiv.org/abs/1801.00615}
  {\path{arXiv:1801.00615}}.

\bibitem{EGH15}
D.~Elfverson, V.~Ginting, P.~Henning,
  \href{http://dx.doi.org/10.1007/s00211-015-0703-z}{On multiscale methods in
  petrov--galerkin formulation}, Numerische Mathematik 131~(4) (2015) 643--682.
\newblock \href {http://dx.doi.org/10.1007/s00211-015-0703-z}
  {\path{doi:10.1007/s00211-015-0703-z}}.
\newline\urlprefix\url{http://dx.doi.org/10.1007/s00211-015-0703-z}

\bibitem{doi:10.1137/17M1147305}
D.~Brown, J.~Gedicke, D.~Peterseim,
  \href{https://doi.org/10.1137/17M1147305}{Numerical homogenization of
  heterogeneous fractional laplacians}, Multiscale Modeling \& Simulation
  16~(3) (2018) 1305--1332.
\newblock \href {http://arxiv.org/abs/https://doi.org/10.1137/17M1147305}
  {\path{arXiv:https://doi.org/10.1137/17M1147305}}, \href
  {http://dx.doi.org/10.1137/17M1147305} {\path{doi:10.1137/17M1147305}}.
\newline\urlprefix\url{https://doi.org/10.1137/17M1147305}

\bibitem{Pet14}
D.~Peterseim, Eliminating the pollution effect in {H}elmholtz problems by local
  subscale correction, Math. Comp. 86~(305) (2017) 1005--1036.
\newblock \href {http://dx.doi.org/10.1090/mcom/3156}
  {\path{doi:10.1090/mcom/3156}}.

\bibitem{GaP15}
D.~Gallistl, D.~Peterseim, Stable multiscale {P}etrov-{G}alerkin finite element
  method for high frequency acoustic scattering, Comp. Meth. Appl. Mech. Eng.
  295 (2015) 1--17.
\newblock \href {http://dx.doi.org/10.1016/j.cma.2015.06.017}
  {\path{doi:10.1016/j.cma.2015.06.017}}.

\bibitem{Brown.Gallistl.Peterseim:2015}
D.~Brown, D.~Gallistl, D.~Peterseim, Multiscale {P}etrov-{G}alerkin method for
  high-frequency heterogeneous {H}elmholtz equations, in: M.~Griebel, M.~A.
  Schweitzer (Eds.), Meshfree Methods for Partial Differential Equations VII,
  Lecture Notes in Computational Science and Engineering, Springer, 2017, pp.
  85--115.

\bibitem{2016arXiv160804243B}
D.~L. {Brown}, D.~{Gallistl}, {Multiscale Sub-grid Correction Method for
  Time-Harmonic High-Frequency Elastodynamics with Wavenumber Explicit Bounds},
  arXiv e-prints (2016) arXiv:1608.04243\href {http://arxiv.org/abs/1608.04243}
  {\path{arXiv:1608.04243}}.

\bibitem{HMP14b}
P.~Henning, A.~M{\aa}lqvist, D.~Peterseim, Two-{L}evel {D}iscretization
  {T}echniques for {G}round {S}tate {C}omputations of {B}ose-{E}instein
  {C}ondensates, SIAM J. Numer. Anal. 52~(4) (2014) 1525--1550.
\newblock \href {http://dx.doi.org/10.1137/130921520}
  {\path{doi:10.1137/130921520}}.

\bibitem{Peterseim:2015}
D.~Peterseim, Variational multiscale stabilization and the exponential decay of
  fine-scale correctors, in: Building bridges: connections and challenges in
  modern approaches to numerical partial differential equations, Vol. 114 of
  Lect. Notes Comput. Sci. Eng., Springer, [Cham], 2016, pp. 341--367.

\bibitem{Hug95}
T.~J.~R. Hughes, Multiscale phenomena: {G}reen's functions, the
  {D}irichlet-to-{N}eumann formulation, subgrid scale models, bubbles and the
  origins of stabilized methods, Comput. Methods Appl. Mech. Engrg. 127~(1-4)
  (1995) 387--401.
\newblock \href {http://dx.doi.org/10.1016/0045-7825(95)00844-9}
  {\path{doi:10.1016/0045-7825(95)00844-9}}.

\bibitem{HFM98}
T.~J.~R. Hughes, G.~R. Feij{\'o}o, L.~Mazzei, J.-B. Quincy, The variational
  multiscale method---a paradigm for computational mechanics, Comput. Methods
  Appl. Mech. Engrg. 166~(1-2) (1998) 3--24.
\newblock \href {http://dx.doi.org/10.1016/S0045-7825(98)00079-6}
  {\path{doi:10.1016/S0045-7825(98)00079-6}}.

\bibitem{HuS07}
T.~J.~R. Hughes, G.~Sangalli, Variational multiscale analysis: the fine-scale
  {G}reen's function, projection, optimization, localization, and stabilized
  methods, SIAM J. Numer. Anal. 45~(2) (2007) 539--557.
\newblock \href {http://dx.doi.org/10.1137/050645646}
  {\path{doi:10.1137/050645646}}.

\bibitem{LaM07}
M.~G. Larson, A.~M{\aa}lqvist, Adaptive variational multiscale methods based on
  a posteriori error estimation: energy norm estimates for elliptic problems,
  Comput. Methods Appl. Mech. Engrg. 196~(21-24) (2007) 2313--2324.
\newblock \href {http://dx.doi.org/10.1016/j.cma.2006.08.019}
  {\path{doi:10.1016/j.cma.2006.08.019}}.

\bibitem{LaM09}
M.~G. Larson, A.~M{\aa}lqvist, An adaptive variational multiscale method for
  convection-diffusion problems, Comm. Numer. Methods Engrg. 25~(1) (2009)
  65--79.
\newblock \href {http://dx.doi.org/10.1002/cnm.1106}
  {\path{doi:10.1002/cnm.1106}}.

\bibitem{2016arXiv160802092G}
D.~Gallistl, D.~Peterseim, Computation of quasi-local effective diffusion
  tensors and connections to the mathematical theory of homogenization,
  Multiscale Model. Simul. 15~(4) (2017) 1530--1552.
\newblock \href {http://dx.doi.org/10.1137/16M1088533}
  {\path{doi:10.1137/16M1088533}}.

\bibitem{2016arXiv160804081K}
R.~Kornhuber, D.~Peterseim, H.~Yserentant, An analysis of a class of
  variational multiscale methods based on subspace decomposition, Math. Comp.
  published electronically.
\newblock \href {http://dx.doi.org/10.1090/mcom/3302}
  {\path{doi:10.1090/mcom/3302}}.

\bibitem{2018arXiv181106319P}
D.~{Peterseim}, D.~{Varga}, B.~{Verf{\"u}rth}, {From Domain Decomposition to
  Homogenization Theory}, arXiv e-prints (2018) arXiv:1811.06319\href
  {http://arxiv.org/abs/1811.06319} {\path{arXiv:1811.06319}}.

\bibitem{Owhadi2014}
H.~Owhadi, L.~Zhang, L.~Berlyand, Polyharmonic homogenization, rough
  polyharmonic splines and sparse super-localization, ESAIM Math. Model. Numer.
  Anal. 48~(2) (2014) 517--552.
\newblock \href {http://dx.doi.org/10.1051/m2an/2013118}
  {\path{doi:10.1051/m2an/2013118}}.

\bibitem{KornhuberYserentant2015}
R.~Kornhuber, H.~Yserentant, Numerical homogenization of elliptic multiscale
  problems by subspace decomposition, Multiscale Model. Simul. 14~(3) (2016)
  1017--1036.
\newblock \href {http://dx.doi.org/10.1137/15M1028510}
  {\path{doi:10.1137/15M1028510}}.

\bibitem{Owhadi2017}
H.~Owhadi, Multigrid with rough coefficients and multiresolution operator
  decomposition from hierarchical information games, SIAM Rev. 59~(1) (2017)
  99--149.
\newblock \href {http://dx.doi.org/10.1137/15M1013894}
  {\path{doi:10.1137/15M1013894}}.

\bibitem{Car99}
C.~Carstensen, Quasi-interpolation and a posteriori error analysis in finite
  element methods, M2AN Math. Model. Numer. Anal. 33~(6) (1999) 1187--1202.
\newblock \href {http://dx.doi.org/10.1051/m2an:1999140}
  {\path{doi:10.1051/m2an:1999140}}.

\bibitem{MR0400739}
P.~Cl{\'e}ment, Approximation by finite element functions using local
  regularization, Rev. Fran\c caise Automat. Informat. Recherche
  Op\'erationnelle S\'er. RAIRO Analyse Num\'erique 9~(R-2) (1975) 77--84.

\bibitem{MR1011446}
L.~R. Scott, S.~Zhang, Finite element interpolation of nonsmooth functions
  satisfying boundary conditions, Math. Comp. 54~(190) (1990) 483--493.
\newblock \href {http://dx.doi.org/10.2307/2008497}
  {\path{doi:10.2307/2008497}}.

\bibitem{MR1706735}
C.~Carstensen, R.~Verf{\"u}rth, Edge residuals dominate a posteriori error
  estimates for low order finite element methods, SIAM J. Numer. Anal. 36~(5)
  (1999) 1571--1587.
\newblock \href {http://dx.doi.org/10.1137/S003614299732334X}
  {\path{doi:10.1137/S003614299732334X}}.

\bibitem{2015arXiv150506931E}
A.~Ern, J.-L. Guermond, Finite element quasi-interpolation and best
  approximation, ESAIM Math. Model. Numer. Anal. 51~(4) (2017) 1367--1385.

\bibitem{10.1007/978-3-319-52389-7_21}
R.~Kornhuber, J.~Podlesny, H.~Yserentant, Direct and iterative methods for
  numerical homogenization, in: C.-O. Lee, X.-C. Cai, D.~E. Keyes, H.~H. Kim,
  A.~Klawonn, E.-J. Park, O.~B. Widlund (Eds.), Domain Decomposition Methods in
  Science and Engineering XXIII, Springer International Publishing, Cham, 2017,
  pp. 217--225.

\bibitem{XuZh01}
J.~Xu, A.~Zhou, A two-grid discretization scheme for eigenvalue problems, Math.
  Comp. 70~(233) (2001) 17--25.
\newblock \href {http://dx.doi.org/10.1090/S0025-5718-99-01180-1}
  {\path{doi:10.1090/S0025-5718-99-01180-1}}.

\bibitem{walt2011numpy}
S.~v.~d. Walt, S.~C. Colbert, G.~Varoquaux, The numpy array: a structure for
  efficient numerical computation, Computing in Science \& Engineering 13~(2)
  (2011) 22--30.

\bibitem{millman2011scipy}
K.~J. Millman, M.~Aivazis, Python for scientists and engineers, Computing in
  Science \& Engineering 13~(2) (2011) 9--12.

\bibitem{dune08:1}
P.~Bastian, M.~Blatt, A.~Dedner, C.~Engwer, R.~K. \noopsort{1}, M.~Ohlberger,
  O.~Sander, A generic grid interface for parallel and adaptive scientific
  computing. {P}art {I}: {A}bstract framework, Computing 82~(2--3) (2008)
  103--119.
\newblock \href {http://dx.doi.org/10.1007/s00607-008-0003-x}
  {\path{doi:10.1007/s00607-008-0003-x}}.

\bibitem{dune08:2}
P.~Bastian, M.~Blatt, A.~Dedner, C.~Engwer, R.~K. \noopsort{2}, R.~Kornhuber,
  M.~Ohlberger, O.~Sander, A generic grid interface for parallel and adaptive
  scientific computing. {P}art {II}: {I}mplementation and tests in {DUNE},
  Computing 82~(2--3) (2008) 121--138.
\newblock \href {http://dx.doi.org/10.1007/s00607-008-0004-9}
  {\path{doi:10.1007/s00607-008-0004-9}}.

\bibitem{blatt2016distributed}
M.~Blatt, A.~Burchardt, A.~Dedner, C.~Engwer, J.~Fahlke, B.~Flemisch,
  C.~Gersbacher, C.~Gr{\"a}ser, F.~Gruber, C.~Gr{\"u}ninger, et~al., The
  distributed and unified numerics environment, version 2.4, Archive of
  Numerical Software 4~(100) (2016) 13--29.

\bibitem{bastian2010generic}
P.~Bastian, F.~Heimann, S.~Marnach, Generic implementation of finite element
  methods in the distributed and unified numerics environment (dune),
  Kybernetika 46~(2) (2010) 294--315.

\bibitem{Varga:2011:CNA:2000306}
K.~Varga, J.~A. Driscoll, Computational Nanoscience: Applications for
  Molecules, Clusters, and Solids, 1st Edition, Cambridge University Press, New
  York, NY, USA, 2011.

\bibitem{arpack1998}
R.~B. Lehoucq, D.~C. Sorensen, C.~Yang, ARPACK users' guide: solution of
  large-scale eigenvalue problems with implicitly restarted Arnoldi methods,
  SIAM, 1998.

\bibitem{BeOlSc2008}
W.~N. Bell, L.~N. Olson, J.~Schroder, \href{http://www.pyamg.org}{{PyAMG:
  Algebraic Multigrid Solvers in Python}}, version 2.1 (2013).
\newline\urlprefix\url{http://www.pyamg.org}

\end{thebibliography}
\end{document}